%% file: secDeriv_FSBP.tex
\documentclass[final,onefignum,onetabnum,letterpaper]{siamonline190516}
\usepackage[scale=0.8]{geometry} % Reduce document margins

\input{secDeriv_FSBP_shared}

\usepackage{lineno}
\ifpdf
\hypersetup{
  pdftitle={Second-derivative FSBP operators},
  pdfauthor={Jan Glaubitz}
}
\fi

\begin{document}

\maketitle

% REQUIRE D
\begin{abstract}
	\input{0_abstract}
\end{abstract}

% REQUIRED
\begin{keywords}
	Summation-by-parts operators, second derivatives, advection-diffusion problems, wave equation, mimetic discretizations, general function spaces 
\end{keywords}

% REQUIRED
\begin{AMS}
	65M12, % Stability and convergence of numerical methods for initial value and initial-boundary value problems involving PDEs
	65M60, % Finite element, Rayleigh-Ritz and Galerkin methods for initial value and initial-boundary value problems involving PDEs
	65M70, % Spectral, collocation and related methods for initial value and initial-boundary value problems involving PDEs
	65D25 % Numerical differentiation
\end{AMS}

% Repository where code can be found (if there is any) 
\begin{Code}
    \url{https://github.com/jglaubitz/2ndDerivativeFSBP}
\end{Code}

% Once the paper is published
\begin{DOI}
	\url{https://doi.org/10.1016/j.jcp.2024.112889}
\end{DOI}

\input{1_introduction} 
\input{2_prelim} 
\input{3_existence}
\input{4_examples}
\input{5_numerics}
\input{6_summary}

%\appendix 
%\input{app_expl_A}

\section*{Acknowledgements} 
JG was supported by the US DOD (ONR MURI) grant \#N00014-20-1-2595. 
JN was supported by the Vetenskapsr\r{a}det Sweden grant 2021-05484 VR. 
P\"O was supported by the Gutenberg Research College, JGU Mainz and by the DFG within SPP 2410, project 525866748 (OE 661/5-1) and under the personal grant 520756621 (OE 661/4-1).

\bibliographystyle{siamplain}
\bibliography{references}

\end{document}

%% file: secDeriv_FSBP_shared.tex
\usepackage[normalem]{ulem} % to strikethrough text
\usepackage{lipsum}
\usepackage{amsfonts}
\usepackage{epstopdf}
\ifpdf
  \DeclareGraphicsExtensions{.eps,.pdf,.png,.jpg}
\else
  \DeclareGraphicsExtensions{.eps}
\fi

\usepackage{datetime}
\newdateformat{monthyeardate}{%
  \monthname[\THEMONTH] \THEDAY, \THEYEAR}

\usepackage{academicons}
\usepackage{xcolor}
\renewcommand{\orcid}[1]{\href{https://orcid.org/#1}{\textcolor[HTML]{A6CE39}{orcid.org/#1}}}

\usepackage{amsmath}
\allowdisplaybreaks
\usepackage{amssymb}
\usepackage{commath}
\usepackage{mathtools}
\usepackage{bbm}

\usepackage{color}
\usepackage{graphicx}
\usepackage[small]{caption}
\usepackage{subcaption}

\usepackage{relsize}
\usepackage{adjustbox}
\usepackage{algorithm}
\usepackage[noend]{algpseudocode}
\usepackage{booktabs}
\usepackage{tikz}

\usepackage{verbatim}

\usepackage{enumitem}
\setlist[enumerate]{leftmargin=.5in}
\setlist[itemize]{leftmargin=.5in}

\newsiamthm{problem}{Problem}
\newsiamremark{remark}{Remark}
\newsiamremark{hypothesis}{Hypothesis} 
\crefname{hypothesis}{Hypothesis}{Hypotheses}
\newsiamremark{example}{Example}
\newsiamthm{claim}{Claim}
\newsiamthm{conjecture}{Conjecture}

\headers{Second-derivative FSBP operators}{J.\ Glaubitz, S.-C.\ Klein, J.\ Nordstr\"om, and P.\ \"Offner}

\title{Summation-by-parts operators for general function spaces: The second derivative
\thanks{
\monthyeardate\today 
\corresponding{Jan Glaubitz} 
}}

\author{
Jan Glaubitz\thanks{Department of Aeronautics and Astronautics \& Laboratory for Information and Decision Systems, Massachusetts Institute of Technology, Cambridge, MA, USA (\email{glaubitz@mit.edu}, \orcid{0000-0002-3434-5563})}
\and 
Simon-Christian Klein\thanks{Department of Mathematics, TU Braunschweig, Braunschweig, Germany (\email{simon-christian.klein@tu-bs.de}, \orcid{0000-0002-8710-9089})}
\and 
{Jan Nordstr\"om\thanks{Department of Mathematics, Link\"oping University, Link\"oping, Sweden (\email{jan.nordstrom@liu.se}, \orcid{0000-0002-7972-6183})}} \thanks{Department of Mathematics and Applied Mathematics, University of Johannesburg, Johannesburg, South Africa}
\and 
Philipp \"Offner\thanks{Institute of Mathematics, Johannes Gutenberg University, Mainz, Germany (\email{poeffner@uni-mainz.de}, \orcid{0000-0002-1367-1917})} 
}

\usepackage{amsopn}

\definecolor{mLightGreen}{HTML}{14B03D}

\DeclareMathOperator{\diag}{diag}
 
\newcommand{\scp}[2]{\left\langle{#1, #2}\right\rangle} 
\newcommand{\Span}{\mathrm{span}}

\renewcommand{\d}{\mathrm{d}}
\newcommand{\intd}{\, \mathrm{d}}

\newcommand{\R}{\mathbb{R}}

%% file: 0_abstract.tex
% 1. Why should people care? 
Many applications rely on solving time-dependent partial differential equations (PDEs) that include second derivatives.
Summation-by-parts (SBP) operators are crucial for developing stable, high-order accurate numerical methodologies for such problems.   
% 2. What is the specific research question? 
Conventionally, SBP operators are tailored to the assumption that polynomials accurately approximate the solution, and SBP operators should thus be exact for them. 
However, this assumption falls short for a range of problems for which other approximation spaces are better suited.
We recently addressed this issue and developed a theory for \emph{first-derivative} SBP operators based on general function spaces, coined function-space SBP (FSBP) operators. 
% 3. What is done here? 
In this paper, we extend the innovation of FSBP operators to accommodate second derivatives. 
% 4. What are the key findings? 
The developed second-derivative FSBP operators maintain the desired mimetic properties of existing polynomial SBP operators while allowing for greater flexibility by being applicable to a broader range of function spaces. 
We establish the existence of these operators and detail a straightforward methodology for constructing them. 
By exploring various function spaces, including trigonometric, exponential, and radial basis functions, we illustrate the versatility of our approach. 
% 5. What are the implications of these findings? 
The work presented here opens up possibilities for using second-derivative SBP operators based on suitable function spaces, paving the way for a wide range of applications in the future. 

%% file: 1_introduction.tex
\section{Introduction} 
\label{sec:introduction} 

% 1. Why should others care? 
% 2. What has been done already? 
By mimicking integration-by-parts on a discrete level, summation-by-parts (SBP) operators can be used in combination with weakly enforced boundary conditions (BCs) to systematically develop energy-stable numerical methods for energy-bounded initial boundary value problems (IBVPs) \cite{svard2014review,fernandez2014review,chen2020review}. 
They are often used to discretize first derivatives, e.g., in first-order hyperbolic conservation laws and ordinary differential equations (ODEs), resulting in high-order accurate and provable stable schemes. 
Examples of schemes based on SBP operators include finite difference (FD) \cite{kreiss1974finite,kreiss1977existence,scherer1977energy,strand1994summation}, essentially non-oscillatory (ENO) and weighted ENO (WENO) \cite{yamaleev2009systematic,fisher2011boundary,carpenter2016entropy}, continuous Galerkin (CG) \cite{abgrall2020analysis,abgrall2021analysis}, discontinuous Galerkin (DG) \cite{gassner2013skew,chen2017entropy}, finite volume (FV) \cite{nordstrom2001finite,nordstrom2003finite}, flux reconstruction (FR) \cite{huynh2007flux,ranocha2016summation,offner2018stability}, and implicit time integration \cite{nordstrom2013summation,linders2020properties,ranocha2021new} methods. 
Second-derivative SBP operators have been the subject of significant research attention as well \cite{mattsson2004summation,mattsson2012summation,del2015generalized}, for instance in the context of the Navier--Stokes equations \cite{mattsson2008stable,svard2008stable}, wave equations \cite{nordstrom2003high,mattsson2006high,mattsson2008stablewave,mattsson2009stable,wang2016high,ranocha2021broad,wang2022energy}, and second-order ODEs \cite{nordstrom2016summation}. 
Moreover, second-derivative SBP operators can be used as building blocks in artificial dissipation operators and filtering procedures to stabilize numerical methods for hyperbolic problems \cite{mattsson2004stable,ranocha2018stability,glaubitz2020shock,lundquist2020stable,nordstrom2021stable}. 

% 3. What is the problem with what has already been done? 
% Motivate general function spaces
The current SBP operators are designed to be exact for polynomials of a specific degree, assuming that polynomials are the best approximation space for the solution. 
However, in some cases, other function spaces provide more accurate approximations. 
Previous studies have highlighted the benefits of using non-polynomial approximation spaces. 
These include exponentially fitted schemes to solve singular perturbation problems \cite{kadalbajoo2003exponentially,kalashnikova2009discontinuous}, DG methods \cite{yuan2006discontinuous} and (W)ENO reconstructions \cite{christofi1996study,iske1996structure,hesthaven2019entropy} based on non-polynomial approximation spaces, radial basis function (RBF) schemes \cite{fornberg2015solving,fornberg2015primer}, and methods based on rational function approximations \cite{nakatsukasa2018aaa,gopal2019solving}. 

\begin{figure}[tb]
	\centering 
	\begin{subfigure}[b]{0.49\textwidth}
		\includegraphics[width=\textwidth]{%
      		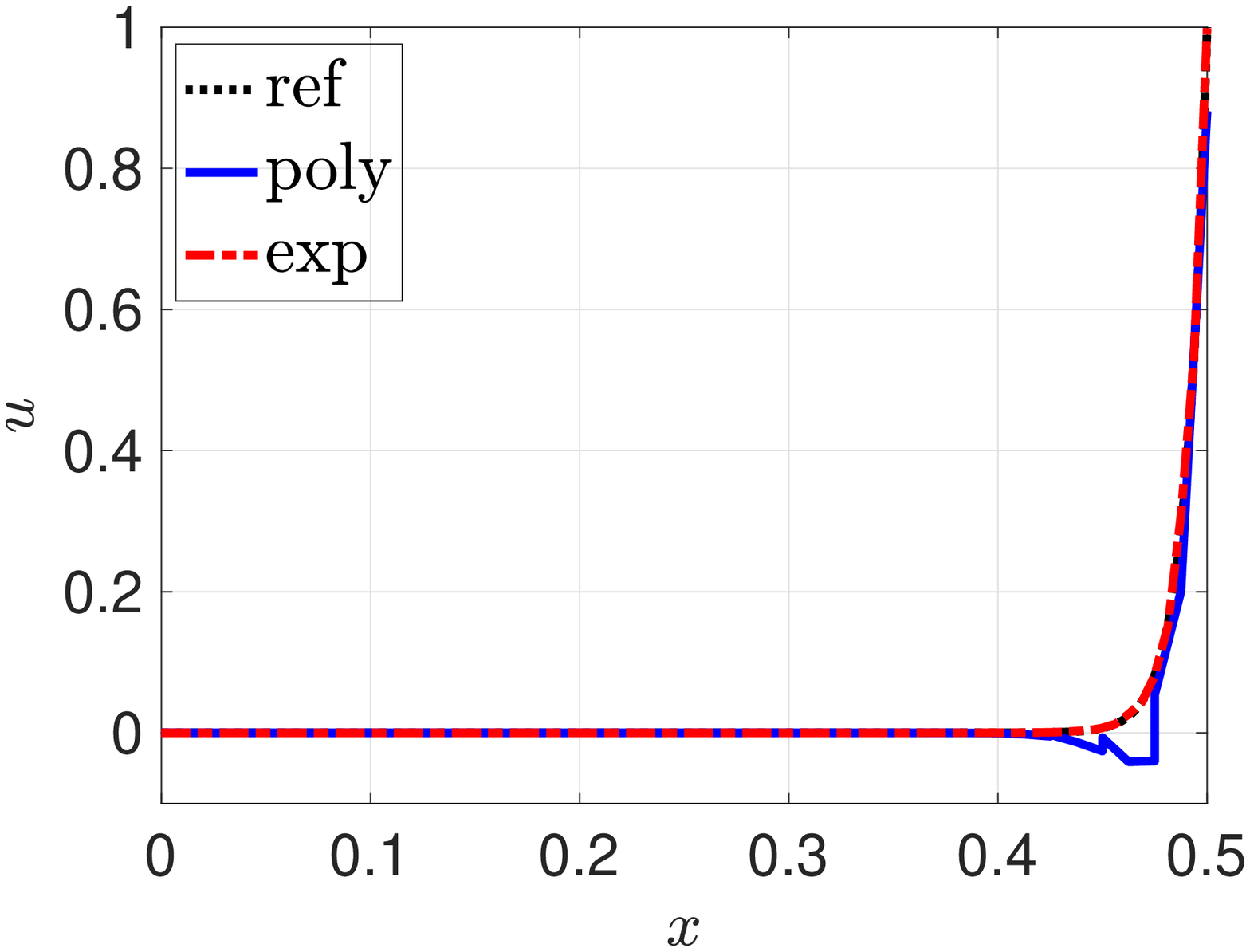} 
    		\caption{Boundary layer problem}
    		\label{fig:into_boundaryLayer}
  	\end{subfigure}%
  	\begin{subfigure}[b]{0.49\textwidth}
		\includegraphics[width=\textwidth]{%
      		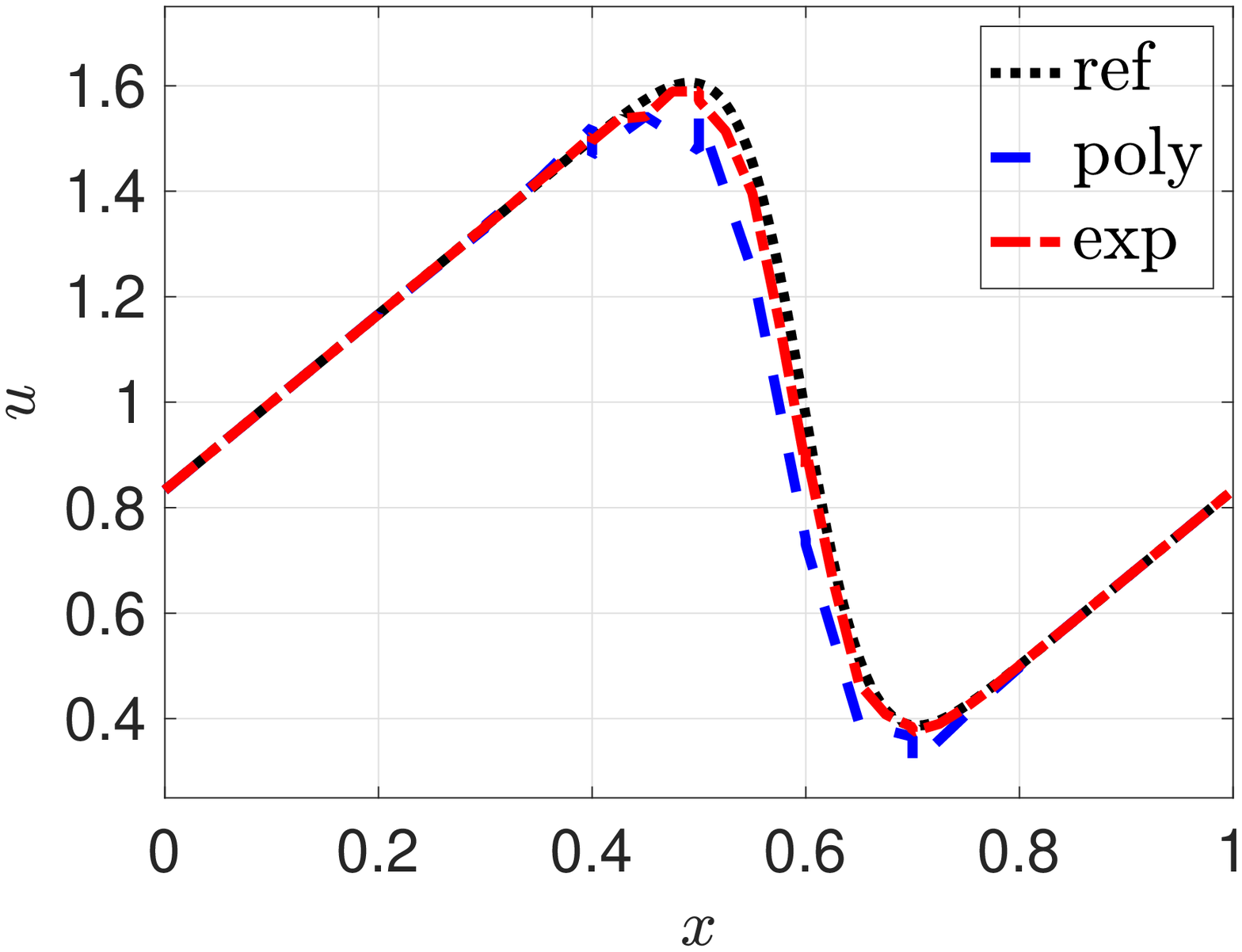} 
    		\caption{The viscous Burgers equation}
    		\label{fig:intro_viscBurgers}
  	\end{subfigure}%
  	\caption{
	Preview of the possible advantage of using non-polynomial approximation spaces for specific problems. 
	Illustrated here are the (numerical) solutions of a boundary layer problem and the viscous Burgers equation, which involve second derivatives. 
	The solutions to both problems feature steep gradients that can be better approximated using exponential rather than traditional polynomial approximation spaces. 
	We obtained the above numerical solutions using the FSBP-SAT scheme with polynomial (``poly") and exponential approximation spaces, $\mathcal{P}_2 =\Span\{1,x,x^2\}$ and $\mathcal{E}_2 = \Span\{1,x,e^x\}$, respectively. 
	See \cref{sub:bLayer,sub:Burgers} for details.  
	}
  	\label{fig:intro}
\end{figure}

% First-derivative FSBP operators
In \cite{glaubitz2022summation,glaubitz2023multi}, we recently developed the theory of function-space SBP (FSBP) operators, which are SBP operators approximating first derivatives designed for general function spaces---beyond just polynomials. 
The main advantage of using non-polynomial approximation spaces is that we can leverage prior knowledge about the behavior of the unknown solution. 
We demonstrated that the mimetic properties and construction strategies of polynomial-based SBP operators mostly carry over to FSBP operators. 
Another advantage is that introducing the FSBP framework into an existing method can lead to provable stability, as shown for RBF methods in \cite{glaubitz2022energy}, when combined with weakly enforced boundary conditions \cite{glaubitz2021stabilizing,glaubitz2021towards}. 
\cref{fig:intro} demonstrates the advantage of FSBP operators over conventional polynomial SBP operators for two specific cases. 
The first one (\cref{fig:into_boundaryLayer}) is a boundary layer solution of the advection-diffusion equation $\partial_t u + \partial_x u = \varepsilon \partial_xx u$, which is described in detail in \cref{sub:bLayer}. 
The second one is the viscous Burgers equation $\partial_t u + \partial_x \left( \frac{u^2}{2} \right) = \varepsilon \partial_{xx} u$, which is described in detail in \cref{sub:Burgers}. 
The solutions to both problems feature steep gradients that can be better approximated using, for instance, exponential rather than traditional polynomial approximation spaces. 
Importantly, both problems include second-derivative terms that should be discretized in an accurate and stable manner. 

%Our Contribution
% 4. What is done here?
Here, we systematically extend the framework of FSBP operators to include second derivatives. 
While equations with second derivatives can often be rewritten as systems of first derivatives, this results in larger systems that can require more computational effort. 
To bypass this, we will approximate second derivatives directly using general function spaces. 
% 5. What are the findings? 
We demonstrate that this yields desirable mimetic properties and stable discretizations. 
A possible approach for constructing polynomial second-derivative SBP operators is to apply the first-derivative operator twice. 
Although only rarely made explicit, this approach is possible because polynomial function spaces are invariant under differentiation (the derivative of a polynomial is again a polynomial). 
However, general function spaces are not necessarily invariant under differentiation.
Hence, special care must be taken to ensure accuracy for the proposed second-derivative FSBP operators. 
This motivates necessary changes in their construction compared to the polynomial case.  
% 6. What are the implications and applications? 
We exemplify the general theory by considering trigonometric, exponential, and radial basis functions. 
Our numerical experiments include the two-dimensional advection-diffusion equation, a boundary layer problem (see \cref{fig:into_boundaryLayer}), and the viscous Burgers equation (see \cref{fig:intro_viscBurgers}), which are solved using FSBP operators corresponding to various non-polynomial approximation spaces. 
% Restrictions 
We restrict the theoretical discussion to one-dimensional problems and constant viscosity distributions in this initial development. 
Future work will address the extensions to multi-dimensional problems, mixed second derivatives, and variable viscosity distributions.

% Outline
\subsection*{Outline} 

We explain the notation and provide preliminaries on first- and second-derivative FSBP operators in \Cref{sec:prelim}. 
In \Cref{sec:existence}, we prove the existence of second-derivative FSBP operators and present a simple procedure for their construction. 
\Cref{sec:examples} exemplifies this procedure for various function spaces, including polynomial, trigonometric, exponential, and radial basis function spaces.
Numerical experiments are showcased in \Cref{sec:numerics}, which include applications of the proposed second-derivative FSBP operators to the two-dimensional advection-diffusion equation, a boundary layer problem, and the viscous Burgers equation. 
We summarize the work in \Cref{sec:summary}.

%% file: 2_prelim.tex
\section{Background and motivation} 
\label{sec:prelim}

We explain the notation and provide some preliminaries on first-and second-derivative SBP operators.

\subsection{Notation}

We use the following notation in the remainder of this work. 
Let \(\Omega = [x_{\min}, x_{\max}]\) denote the physical domain for solving a specified time-dependent PDE. 
We focus on SBP operators applicable to a generic element \(\Omega_{\text{ref}} = [x_L, x_R]\). 
In approaches such as pseudo-spectral or global FD, \(\Omega_{\text{ref}}\) encompasses the entire physical domain, that is, \(\Omega_{\text{ref}} = \Omega\). 
Conversely, in FE and multi-block FD methods, where \(\Omega\) is subdivided into several smaller, non-overlapping elements \(\Omega_j\) satisfying \(\bigcup_{j} \Omega_j = \Omega\), any \(\Omega_j\) can serve as \(\Omega_{\text{ref}}\). 
Additionally, SBP operators are precomputed and stored, ready for subsequent use in the numerical treatment of a given PDE.
Consequently, the computational cost of precalculating these operators is relatively minor compared to the overall computational effort.

Let $\mathbf{x} = [x_1,\dots,x_N]$ be a vector of grid points on the reference element $\Omega_{\rm ref}$, where $x_L \leq x_1 < \dots < x_N \leq x_R$. 
If $f$ is a two-times continuously differentiable function on $\Omega_{\rm ref}$, then 
\begin{equation} 
\begin{aligned}
	\mathbf{f} & = [ f(x_1), \dots, f(x_N) ]^T, \\ 
	\mathbf{f'} & = [ f'(x_1), \dots, f'(x_N) ]^T, \\ 
	\mathbf{f''} & = [ f''(x_1), \dots, f''(x_N) ]^T 
\end{aligned}
\end{equation}
denote the nodal values of $f$ and first two derivatives at the grid points $\mathbf{x}$.

\subsection{The first derivative} 
\label{sub:first}

Consider the scalar linear advection equation 
\begin{equation}\label{eq:linear_adv}
	\partial_t u + \alpha \partial_x u = 0
\end{equation} 
on the reference element $\Omega_{\rm ref} = [x_L,x_R]$.
The exact solutions of \cref{eq:linear_adv} satisfy  
\begin{equation}\label{eq:linear_adv_energy} 
	\frac{\d}{\d t} \| u \|_{L^2}^2 
		= - \alpha \int_{x_L}^{x_R} 2 u ( \partial_x u ) \intd x 
		= -\alpha u^2 \big|_{x=x_L}^{x_R}, 
\end{equation}  
where $u^2 \big|_{x=x_L}^{x_R} = u(t,x_R)^2 - u(t,x_L)^2$.  
We semi-discretize \cref{eq:linear_adv} as 
\begin{equation}\label{eq:linear_adv_disc} 
	\frac{\d}{\d t} \mathbf{u} + \alpha D_1 \mathbf{u} = 0, 
\end{equation} 
where $\mathbf{u} = [u_1,\dots,u_N]^T$ denotes the vector of nodal values of the numerical solution at the grid points and $D_1$ is a discrete derivative operator approximating the continuous first-derivative $\partial_x$ on the grid points. 
Depending on how BCs are enforced, \cref{eq:linear_adv_disc} can further include simultaneous approximation terms (SATs)---or some other terms---to weakly enforce the BC \cite{svard2014review,fernandez2014review,chen2020review}. 
Now suppose it is reasonable to approximate the solution $u$ of \cref{eq:linear_adv} on $[x_L,x_R]$ using the function space $\mathcal{F} \subset C^1([x_L,x_R])$. 
In this case, the discrete first-derivative operator $D_1$ should be exact for all nodal vectors $\mathbf{u}$ that correspond to a function $u \in \mathcal{F}$. 
To this end, we recall the definition of FSBP operators, introduced in \cite{glaubitz2022summation}. 

\begin{definition}[First-derivative FSBP operators]\label{def:FSBP} 
	Let $\mathcal{F} \subset C^1([x_L,x_R])$. 
	An operator $D_1 = P^{-1} Q$, approximating $\partial_x$, is an \emph{$\mathcal{F}$-exact first-derivative FSBP operator} if 
	\begin{enumerate} 
		\item[(i)] 
		$D_1 \mathbf{f} = \mathbf{f'}$ for all $f \in \mathcal{F}$, 
		
		\item[(ii)] 
		$P$ is a symmetric positive definite matrix, and 
		
		\item[(iii)] 
		$Q + Q^T = B = \diag(-1,0,\dots,0,1)$,
		
	\end{enumerate}  
	where we have assumed $x_1 = x_L$ and $x_N = x_R$ in (iii).
\end{definition} 

Relation (i) ensures that $D_1$ accurately approximates $\partial_x$ by requiring $D_1$ to be exact for all functions from $\mathcal{F}$. 
Condition (ii) guarantees that $P$ induces a discrete inner product and norm, which are given by $\scp{\mathbf{u}}{\mathbf{v}}_P = \mathbf{u}^T P \mathbf{v}$ and $\|\mathbf{u}\|^2_P = \mathbf{u}^T P \mathbf{u}$ for $\mathbf{u},\mathbf{v} \in \R^N$, respectively. 
Relation (iii) encodes the SBP property, which allows us to mimic integration-by-parts (IBP) on a discrete level. 
In particular, (iii) ensures that the numerical solution given by \cref{eq:linear_adv_disc} mimics \cref{eq:linear_adv_energy}, opening up for energy-stability. 
See \cite{glaubitz2022summation,glaubitz2023multi} for more details on FSBP operators, their mimetic properties, and possible generalizations.

\subsection{The second derivative} 
\label{sub:second}

Now consider the heat equation  
\begin{equation}\label{eq:heat}
	\partial_t u = \beta \partial_{xx} u 
\end{equation} 
on the reference element $\Omega_{\rm ref} = [x_L,x_R]$.
While it is possible to rewrite \cref{eq:heat} as a system of first derivatives, solving this system can require more computational effort. 
To avoid this issue, we approximate the second derivative in \cref{eq:heat} directly. 
Suppose it is reasonable to approximate the solution of \cref{eq:heat} on $\Omega_{\rm ref} = [x_L,x_R]$ using the function space $\mathcal{F} \subset C^2([x_L,x_R])$. 
Then, a discrete second-derivative operator $D_2$ to approximate $\partial_{xx}$ should be exact for $\mathcal{F}$, i.e., $D_2 \mathbf{f} = \mathbf{f''}$ for all $f \in \mathcal{F}$. 
Moreover, IBP implies that exact solutions of \cref{eq:heat} satisfy 
\begin{equation}\label{eq:heat_energy} 
	\frac{\d}{\d t} \| u \|_{L^2}^2 
		= 2 \beta \int_{x_L}^{x_R} u ( \partial_{xx} u ) \intd x 
		= 2 \beta u (\partial_x u) \big|_{x=x_L}^{x_R} - 2 \beta \| \partial_x u \|_{L^2}^2,  
\end{equation} 
and $D_2$ should yield a numerical solution that mimics \cref{eq:heat_energy}. 
We can fully mimic \cref{eq:heat_energy} by choosing $D_2 = P^{-1}( BS - D_1^T B D_1 )$, where $D_1$ is a consistent approximation of $\partial_x$, $S$ includes a consistent approximation of $\partial_x$ at the boundary in the sense that $B S = B D_1$, and $B = \diag(-1,0,\dots,0,1)$.
The energy method applied to the semi-discrete equation $\frac{\d}{\d t} \mathbf{u} = \beta D_2 \mathbf{u}$ then leads to 
\begin{equation}\label{eq:heat_energy_discr_aux1} 
	\frac{\d}{\d t} \| \mathbf{u} \|_P^2 
		= 2 \beta \mathbf{u}^T B S \mathbf{u} - 2 \beta (D_1 \mathbf{u})^T P (D_1 \mathbf{u}).
\end{equation} 
Equation \cref{eq:heat_energy_discr_aux1} exactly mimics \cref{eq:heat_energy}, which motivates the following definition of second-derivative SBP operators for general function spaces. 

\begin{definition}[Second-derivative FSBP operators]\label{def:2nd_FSBP}
	Let $\mathcal{F} \subset C^2([x_L,x_R])$ and $D_2 = P^{-1}( BS - D_1^T P D_1 )$ be an approximation of $\partial_{xx}$, where $D_1$ is a consistent approximation of $\partial_x$, $S$ satisfies $B S = B D_1$, and $B = \diag(-1,0,\dots,0,1)$. 
	The operator $D_2$ is an \emph{$\mathcal{F}$-exact second-derivative FSBP operator} if 
	\begin{enumerate} 
		\item[(i)] 
		$D_2 \mathbf{f} = \mathbf{f''}$ for all $f \in \mathcal{F}$ and  
		
		\item[(ii)] 
		$P$ is a symmetric positive definite matrix. 
		
	\end{enumerate} 
\end{definition}
  
Similar to before, (i) ensures that $D_2$ accurately approximates $\partial_{xx}$ by requiring $D_2$ to be exact for all functions from $\mathcal{F}$ and (ii) guarantees that $P$ induces a discrete inner product and norm. 
A few remarks are in order. 

\begin{remark} 
	Strictly speaking, it is not necessary to exactly mimic IBP---as in \cref{eq:heat_energy_discr_aux1}---to get an energy estimate. 
	Consider $D_2 = P^{-1}( BS - A )$ approximating $\partial_{xx}$. 
	Then, the energy method applied to the semi-discrete equation $\frac{\d}{\d t} \mathbf{u} = \beta D_2 \mathbf{u}$ yields 
	\begin{equation}\label{eq:heat_energy_discr_aux2} 
		\frac{\d}{\d t} \| \mathbf{u} \|_P^2 
			= 2 \beta \mathbf{u}^T B D_1 \mathbf{u} - \beta \mathbf{u}^T (A + A^T) \mathbf{u}.
	\end{equation} 
	To get an energy estimate, it suffices that the symmetric part of $A$, $A + A^T$, is positive semi-definite, assuming that the BCs are enforced weakly using appropriate SATs.\footnote{ 
	An energy estimate for the wave equation $\partial_{tt} u = \partial_{xx} u$ can be obtained when $A$ additionally is symmetric \cite[Appendix A]{mattsson2004summation}.
}
	In the context of FD-SBP operators, this generalization is necessary to construct compact second-derivative FD-SBP operators \cite{mattsson2004summation,mattsson2008stable,mattsson2012summation,del2015generalized}.   
	However, as we construct our second-derivative FSBP operators in a spectral element setting, we found no advantage in replacing $D_1^T P D_1$ with a more general $A$.  
\end{remark}

\begin{remark} 
	When considering an advection--diffusion equation $\partial_t u + \alpha \partial_x u = \beta \partial_{xx} u$, an energy estimate is obtained when the second-derivative operator $D_2$ has the same norm matrix $P$ as the first-derivative operator $D_1 = P^{-1} Q$, approximating $\partial_x$. 
\end{remark}

%% file: 3_existence.tex
\section{Existence and construction of second-derivative FSBP operators} 
\label{sec:existence} 

We derive necessary and sufficient conditions for the existence of second-derivative FSBP operators. 
We start by connecting the existence of $\mathcal{F}$-exact second-derivative FSBP operators to the existence of $( \mathcal{F} \bigoplus \mathcal{F}')$-exact first derivative operators in \cref{sub:existence_second}, before addressing the existence of the later ones in \cref{sub:existence_first}.
Finally, we provide a simple construction procedure in \cref{sub:construction}.

\subsection{Existence of $\mathcal{F}$-exact second-derivative FSBP operators}
\label{sub:existence_second}

Consider an operator $D_2 = P^{-1}( BS - D_1^T P D_1 )$ approximating $\partial_{xx}$ with symmetric positive define $P$ and $B = \diag(-1,0,\dots,0,1)$. 
We investigate conditions on $D_1$ for which $D_2$ is an $\mathcal{F}$-exact second-derivative FSBP operator. 
To this end, recall that for $D_2 = P^{-1}( BS - D_1^T P D_1 )$ to be an $\mathcal{F}$-exact second-derivative operator, (i) in \cref{def:2nd_FSBP} must hold, which is equivalent to 
\begin{equation}\label{eq:exactness_1} 
	P^{-1}( BS - D_1^T P D_1 ) \mathbf{f} = \mathbf{f''} \quad \forall f \in \mathcal{F}.
\end{equation} 
Since $B S = B D_1$, where $D_1 = P^{-1} Q$ is a first-derivative FSBP operator, and $B = Q + Q^T = P D_1 + D_1^T P$, \cref{eq:exactness_1} becomes 
\begin{equation}\label{eq:exactness_2} 
	D_1 D_1 \mathbf{f} = \mathbf{f''} \quad \forall f \in \mathcal{F}.
\end{equation}
Recall that if $D_1$ is $\mathcal{F}$-exact, then $D_1 \mathbf{f} = \mathbf{f'}$ for all $f \in \mathcal{F}$ and \cref{eq:exactness_2} becomes 
\begin{equation}\label{eq:exactness_3} 
	D_1 \mathbf{f'} = \mathbf{f''} \quad \forall f \in \mathcal{F}. 
\end{equation} 
However, in general, $f' \not\in \mathcal{F}$. 
Consequently, $D_1 \mathbf{f'} = \mathbf{f''}$ only holds if $\mathcal{F}' \subset \mathcal{F}$ (the function space is invariant under differentiation).
While $\mathcal{F}' \subset \mathcal{F}$ is satisfied for polynomials (the derivative of a polynomial is a polynomial of lower degree), it does not necessarily hold for non-polynomial function spaces. 
We discuss such an example in \cref{sub:examples_kernel}. 

To circumvent this issue, we necessitate $D_1$ to be an $(\mathcal{F} \bigoplus \mathcal{F}')$-exact first-derivative FSBP operator. 
Here, $\mathcal{F} \bigoplus \mathcal{F}' = \{ \, f_1 + f_2 \mid f_1 \in \mathcal{F}, f_2 \in \mathcal{F}' \, \}$ is the direct sum of $\mathcal{F}$ and $\mathcal{F}'$, which is the smallest linear function space including both $\mathcal{F}$ and $\mathcal{F}'$. 
In particular, $D_1$ then is exact for all functions from $\mathcal{F}$ and $\mathcal{F}'$ simultaneously.
In this case, $D_1 D_1 \mathbf{f} = D_1 \mathbf{'f} = \mathbf{f''}$ for all $f \in \mathcal{F}$ and the second-derivative FSBP operator $D_2 = P^{-1}( B D_1 + D_1^T P D_1 )$ is $\mathcal{F}$-exact. 
We are now positioned to characterize the existence of $\mathcal{F}$-exact second-derivative operators. 

\begin{theorem}\label{thm:existence_second}
	Let $\mathcal{F} \subset C^2([x_L,x_R])$ and $D_1 = P^{-1} Q$ be a first-derivative FSBP operator with $Q + Q^T = B$. 
	The operator $D_2 = P^{-1}( B D_1 - D_1^T P D_1 )$ approximating $\partial_{xx}$ is an $\mathcal{F}$-exact second-derivative operator if and only if $D_1$ is $(\mathcal{F} \bigoplus \mathcal{F}')$-exact. 
\end{theorem}

\cref{thm:existence_second} immediately follows from the above discussion, and we thus omit its proof.

\subsection{Existence of $(\mathcal{F} \bigoplus \mathcal{F}')$-exact first-derivative FSBP operators}
\label{sub:existence_first}

\cref{thm:existence_second} essentially states that an $\mathcal{F}$-exact second-derivative FSBP operator exists if and only if there exists an $(\mathcal{F} \bigoplus \mathcal{F}')$-exact first-derivative FSBP operator. 
We now address the existence of such first-derivative operators. 
To this end, we can draw upon \cite{glaubitz2022summation} (also refer to \cite{glaubitz2023multi} for the multidimensional case), in which we demonstrated that an $\mathcal{F}$-exact first-derivative FSBP operator $D_1 = P^{-1} Q$ exists if and only if there exists a positive and $(\mathcal{F}^2)'$-exact quadrature. Here, $(\mathcal{F}^2)' = \{ \, f' g + f g' \mid f,g \in \mathcal{F} \, \}$ is the space of all functions corresponding to the derivative of the product of two functions from $\mathcal{F}$. 
Such quadratures were constructed in \cite{glaubitz2021stable,glaubitz2022constructing} under relatively mild conditions, utilizing a least-squares approach.

By replacing $\mathcal{F}$ with $\mathcal{F} \bigoplus \mathcal{F}'$, we can deduce that an $(\mathcal{F} \bigoplus \mathcal{F}')$-exact first-derivative FSBP operator exists if and only if there exists a positive quadrature that is exact for $( [\mathcal{F} \bigoplus \mathcal{F}']^2 )'$. 
It is important to note that the function space $\mathcal{F} \bigoplus \mathcal{F}'$ can be larger than $\mathcal{F}$. 
Consequently, finding a positive quadrature that is exact for $\mathcal{F} \bigoplus \mathcal{F}'$ may necessitate more points than one that is exact for $\mathcal{F}$. 
This, in turn, implies that the first- and second-derivative FSBP operators may require more points, potentially reducing their efficiency. 
This issue can be circumvented if $\mathcal{F}$ is invariant under differentiation since then $\mathcal{F}' \subset \mathcal{F}$ and $\mathcal{F} \bigoplus \mathcal{F}' = \mathcal{F}$. 
For this reason, we advocate using function spaces that are invariant under differentiation whenever possible. 
Polynomials naturally fulfill this condition, and we explore examples of other non-polynomial function spaces possessing this property in \cref{sec:examples}.

\subsection{Construction of first- and second-derivative FSBP operators}
\label{sub:construction} 

Let an arbitrary function space $\mathcal{F} \subset C^2([x_L,x_R])$ be given. 
We now provide a straightforward procedure to construct $\mathcal{F}$-exact second-derivative FSBP operators. 
In \cref{sub:existence_second}, we have showed that an $\mathcal{F}$-exact second-derivative FSBP operators is given by $D_2 = P^{-1}( B D_1 - D_1^T P D_1 )$ if $D_1 = P^{-1} Q$ is an $(\mathcal{F} \bigoplus \mathcal{F}')$-exact first-derivative FSBP operator. 
To construct such a first-derivative FSBP operator $D_1$, we can follow \cite{glaubitz2022summation} (also see \cite{glaubitz2023multi} for the multi-dimensional case and \cite{fernandez2014generalized,hicken2016multidimensional} for polynomial SBP operators), where the following procedure was established: 

First, we construct a positive quadrature that is exact for $( [ \mathcal{F} \bigoplus \mathcal{F}' ]^2 )'$. 
We can always use the least-squares approach \cite{huybrechs2009stable,glaubitz2021stable,glaubitz2022constructing} to find such quadratures. 
In certain situations, more efficient Gaussian quadratures \cite{ma1996generalized,bremer2010nonlinear,huybrechs2022computation} can also be used, requiring fewer grid points.  
We then choose the positive weights $\mathbf{p}$ of this quadrature as the diagonal elements of the norm matrix $P$, i.e., $P = \diag(\mathbf{p})$.

Next, we get $Q$ by re-writing it as $Q = Q_A + B/2$, where $Q_A$ is the anti-symmetric part of $Q$, and recovering $Q_A$ by solving 
\begin{equation}\label{eq:exactness_QA} 
	Q_A V = P V' - \frac{1}{2} B V.
\end{equation}
Here, $V = [\mathbf{g_1},\dots,\mathbf{g_L}]$ and $V' = [\mathbf{g_1'},\dots,\mathbf{g_L'}]$ are the Vandermonde-like matrices for the basis elements $g_1,\dots,g_L$ of $\mathcal{F} \bigoplus \mathcal{F}'$ and their derivatives, respectively. 
One can solve the matrix equation \cref{eq:exactness_QA} by recasting it as a linear system 
\begin{equation}\label{eq:LS_q} 
 	C \mathbf{q} = \mathbf{y}, 
\end{equation}
where $\mathbf{q}$ denotes the vector that contains the strictly lower part of $Q_A$. 
It is computationally convenient to select the unique least-squares solution of \cref{eq:LS_q}, which is the solution with minimal Euclidean norm $\|\cdot\|_{\ell^2}$ among all possible solutions \cite{golub2012matrix}.
The strictly upper part of $Q_A$ is then obtained by $(Q_A)_{j,i} = - (Q_A)_{i,j}$, and the diagonal elements are set to zero. 

Finally, we get an $(\mathcal{F} \bigoplus \mathcal{F}')$-exact first-derivative FSBP operator as $D_1 = P^{-1} Q$, and an $\mathcal{F}$-exact second-derivative FSBP operators as $D_2 = P^{-1}( B D_1 - D_1^T P D_1 )$.

%% file: 4_examples.tex
\section{Examples} 
\label{sec:examples}

We now exemplify the above construction procedure for an array of function spaces, defined on the reference element $\Omega_{\rm ref} = [-1,1]$. 
Our selection includes a polynomial, a trigonometric, and an exponential function space. 
We also illustrate the procedure for a function space incorporating a Gaussian RBF, for which $\mathcal{F}' \not\subset \mathcal{F}$.

\subsection{Polynomial SBP operators} 
\label{sub:examples_poly}

We build a polynomial second-derivative SBP operator on Gauss--Lobatto points to demonstrate that our FSBP theory for general function spaces recover existing polynomial SBP operators. 
To this end, consider the space of polynomials up to degree $d$, 
\begin{equation}\label{eq:expl_poly} 
	\mathcal{F} 
		= \mathcal{P}_d 
		= \Span\{ \, x^k \mid k=0,1,\dots,d \, \}.
\end{equation} 
Observe that $\mathcal{P}_{d}' = \mathcal{P}_{d-1} \subset \mathcal{P}_{d}$ and therefore $\mathcal{F} \bigoplus \mathcal{F}' = \mathcal{F}$. 
A $\mathcal{P}_d$-exact first-derivative SBP operator $D_1 = P^{-1} Q$ on $d+1$ Gauss--Lobatto points is readily available in the literature \cite{gassner2013skew,ranocha2016summation}; 
It can be found by evaluating the derivatives of the corresponding Lagrange basis functions at the Gauss--Lobatto grid points. 
For $d=2$, The associated grid points, quadrature weights (diagonal entries of $P$), and first-derivative SBP operator on $[-1,1]$ are 
\begin{equation} 
	\mathbf{x} = 
	\begin{bmatrix} 
		-1 \\ 0 \\ 1
	\end{bmatrix}, \quad 
	\mathbf{p} = 
	\begin{bmatrix} 
		1/3 \\ 4/3 \\ 1/3
	\end{bmatrix}, \quad 
	D_1 = 
	\begin{bmatrix} 
		-3/2 & 2 & -1/2 \\ 
		-1/2 & 0 & 1/2 \\ 
		1/2 & -2 & 3/2 
	\end{bmatrix}. 
\end{equation} 
Finally, we get an $\mathcal{P}_d$-exact second-derivative SBP operator as $D_2 = P^{-1}( B D_1 - D_1^T B D_1)$. 
For $d=2$ and three Gauss--Lobatto points, as considered above, this yields 
\begin{equation} 
	D_2 = 
	\begin{bmatrix} 
		1 & -2 & 1 \\ 
		1 & -2 & 1 \\ 
		1 & -2 & 1 
	\end{bmatrix}.
\end{equation} 
This demonstrates that our FSBP theory generalizes the existing SBP framework to general function spaces while being able to recover existing polynomial SBP operators.

\subsection{Trigonometric FSBP operators} 
\label{sub:examples_trig}

Consider the trigonometric function space 
\begin{equation}\label{eq:expl_trig} 
	\mathcal{F} 
		= \mathcal{T}_d 
		= \Span\{ \, 1, \sin( k \pi x ), \cos( k \pi x ) \mid k=1,\dots,d \, \}. 
\end{equation}
We first construct an $(\mathcal{F} \bigoplus \mathcal{F}')$-exact first-derivative FSBP operator. 
To this end, note that 
\begin{equation}
	\mathcal{T}_d' = \Span\{ \, \sin( k \pi x ), \cos( k \pi x ) \mid k=1,\dots,d \, \} \subset \mathcal{T}_d
\end{equation}
and therefore $\mathcal{F} \bigoplus \mathcal{F}' = \mathcal{F}$. 
The construction of an $\mathcal{T}_d$-exact first-derivative FSBP operator $D_1 = P^{-1} Q$ was described in \cite{glaubitz2022summation} using the composite trapezoidal rule with $N = 2d + 2$ equidistant points. 
For $d=1$, the associated grid points, quadrature weights, and first-derivative FSBP operator on $[-1,1]$ are 
\begin{equation} 
	\mathbf{x} = 
	\begin{bmatrix} 
		-1 \\ -1/3 \\ 1/3 \\ 1
	\end{bmatrix}, \quad 
	\mathbf{p} = 
	\begin{bmatrix} 
		1/3 \\ 2/3 \\ 2/3 \\ 1/3
	\end{bmatrix}, \quad 
	D_1 = 
	\begin{bmatrix} 
		-1.50 & 1.81 & -1.81 & 1.50 \\ 
   		-0.91 & 0 & 1.81 & -0.91 \\ 
    		0.91 & -1.81 & 0 & 0.91 \\ 
   		-1.50 & 1.81 & -1.81 & 1.50
	\end{bmatrix}, 
\end{equation} 
where we have rounded the numbers to the second decimal place for $D_1$. 
Finally, we get an $\mathcal{T}_d$-exact second-derivative FSBP operator as $D_2 = P^{-1}( B D_1 - D_1^T B D_1 )$. 
For $d=1$, this yields 
\begin{equation} 
	D_2 = 
	\begin{bmatrix} 
		-3.29 & 3.29 & 3.29 & -3.29 \\
    		4.37 & -6.58 & 3.29 & -1.08 \\
   		-1.08 & 3.29 & -6.58 & 4.37 \\
   		-3.29 & 3.29 & 3.29 & -3.29  
	\end{bmatrix},
\end{equation} 
where we again rounded the numbers to the second decimal place for $D_2$.

\subsection{Exponential SBP operators} 
\label{sub:examples_exp} 

Consider the combined polynomial and exponential function space 
\begin{equation}\label{eq:expl_exp} 
	\mathcal{F} 
		= \mathcal{E}_d 
		= \Span\{ \, 1,x, \dots, x^{d-1}, e^x \,\}.
\end{equation}
We again first construct an $(\mathcal{F} \bigoplus \mathcal{F}')$-exact first-derivative FSBP operator. 
Observe that 
\begin{equation}
	\mathcal{E}_d' 
		= \Span\{ \, 1,x, \dots, x^{d-2}, e^x \,\} 
		\subset \mathcal{E}_d
\end{equation}
and therefore $\mathcal{F} \bigoplus \mathcal{F}' = \mathcal{F}$. 
The construction of an $\mathcal{E}_d$-exact first-derivative FSBP operator was described in \cite{glaubitz2022summation} using a least-squares quadrature with equidistant points. 
For $d=2$, the associated grid points, quadrature weights, and first-derivative FSBP operator on $[-1,1]$ are 
\begin{equation}\label{eq:examples_exp_D1} 
	\mathbf{x} = 
	\begin{bmatrix} 
		-1 \\ -1/2 \\ 0 \\ 1/2 \\ 1
	\end{bmatrix}, \quad 
	\mathbf{p} = 
	\begin{bmatrix} 
		0.14 \\ 0.77 \\ 0.19 \\ 0.75 \\ 0.15
	\end{bmatrix}, \quad 
	D_1 = 
	\begin{bmatrix} 
		-3.64 & 4.97 & -0.48 & -1.38 & 0.53 \\
   		-0.88 & 0 & 0.41 & 0.72 & -0.24 \\
    		0.35 & -1.65 & 0 & 1.56 & -0.25 \\
    		0.25 & -0.74 & -0.40 & 0 & 0.88 \\
   		-0.50 & 1.27 & 0.33 & -4.5 & 3.39 
	\end{bmatrix}. 
\end{equation} 
In \cref{eq:examples_exp_D1}, we have rounded the numbers to the second decimal place for $\mathbf{p}$ and $D_1$. 
Finally, we get an $\mathcal{E}_d$-exact second-derivative FSBP operator as $D_2 = P^{-1}( B D_1 - D_1^T B D_1 )$. 
For $d=2$, this yields 
\begin{equation} 
	D_2 = 
	\begin{bmatrix} 
		8.07 & -15.59 & 4.53 & 5.42 & -2.43 \\
    		3.66 & -5.91 & 0.06 & 2.95 & -0.77 \\
    		0.72 & 0.25 & -1.55 & -0.51 & 1.09 \\
   		-0.84 & 3.03 & -0.13 & -5.47 & 3.41 \\ 
   		-2.02 & 4.61 & 3.68 & -13.13 & 6.85 \\ 
	\end{bmatrix},
\end{equation} 
where we again rounded the numbers to the second decimal place for $D_2$.

\subsection{Radial basis function SBP operators} 
\label{sub:examples_kernel}

All of the above function spaces satisfied $\mathcal{F}' \subset \mathcal{F}$, simplifying the construction of the first-derivative FSBP operator.  
We conclude our examples with a function space for which $\mathcal{F}' \not\subset \mathcal{F}$.
To this end, consider the three-dimensional function space 
\begin{equation}\label{eq:expl_RBF} 
	\mathcal{F} 
		= \mathcal{R}_3 
		= \Span\{ \, 1, x, e^{x^2} \, \}, 
\end{equation}
containing a Gaussian radial basis function (RBF) centered at zero. 
We first construct an $(\mathcal{F} \bigoplus \mathcal{F}')$-exact first-derivative FSBP operator. 
In contrast to the previous function spaces, we have 
\begin{equation}
	\mathcal{R}_3' = \Span\{ \, 1, x e^{x^2} \, \} \not\subset \mathcal{R}_3,  
\end{equation}
and hence 
\begin{equation}
	\mathcal{R}_3 \bigoplus \mathcal{R}_3' 
		= \Span\{ \, 1, x, e^{x^2}, x e^{x^2} \, \},  
\end{equation} 
which is a larger space than $\mathcal{R}_3$. 
Importantly, if we constructed a $\mathcal{R}_3$-exact first-derivative FSBP operator $D_1$, the resulting second-derivative FSBP operator $D_2 = P^{-1}( B D_1 - D_1^T B D_1 )$ would \emph{not} be $\mathcal{R}_3$-exact. 
See the discussion in \cref{sub:existence_second}. 
To circumvent this issue, and get an $\mathcal{R}_3$-exact second-derivative FSBP operator, we instead construct an $(\mathcal{R}_3 \bigoplus \mathcal{R}_3')$-exact first-derivative FSBP operator. 
To this end, we again use a least-squares quadrature with equidistant points. 
The associated grid points, quadrature weights, and first-derivative FSBP operator on $[-1,1]$ are 
\begin{equation} 
	\mathbf{x} = 
	\begin{bmatrix} 
		-1 \\ -1/2 \\ 0 \\ 1/2 \\ 1
	\end{bmatrix}, \quad 
	\mathbf{p} = 
	\begin{bmatrix} 
		0.20 \\ 0.58 \\ 0.44 \\ 0.58 \\ 0.20
	\end{bmatrix}, \quad 
	D_1 = 
	\begin{bmatrix} 
		-2.45 & 3.13 & -0.57 & -0.45 & 0.35 \\
   		-1.11 & 0 & 1.16 & 0.10 & -0.16 \\
    		0.27 & -1.54 & 0 & 1.54 & -0.27 \\
    		0.16 & -0.10 & -1.16 & 0 & 1.11 \\
   		-0.35 & 0.45 & 0.57 & -3.13 & 2.45 \\
	\end{bmatrix}, 
\end{equation} 
where we have rounded the numbers to the second decimal place for $\mathbf{p}$ and $D_1$. 
Finally, we get an $\mathcal{R}_3$-exact second-derivative FSBP operator as $D_2 = P^{-1}( B D_1 - D_1^T B D_1)$, yielding 
\begin{equation} 
	D_2 = 
	\begin{bmatrix} 
		2.17 & -6.56 & 5.77 & -0.53 & -0.85 \\
    		3.10 & -5.34 & 0.42 & 2.79 & -0.97 \\
    		1.39 & 0.56 & -3.89 & 0.56 & 1.39 \\
   		-0.97 & 2.79 & 0.42 & -5.34 & 3.10 \\
   		-0.85 & -0.53 & 5.77 & -6.56 & 2.17 \\
	\end{bmatrix},
\end{equation} 
where we again rounded the numbers to the second decimal place for $D_2$.

\subsection{Nullspaces of FSBP operators}

We briefly investigate the nullspace of the above FSBP operators.  

\begin{definition}[The nullspace]
	Let $L: X \to Y$ be a linear operator between the two vector spaces $X$ and $Y$. 
	The \emph{nullspace} of $L$ is the set of all vectors $x \in X$ that are mapped to zero by $L$, which we denote by $\operatorname{null}(L) = \{ \, x \in X \mid Lx = 0 \, \}$.
\end{definition}

The nullspaces of the continuous first- and second-derivative operators $\partial_x, \partial_{xx}: \mathcal{F} \to \mathcal{F}$, where we have restricted them to a function space $\mathcal{F} \subset C^2$, are 
\begin{equation} 
\begin{aligned} 
	\operatorname{null}( \partial_x ) 
		& = \{ \, f \in \mathcal{F} \mid \partial_x f \equiv 0 \, \}, \\ 
	\operatorname{null}( \partial_{xx} ) 
		& = \{ \, f \in \mathcal{F} \mid \partial_{xx} f \equiv 0 \, \}.  
\end{aligned}
\end{equation} 
If $\Span\{ 1 \} \subset \mathcal{F}$, then $\operatorname{null}( \partial_x ) = \Span\{ 1 \}$, and if $\Span\{ 1, x \} \subset \mathcal{F}$, then $\operatorname{null}( \partial_{xx} ) = \Span\{ 1, x \}$. 
It is often desirable---and sometimes even necessary \cite{svard2019convergence}---for a discrete derivative operator to have, on a discrete level, the same nullspace as its continuous counterpart. 
In this case, we say that the discrete operator is \emph{nullspace consistent}. 
For polynomial SBP operators, nullspace consistency has been studied in, for instance, \cite{linders2019convergence,ranocha2019some,svard2019convergence,ranocha2020discrete,ranocha2021new,linders2022eigenvalue}. 
Here, we report on the nullspaces in the examples of the above-derived first- and second-derivative FSBP operators.

\subsubsection*{Polynomial SBP operators}

We start by considering traditional polynomial SBP operators that are exact for polynomials up to degree two, described in \cref{sub:examples_poly}. 
In this case, $\mathcal{F} = \Span\{ 1, x, x^2 \}$ and the nullspace of the continuous first- and second-derivative operator, $\partial_x: \mathcal{F} \to \mathcal{F}$ and $\partial_{xx}: \mathcal{F} \to \mathcal{F}$, is $\operatorname{null}( \partial_x ) = \{ 1 \}$ and $\operatorname{null}( \partial_{xx} ) = \{ 1, x \}$, respectively. 
The nullspaces of the corresponding SBP operators are $\operatorname{null}( D_1 ) = \{ \mathbf{1} \}$ and $\operatorname{null}( D_2 ) = \{ \mathbf{1}, \mathbf{x} \}$. 
This implies that both operators are nullspace consistent.

\subsubsection*{Trigonometric FSBP operators}

Consider the trigonometric FSBP operators that are exact for $\mathcal{F} = \Span\{ 1, \sin(\pi x), \cos(\pi x) \}$, described in \cref{sub:examples_trig}. 
In this case, the nullspaces of the continuous first- and second-derivative operator, $\partial_x, \partial_{xx}: \mathcal{F} \to \mathcal{F}$ and $\partial_{xx}$, are $\operatorname{null}( \partial_x ) = \operatorname{null}( \partial_{xx} ) = \{ 1 \}$. 
The nullspaces of the corresponding SBP operators are $\operatorname{null}( D_1 ) = \{ \mathbf{1} \}$ and $\operatorname{null}( D_2 ) = \{ \mathbf{1}, \mathbf{v} \}$ with $\mathbf{v} = [ 0, 0.22, 0.77, 1 ]^T$, where we have rounded the numbers to the second decimal place.\footnote{We computed the null space using the command ``null" in Matlab.}
This implies that the first-derivative trigonometric FSBP operator is nullspace consistent, while the second-derivative operator is not.

\subsubsection*{Exponential SBP operators}

Consider the exponential FSBP operators that are exact for $\mathcal{F} = \Span\{ 1, x, e^x \}$, described in \cref{sub:examples_exp}. 
In this case, the nullspaces of the continuous first- and second-derivative operator, $\partial_x, \partial_{xx}: \mathcal{F} \to \mathcal{F}$ and $\partial_{xx}$, are $\operatorname{null}( \partial_x ) = \{ 1 \}$ and $\operatorname{null}( \partial_{xx} ) = \{ 1, x \}$. 
The nullspaces of the corresponding FSBP operators are $\operatorname{null}( D_1 ) = \{ \mathbf{1} \}$ and $\operatorname{null}( D_2 ) = \{ \mathbf{1}, \mathbf{x} \}$, which are the discrete analogues to the continuous case. 
Hence, the first- and second-derivative exponential FSBP operators are both nullspace consistent.

\subsubsection*{Radial basis function SBP operators}
\cref{sub:examples_kernel}

Finally, we consider the Gaussian RBF FSBP operators that are exact for $\mathcal{F} = \Span\{ 1, x, e^{x^2} \}$, described in \cref{sub:examples_kernel}. 
In this case, the nullspaces of the continuous first- and second-derivative operator, $\partial_x, \partial_{xx}: \mathcal{F} \to \mathcal{F}$ and $\partial_{xx}$, are again $\operatorname{null}( \partial_x ) = \{ 1 \}$ and $\operatorname{null}( \partial_{xx} ) = \{ 1, x \}$. 
The nullspaces of the corresponding FSBP operators are $\operatorname{null}( D_1 ) = \{ \mathbf{1} \}$ and $\operatorname{null}( D_2 ) = \{ \mathbf{1}, \mathbf{x} \}$, which are the discrete analogues to the continuous case. 
Hence, the first- and second-derivative Gaussian RBF FSBP operators are again both nullspace consistent. 

%% file: 5_numerics.tex
\section{Numerical tests} 
\label{sec:numerics}

We demonstrate the versatility of the developed second-derivative FSBP operators by examining their application across various function spaces in selected problems. 
While these problems were chosen to exemplify the utility of certain function spaces, we do not intend to claim that these spaces are the definitive optimal choices for these problems. 
Instead, our primary focus is on validating the efficacy and adaptability of the proposed second-derivative FSBP operators when applied to a diverse array of non-polynomial approximation spaces.
For all subsequent numerical tests, we used the explicit strong stability preserving (SSP) Runge--Kutta (RK) method of third order using three stages (SSPRK(3,3)) \cite{shu1988total} with time step size scaling $1/\Delta t \sim \lambda_{\max}/\Delta x + \varepsilon/ \Delta x^2$, where $\lambda_{\max}$ is the largest characteristic velocity, $\varepsilon$ is the diffusivity, and $\Delta x$ is the spatial resolution. 
For simplicity, we only consider FSBP operators on equidistant grid points. 
The MATLAB code used to generate the numerical tests presented here is open access and can be found on GitHub (\url{https://github.com/jglaubitz/2ndDerivativeFSBP}).

\subsection{The linear advection--diffusion equation: The single-block case} 
\label{sub:periodic}

Consider the linear advection-diffusion equation 
\begin{equation}\label{eq:periodc_advDif} 
	\partial_t u + a \partial_x u = \varepsilon \partial_{xx} u  
\end{equation} 
on $\Omega = [-1,1]$ with periodic boundary conditions and initial data $u(x,0) = \cos(4\pi x) + 0.75 \sin(40 \pi x)$. 
The analytical solution of \cref{eq:periodc_advDif} is 
\begin{equation}
	u(x,t) = e^{- \varepsilon (4 \pi)^2 t} \cos(4\pi [x-at]) + 0.75 e^{- \varepsilon (40 \pi)^2 t} \sin(40 \pi [x-at]). 
\end{equation} 
Let $D_1 = P^{-1} Q$ and $D_2 = P^{-1}( BS - D_1^T P D_1 )$ be 'global' first- and second-derivative FSBP operators on the physical domain $\Omega$. 
In this case, we use a pseudo-spectral approach in which the reference element $\Omega_{\rm ref}$, on which the FSBP operators are constructed, is equal to $\Omega$.
Then, the single-block FSBP-SAT semi-discretization of \cref{eq:periodc_advDif} is 
\begin{equation}\label{eq:periodc_advDif_discr} 
	\mathbf{u}_t + a D_1 \mathbf{u} 
		= \varepsilon D_2 \mathbf{u} + P^{-1} \mathbb{S}.
\end{equation}
Here, $\mathbf{u} = [u_1,\dots,u_N]^T$ denotes the vector of the nodal values of the numerical solution at the grid points $\mathbf{x} = [x_1,\dots,x_N]^T$.
Furthermore, $\mathbb{S} = \mathbb{S}_L + \mathbb{S}_R$, where $\mathbb{S}_L$ and $\mathbb{S}_R$ are \emph{simultaneous approximation terms} (SATs) at the left and right boundary of the single block, respectively, that weakly enforce the boundary condition. 
We follow \cite{carpenter2010revisiting,gong2011interface} and choose the SATs according to the Baumann--Oden method \cite{baumann1999discontinuous} as 
\begin{equation}\label{eq:periodc_advDif_SATs} 
\begin{aligned}
	\mathbb{S}_L 
		& = \sigma_1^{L} \mathbf{e}_L [ u_1 - u_N ] 
			+ \sigma_2^{L} \mathbf{e}_L [ (D_1 \mathbf{u})_1 - (D_1 \mathbf{u})_N ] 
			+ \sigma_3^{L} D_1^T \mathbf{e}_L [ u_1 - u_N ], \\ 
	\mathbb{S}_R 
		& = \sigma_1^{R} \mathbf{e}_R [ u_N - u_1 ] 
			+ \sigma_2^{R} \mathbf{e}_R [ (D_1 \mathbf{u})_N - (D_1 \mathbf{u})_1 ] 
			+ \sigma_3^{R} D_1^T \mathbf{e}_R [ u_N - u_1 ], 
\end{aligned}
\end{equation}
where $\mathbf{e}_L = [1,0,\dots,0]^T$ and $\mathbf{e}_R = [0,\dots,0,1]^T$.
It was proved in \cite{carpenter2010revisiting} that the FSBP-SAT semi-discretization \cref{eq:periodc_advDif_discr} with SATs as in \cref{eq:periodc_advDif_SATs} is stable if $\sigma_1^R \leq a/2$, $\sigma_1^L = \sigma_1^R -a$, $\sigma_2^R \in \R$, $\sigma_2^L = \varepsilon + \sigma_2^R$, $\sigma_3^R = -\varepsilon - \sigma_2^R$, and $\sigma_3^L = -\sigma_2^R$. 
The two free parameters are $\sigma_1^R$ and $\sigma_2^R$. 
We chose $\sigma_1^R = 0$ and $\sigma_2^R = -\varepsilon/2$ for the subsequent numerical tests.
We also experimented with other choices; 
However, for brevity, we do not present these here. 

\begin{remark} 
	The idea behind SATs is to \emph{simultaneously approximate} the equation and the boundary conditions---as well as inter-element coupling in a multi-block approach---by ``pulling" the numerical solution towards the boundary data. 
	While there are various methods to incorporate boundary conditions, we emphasize that the SBP-SAT framework is versatile and can be applied to ensure stability across various schemes. 
	See \cite{svard2014review,fernandez2014review,abgrall2020analysis,abgrall2021analysis} and references therein. 
	In fact, one can also mix different types of discretization schemes in the same solver provided that they are based on SBP operators, as demonstrated in \cite{nordstrom2006stable,gong2007stable,nordstrom2009hybrid}.
	Another notable strength of the SBP-SAT framework is its ability to incorporate any continuous BC that contributes to an energy estimate in a stable way. 
\end{remark}

\begin{figure}[tb]
	\centering 
	\begin{subfigure}[b]{0.49\textwidth}
		\includegraphics[width=\textwidth]{%
      		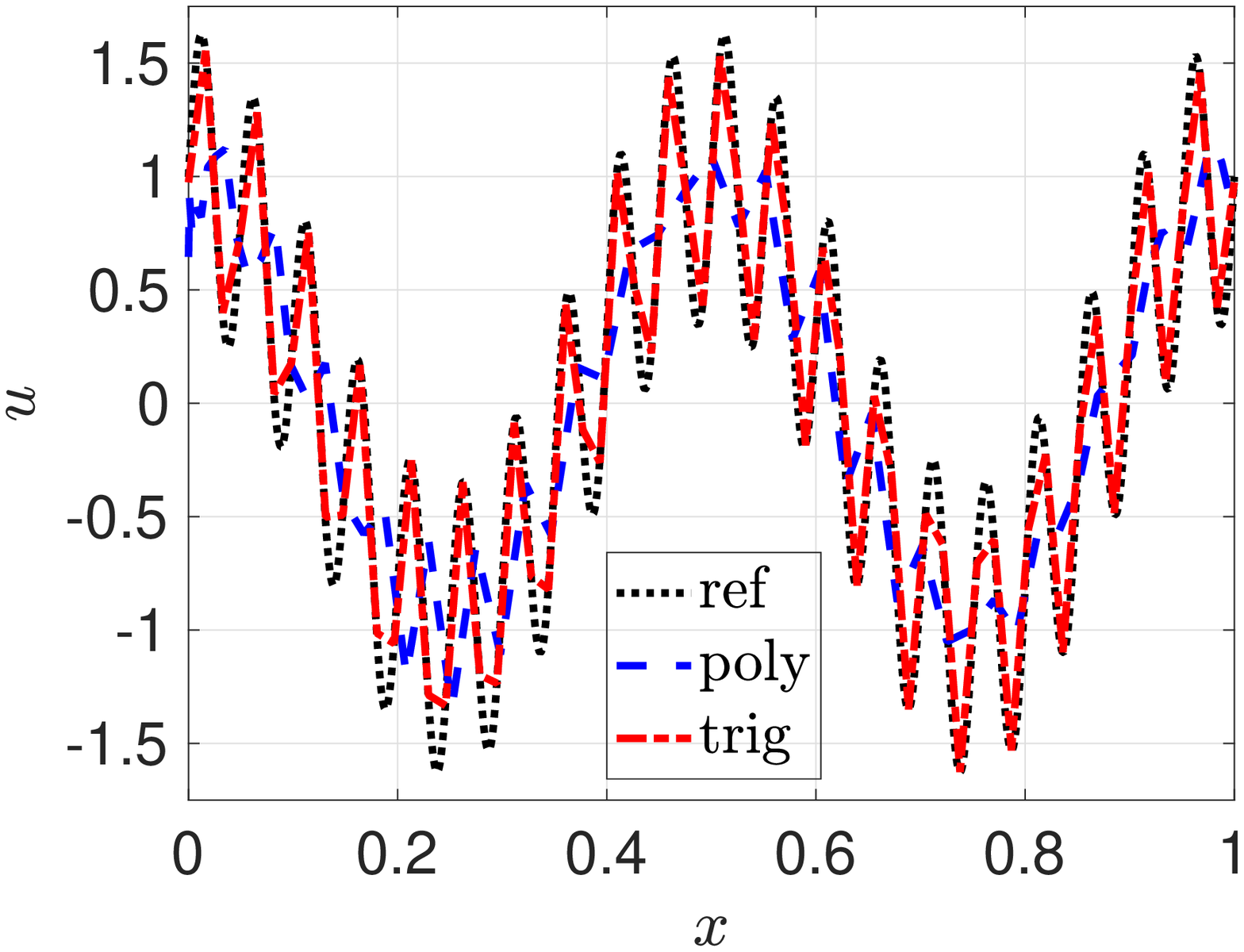} 
    		\caption{Solutions at $t=1$}
    		\label{fig:linear_advDiff_periodic_sol_d30_t1_eps5}
  	\end{subfigure}%
	\begin{subfigure}[b]{0.49\textwidth}
		\includegraphics[width=\textwidth]{%
      		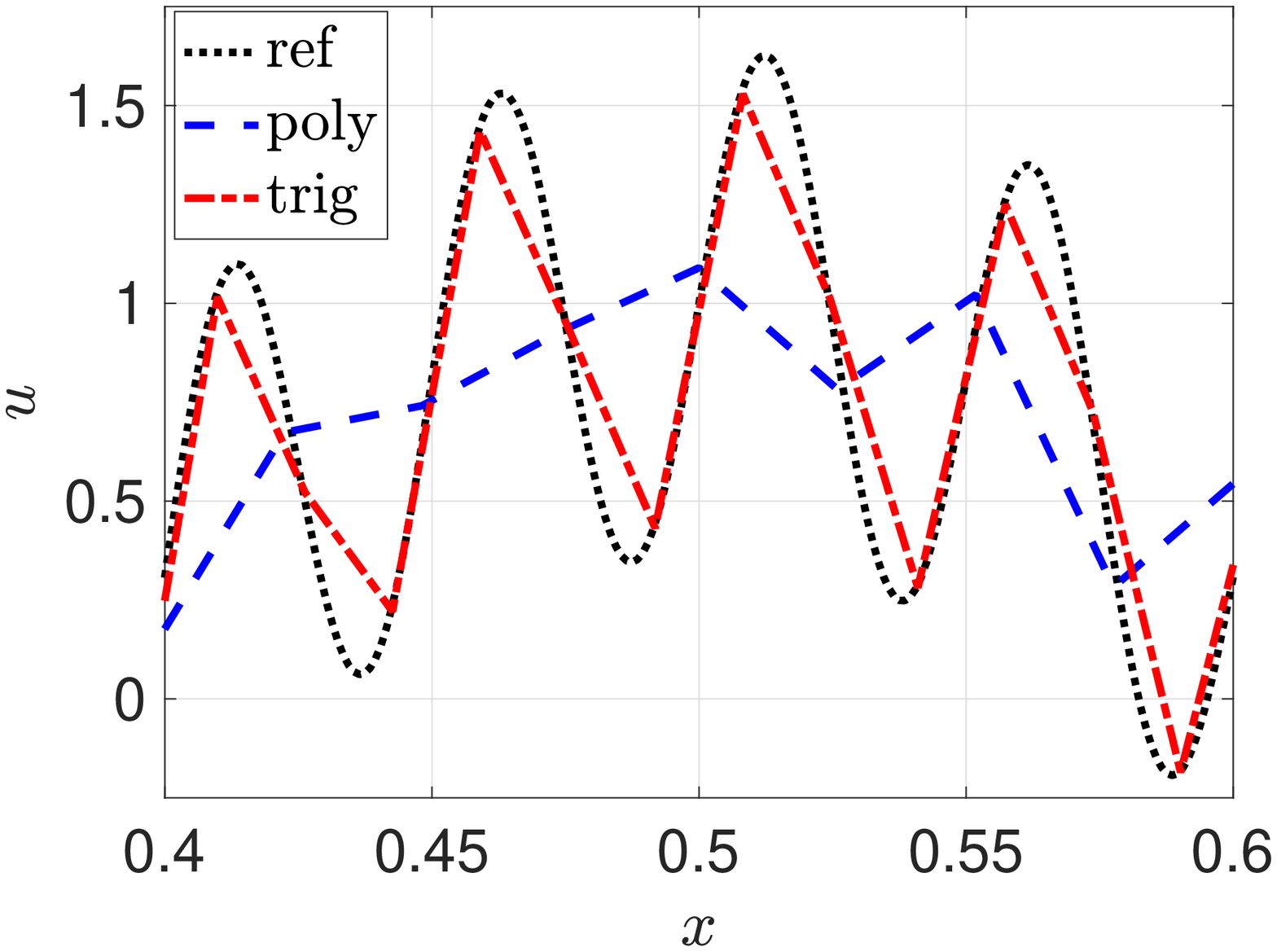} 
    		\caption{Zoomed-in solutions at $t=1$}
		    		\label{fig:linear_advDiff_periodic_solZoomed_d30_t1_eps5}
  	\end{subfigure}%
	\\
  	\begin{subfigure}[b]{0.49\textwidth}
		\includegraphics[width=\textwidth]{%
      		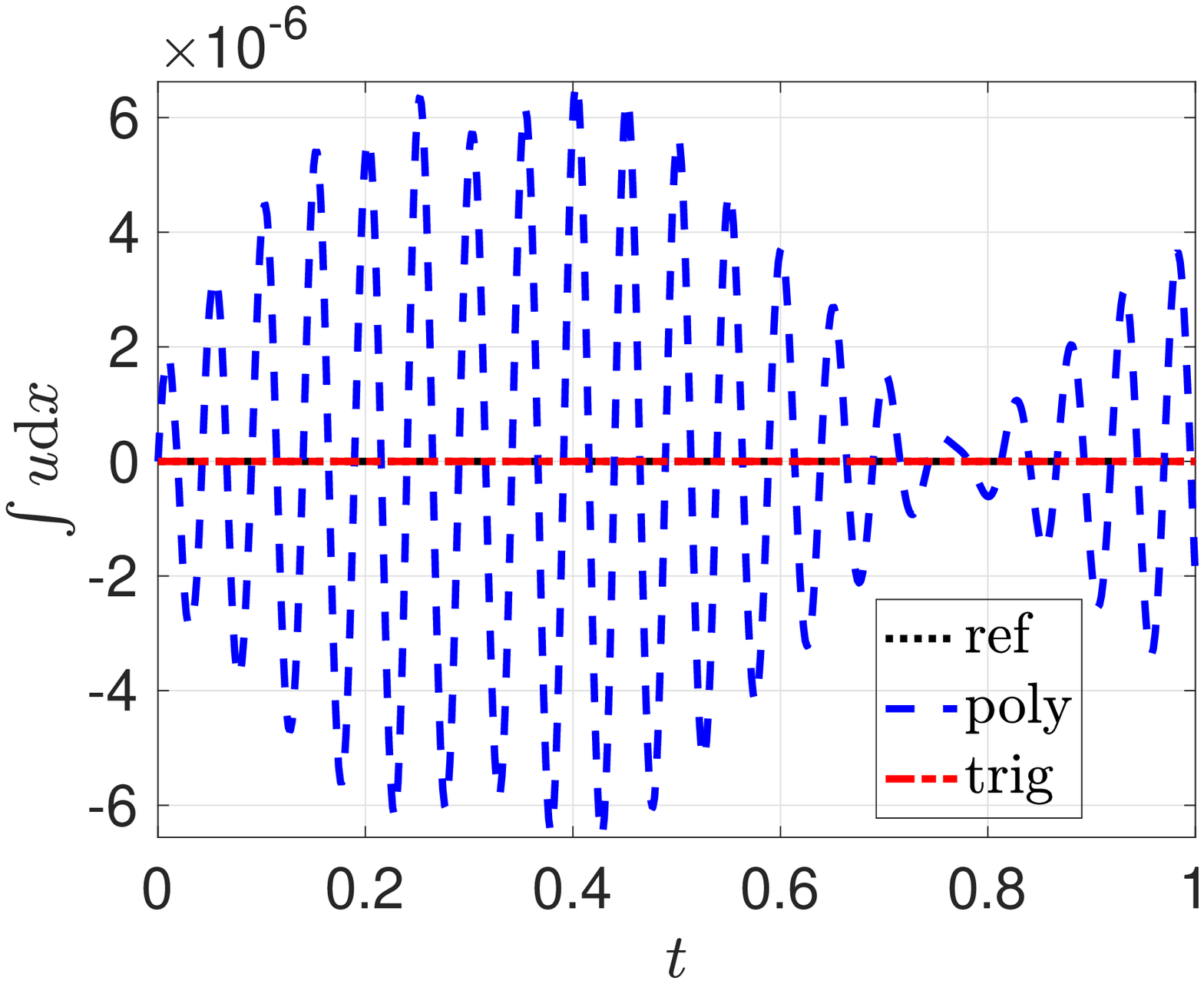} 
    		\caption{Mass over time}
    		\label{fig:linear_advDiff_periodic_mass_d30_t1_eps5}
  	\end{subfigure}%
  	\begin{subfigure}[b]{0.49\textwidth}
		\includegraphics[width=\textwidth]{%
      		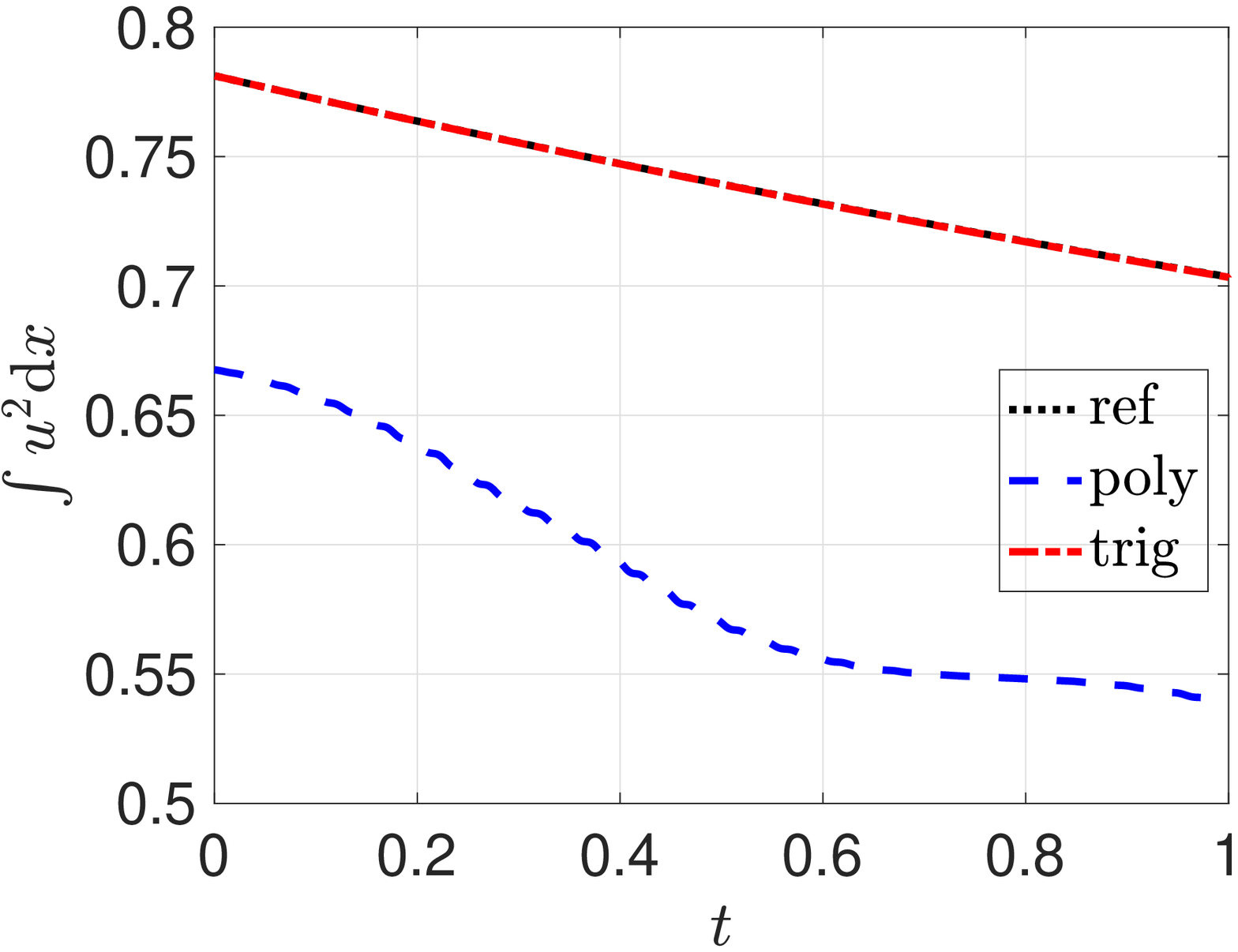} 
    		\caption{Energy over time}
    		\label{fig:linear_advDiff_periodic_energy_d30_t1_eps5}
  	\end{subfigure}%
  	\caption{
  	(Numerical) solutions of the linear advection--diffusion equation \cref{eq:periodc_advDif} for $a=1$ and $\varepsilon = 10^{-5}$ with periodic initial and boundary data at time $t=1$ and mass and energy over time. 
	The numerical solutions correspond to the FSBP-SAT method \cref{eq:periodc_advDif_discr} using the polynomial (``poly") and trigonometric (``trig") approximation space, $\mathcal{P}_{60}$ and $\mathcal{T}_{30}$. 
	Both spaces have the same dimension of $K=61$ and use $N=61$ (Gauss--Lobatto) and $N=62$ (equidistant) grid points, respectively. 
	The mass and energy profiles of the reference and numerical solutions using the trigonometric FSBP operators align. 
  	}
  	\label{fig:linear_advDiff_periodic}
\end{figure}   

\cref{fig:linear_advDiff_periodic} provides a comparison between the exact solution and numerical solutions of the periodic linear advection-diffusion problem \cref{eq:periodc_advDif} for $a=1$ and $\varepsilon = 10^{-5}$ at time $t=1$. 
The numerical solutions correspond to the FSBP-SAT method \cref{eq:periodc_advDif_discr} using the polynomial (``poly") and trigonometric (``trig") approximation space, $\mathcal{P}_{60}$ and $\mathcal{T}_{30}$. 
Both spaces have the same dimension of $K=61$ and use $N=61$ (Gauss--Lobatto) and $N=62$ (equidistant) grid points, respectively. 
Notably, the initial condition includes the term ``$\sin(40\pi x)$", while the basis functions with the highest frequency in $\mathcal{T}_{30}$ are ``$\sin(30\pi x)$" and ``$\cos(30\pi x)$." 
Hence, the trigonometric approximation space $\mathcal{T}_{30}$ does not exactly represent the initial condition or solution. 
\cref{fig:linear_advDiff_periodic_sol_d30_t1_eps5} 
shows the corresponding numerical solutions.  
\cref{fig:linear_advDiff_periodic_mass_d30_t1_eps5,fig:linear_advDiff_periodic_energy_d30_t1_eps5} illustrate the corresponding mass ($\int u \intd x$) and energy ($\int u^2 \intd x$) over time. 
The trigonometric FSBP-SAT scheme is significantly more accurate than the polynomial-based SBP-SAT scheme.

\subsection{The linear advection--diffusion equation: The multi-block cases} 
\label{sub:kernel}

We again consider the linear advection--diffusion equation \cref{eq:periodc_advDif} on the global computational domian $\Omega = [-1,1]$ with periodic boundary conditions and initial condition $u(x,0) = \cos(4\pi x) + 2 \sin(10 \pi x)$. 
The analytical solution is  
\begin{equation}
	u(x,t) = e^{- \varepsilon (4 \pi)^2 t} \cos(4\pi [x-at]) + 2 e^{- \varepsilon (10 \pi)^2 t} \sin(10 \pi [x-at]). 
\end{equation}
This time, we use a multi-block FSBP-SAT semi-discretization on $I$ uniform blocks,
\begin{equation}\label{eq:advDif_discr} 
	\mathbf{u}_t^{(i)} + a D_1 \mathbf{u}^{(i)} 
		= \varepsilon D_2 \mathbf{u}^{(i)} + P^{-1} \mathbb{S}^{(i)}, \quad i=1,\dots,I,
\end{equation}
where $\mathbf{u}^{(i)}$ corresponds to the numerical solution on the $i$-th block and $\mathbb{S}^{(i)} = \mathbb{S}_L^{(i)} + \mathbb{S}_R^{(i)}$. 
In this case, we construct the FSBP operators on the reference element $\Omega_{\rm ref} = [-1,1]$, where all computations are carried out before it is diffeomorphically mapped to the different blocks that partition $\Omega$. 
The multi-block SATs are 
\begin{equation}\label{eq:advDif_SATs}  
\begin{aligned}
	\mathbb{S}_L^{(i)} 
		& = \sigma_1^{L} \mathbf{e}_L [ u_1^{(i)} - u_N^{(i-1)} ] 
			+ \sigma_2^{L} \mathbf{e}_L [ (D_1 \mathbf{u}^{(i)})_1 - (D_1 \mathbf{u}^{(i-1)})_N ] 
			+ \sigma_3^{L} D_1^T \mathbf{e}_L [ u_1^{(i)} - u_N^{(i-1)} ], \\ 
	\mathbb{S}_R^{(i)} 
		& = \sigma_1^{R} \mathbf{e}_R [ u_N^{(i)} - u_1^{(i+1)} ] 
			+ \sigma_2^{R} \mathbf{e}_R [ (D_1 \mathbf{u}^{(i)})_N - (D_1 \mathbf{u}^{(i+1)})_1 ] 
			+ \sigma_3^{R} D_1^T \mathbf{e}_R [ u_N^{(i)} - u_1^{(i+1)} ], 
\end{aligned}
\end{equation}
for $i=1,\dots,I$, where $u^{(0)} = u^{(I)}$ and $u^{(I+1)} = u^{(1)}$ (to weakly enforce the periodic boundary conditions).
Similar to before, we chose $\sigma_1^R = 0$, $\sigma_2^R = -\varepsilon/2$, $\sigma_1^L = \sigma_1^R -a$, $\sigma_2^L = \varepsilon + \sigma_2^R$, $\sigma_3^R = -\varepsilon - \sigma_2^R$, and $\sigma_3^L = -\sigma_2^R$.

\begin{figure}[tb]
	\centering 
	\begin{subfigure}[b]{0.49\textwidth}
		\includegraphics[width=\textwidth]{%
      		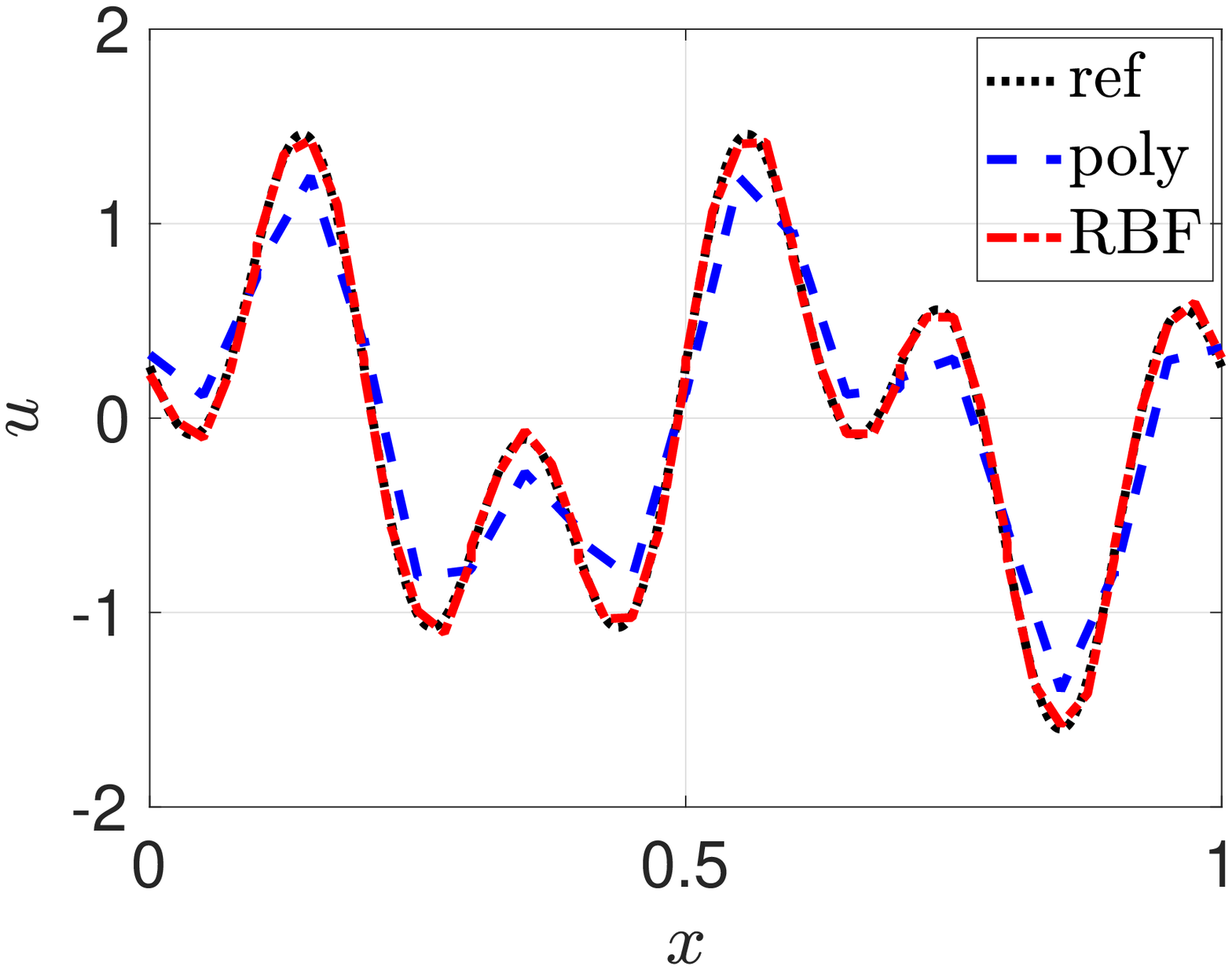} 
    		\caption{Solutions for $I=10$ blocks}
    		\label{fig:linear_AdvDiff_kernel_sol_I10}
  	\end{subfigure}%
  	\begin{subfigure}[b]{0.49\textwidth}
		\includegraphics[width=\textwidth]{%
      		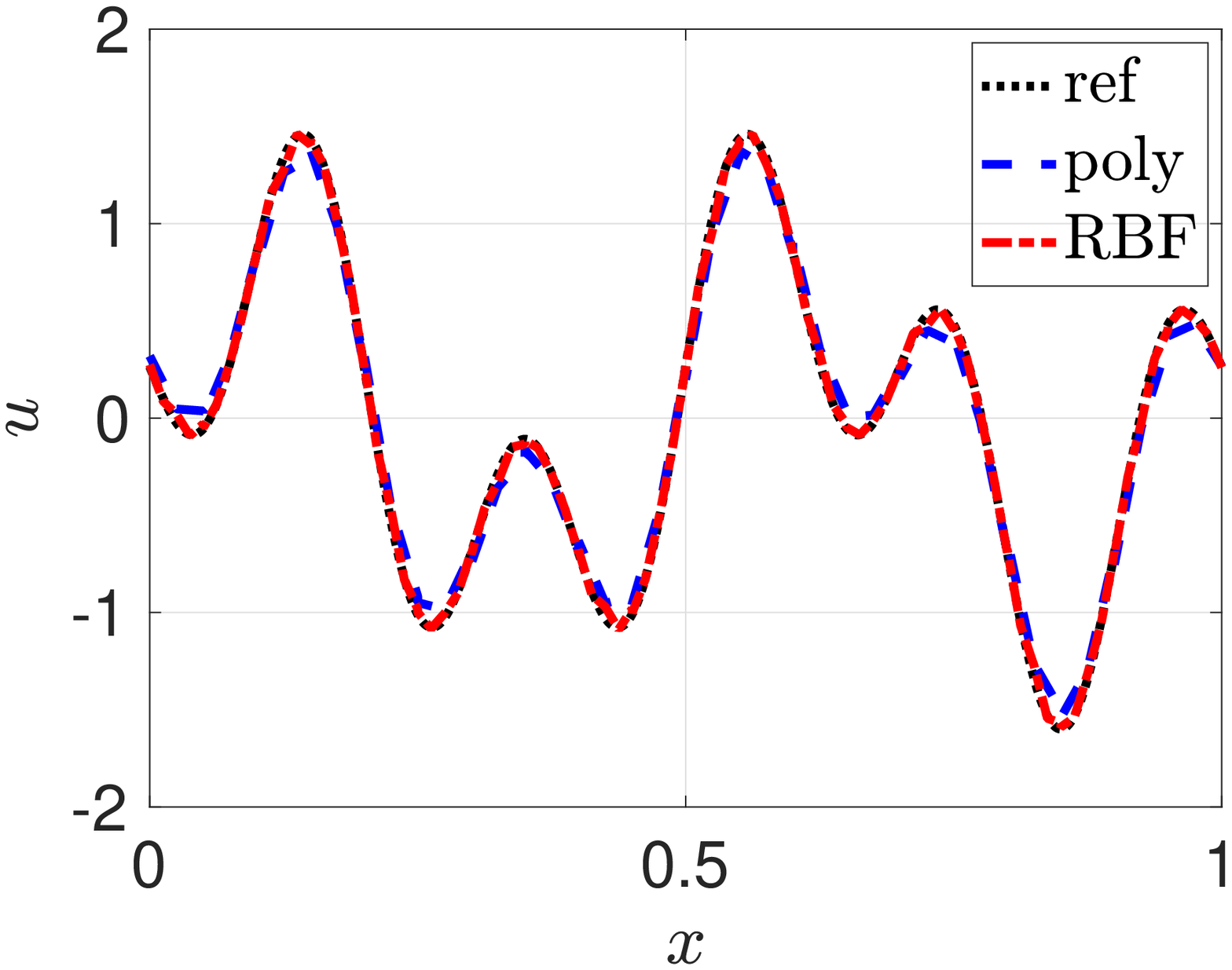} 
    		\caption{Solutions for $I=20$ blocks}
    		\label{fig:linear_AdvDiff_kernel_sol_I20}
  	\end{subfigure}%
	\\
	\begin{subfigure}[b]{0.49\textwidth}
		\includegraphics[width=\textwidth]{%
      		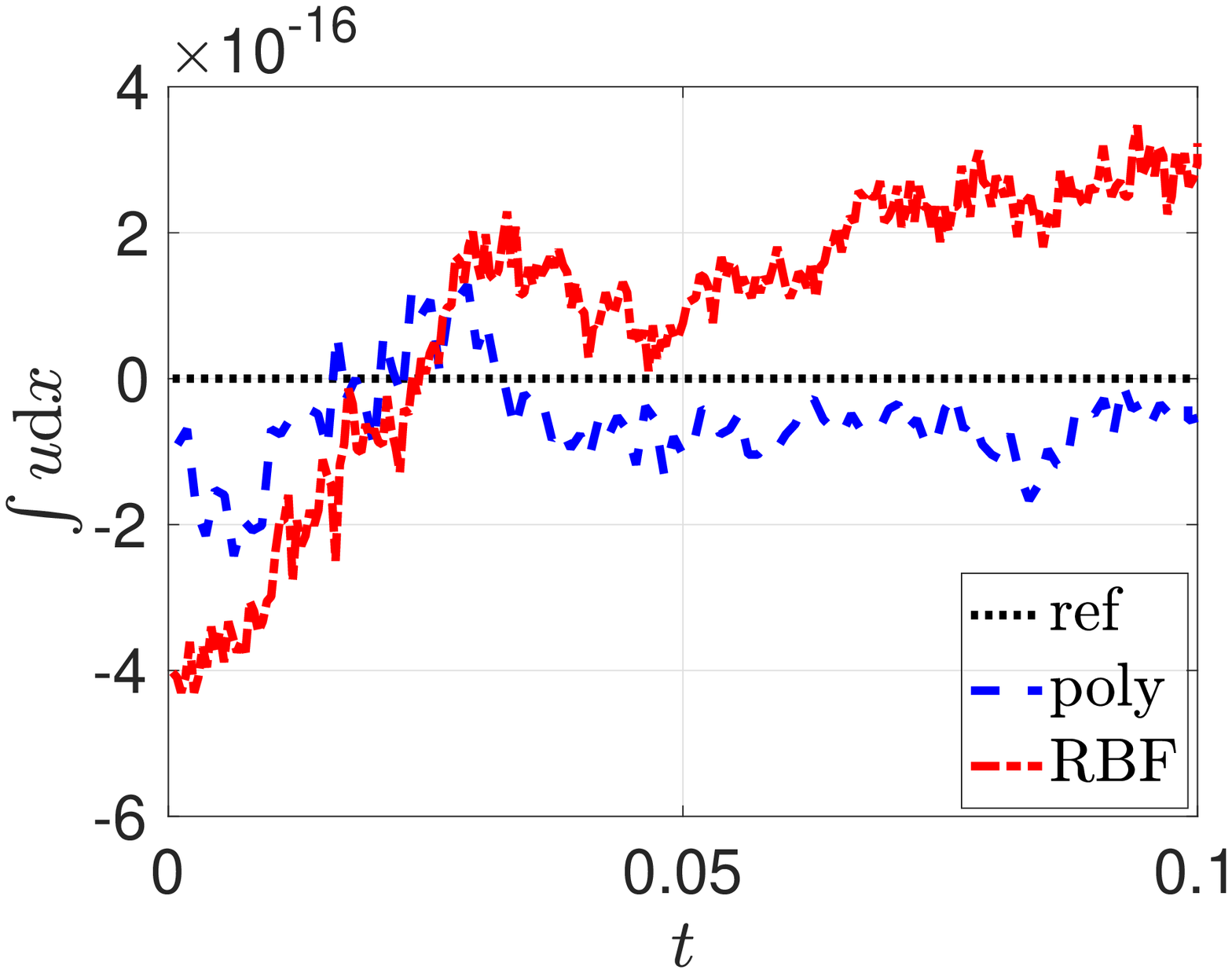} 
    		\caption{Mass over time for $I=10$}
    		\label{fig:linear_AdvDiff_kernel_mass_I10}
  	\end{subfigure}%
  	\begin{subfigure}[b]{0.49\textwidth}
		\includegraphics[width=\textwidth]{%
      		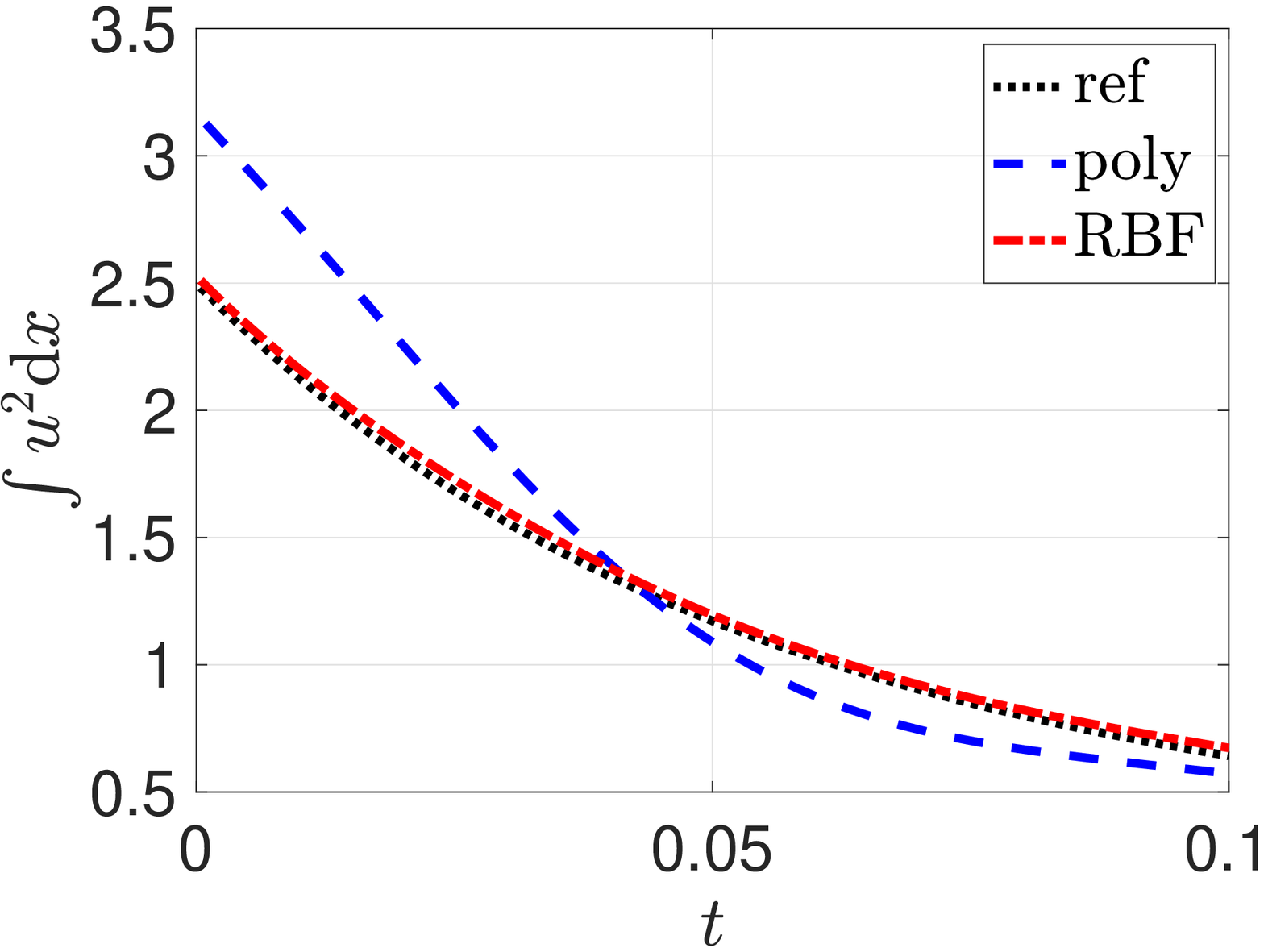} 
    		\caption{Energy over time for $I=10$}
    		\label{fig:linear_AdvDiff_kernel_energy_I10}
  	\end{subfigure}%
  	\caption{
  	(Numerical) solutions of the linear advection--diffusion equation \cref{eq:periodc_advDif} for $a=1$ and $\varepsilon = 10^{-2}$ with periodic initial and boundary data at time $t=0.1$ and mass and energy over time. 
	We used a multi-block FSBP-SAT scheme with a polynomial and RBF approximation space, $\mathcal{P}_{2} = \Span\{1,x,x^2\}$ and $\mathcal{R}_{3} = \Span\{1,x,e^{x^2}\}$, respectively. 
  	}
  	\label{fig:linear_AdvDiff_kernel}
\end{figure} 

\begin{figure}[tb]
	\centering 
	\begin{subfigure}[b]{0.33\textwidth}
		\includegraphics[width=\textwidth]{%
      		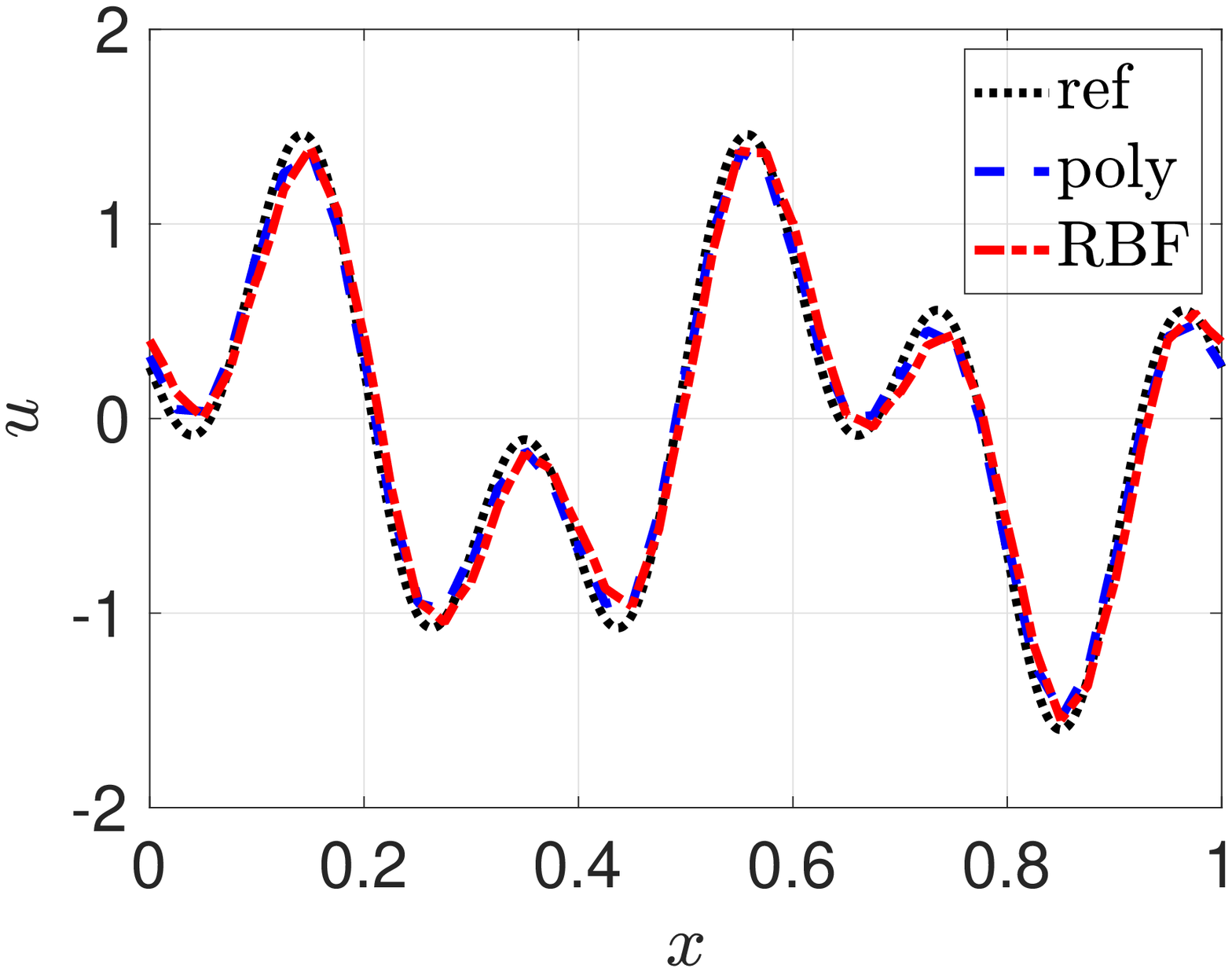} 
    		\caption{Solutions for $\alpha = 1/2$}
    		\label{fig:linear_AdvDiff_kernel_sol_I10_shapeParam05}
  	\end{subfigure}%
	\begin{subfigure}[b]{0.33\textwidth}
		\includegraphics[width=\textwidth]{%
      		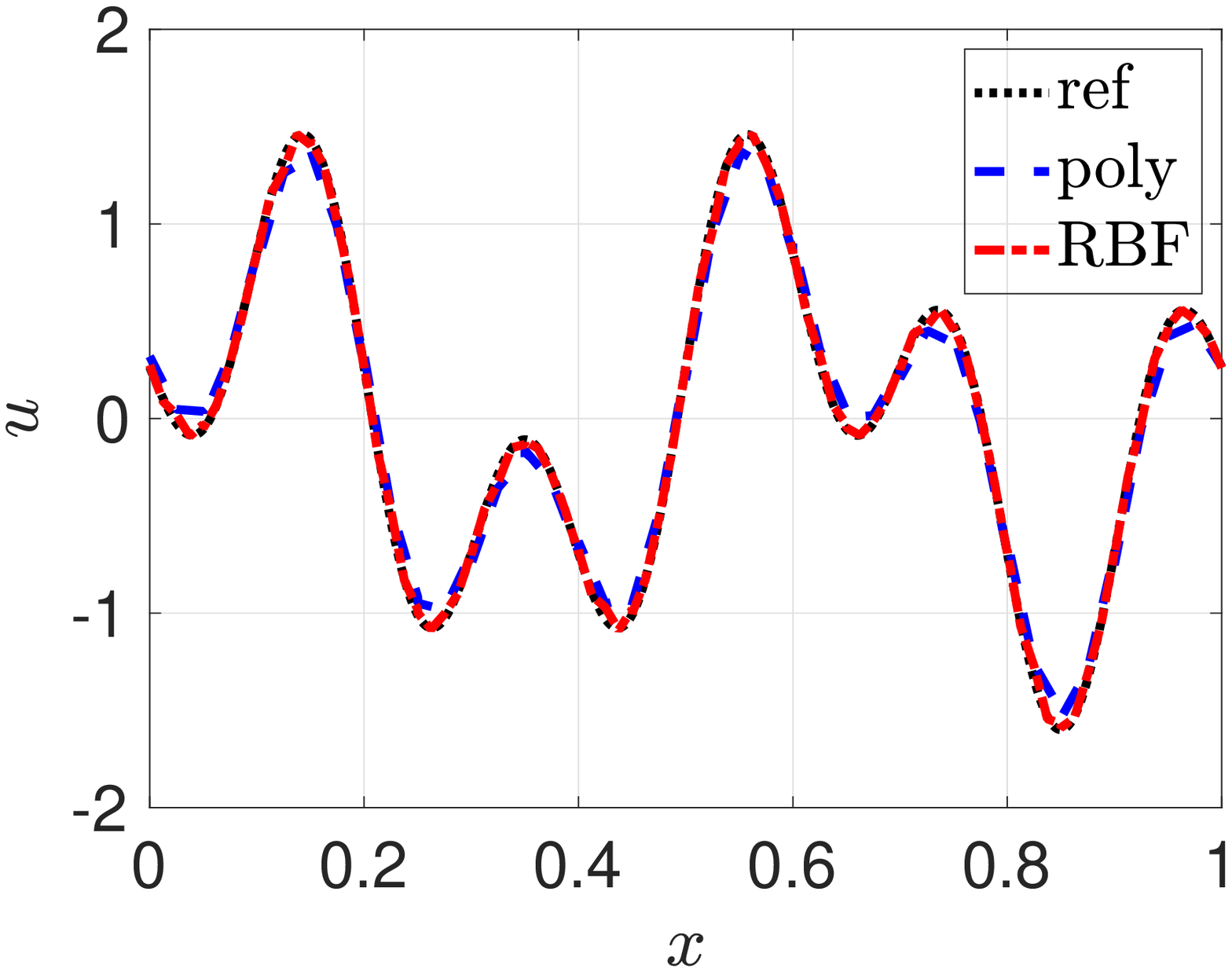} 
    		\caption{Solutions for $\alpha = 1$}
    		\label{fig:linear_AdvDiff_kernel_sol_I10_shapeParam1}
  	\end{subfigure}%
	\begin{subfigure}[b]{0.33\textwidth}
		\includegraphics[width=\textwidth]{%
      		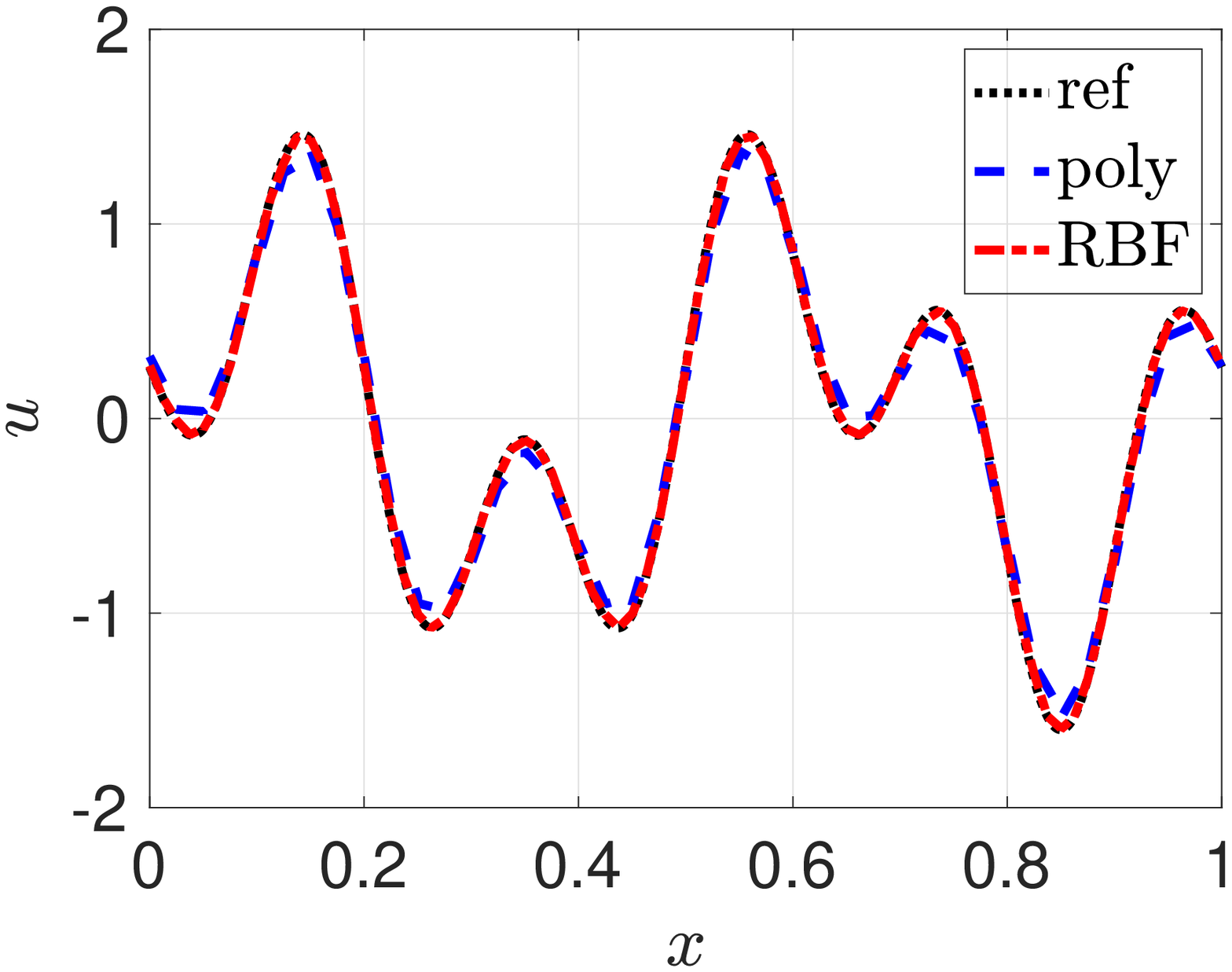} 
    		\caption{Solutions for $\alpha = 2$}
    		\label{fig:linear_AdvDiff_kernel_sol_I10_shapeParam2}
  	\end{subfigure}%
	\\
  	\begin{subfigure}[b]{0.33\textwidth}
		\includegraphics[width=\textwidth]{%
      		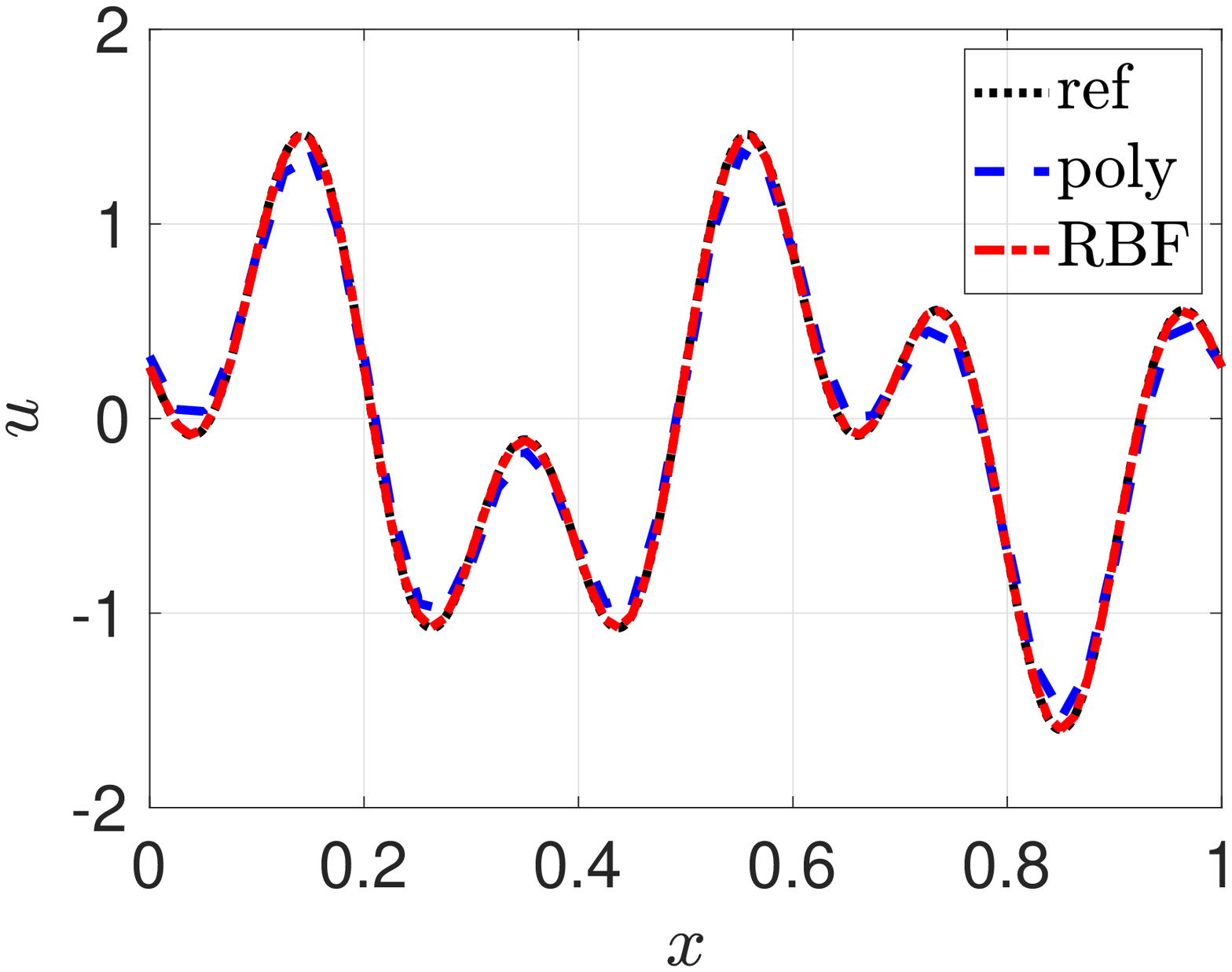} 
    		\caption{Solutions for $\alpha = 4$}
    		\label{fig:linear_AdvDiff_kernel_sol_I10_shapeParam4}
  	\end{subfigure}%
	\begin{subfigure}[b]{0.33\textwidth}
		\includegraphics[width=\textwidth]{%
      		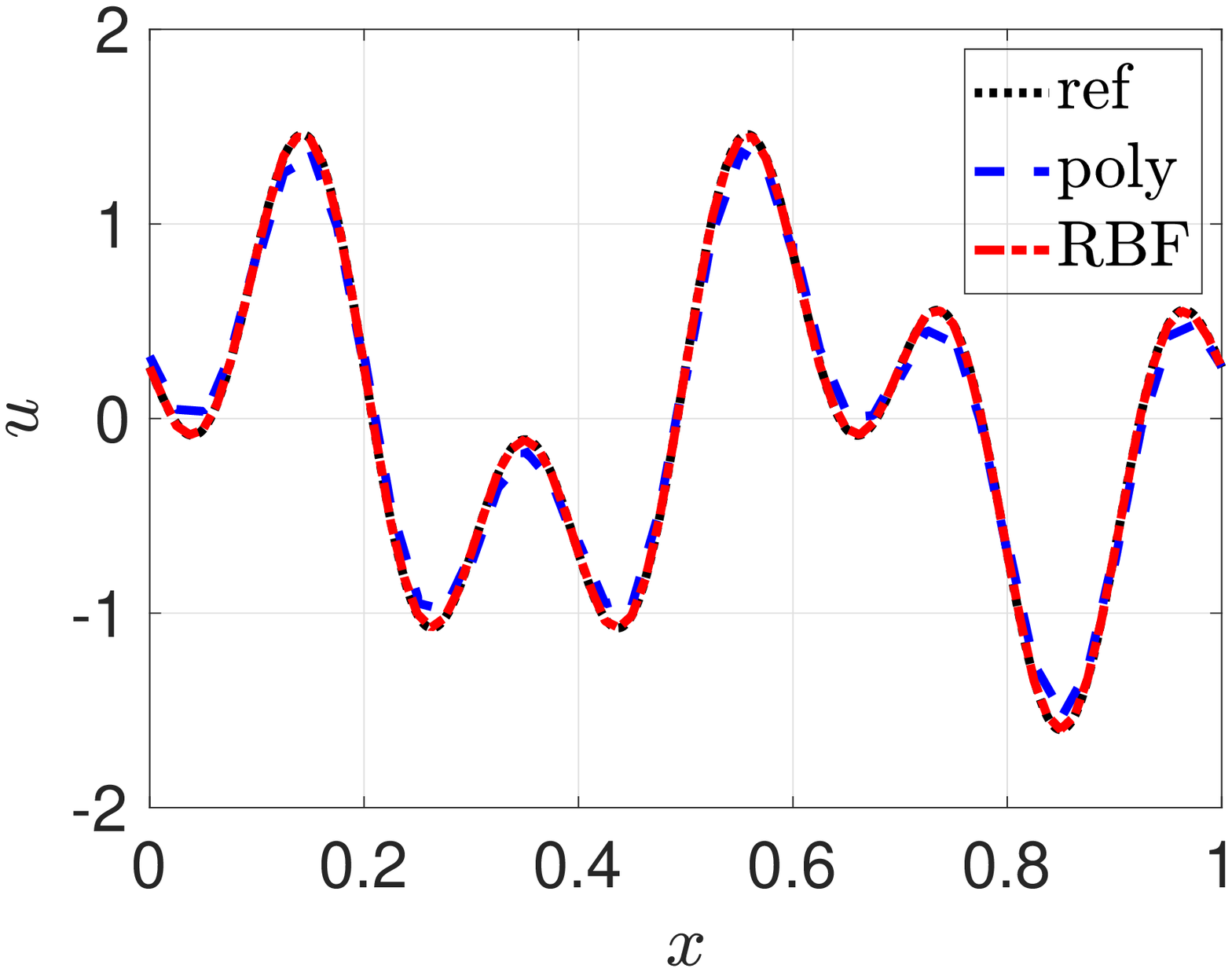} 
    		\caption{Solutions for $\alpha = 8$}
    		\label{fig:linear_AdvDiff_kernel_sol_I10_shapeParam8}
  	\end{subfigure}%
  	\begin{subfigure}[b]{0.33\textwidth}
		\includegraphics[width=\textwidth]{%
      		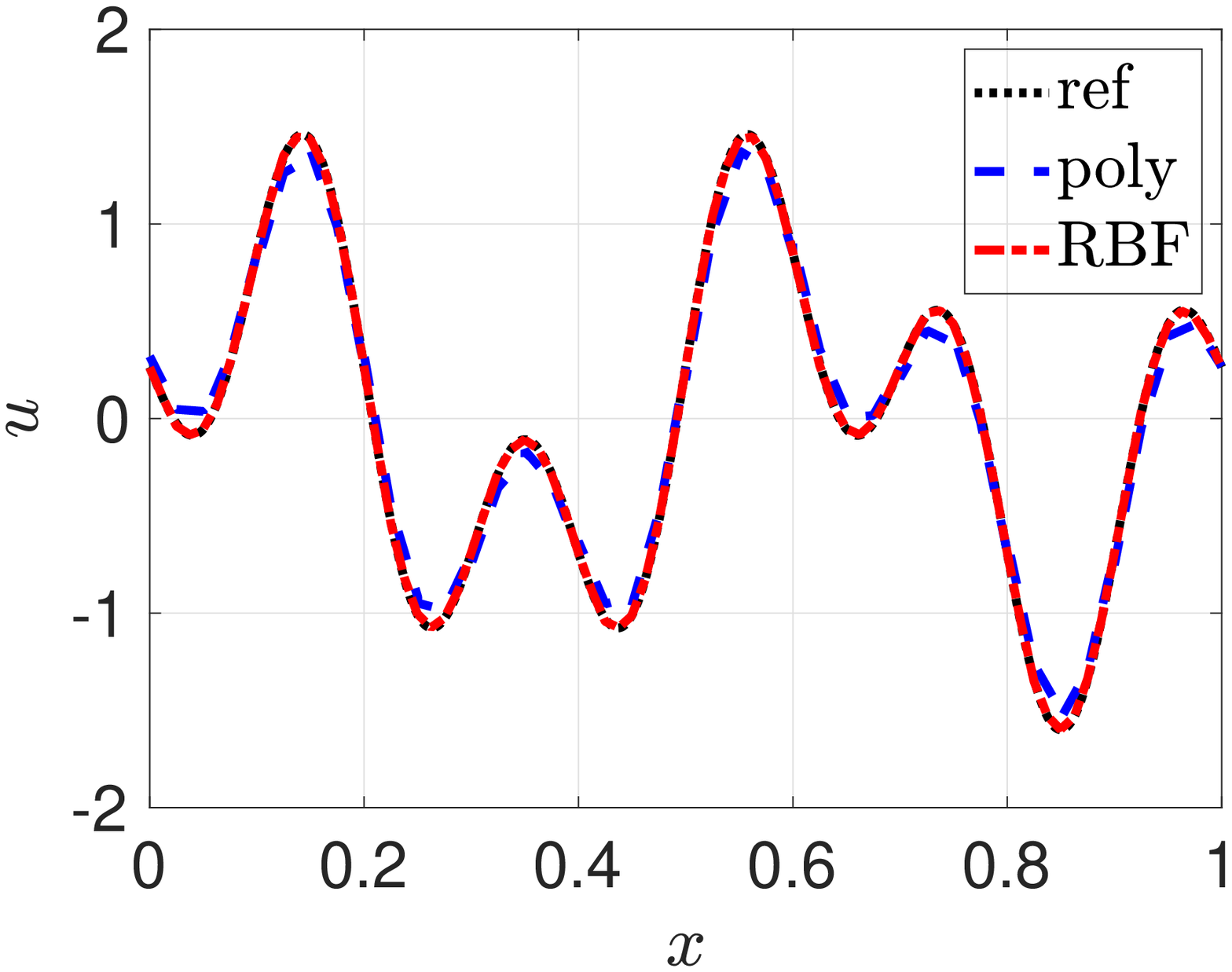} 
    		\caption{Solutions for $\alpha = 16$}
    		\label{fig:linear_AdvDiff_kernel_sol_I10_shapeParam16}
  	\end{subfigure}%
  	\caption{
  	(Numerical) solutions of the linear advection--diffusion equation \cref{eq:periodc_advDif} for $a=1$ and $\varepsilon = 10^{-2}$ with periodic initial and boundary data at time $t=0.1$ and mass and energy over time. 
	We used a multi-block FSBP-SAT scheme with $I=10$ blocks for the polynomial approximation space $\mathcal{P}_{2} = \Span\{1,x,x^2\}$ and the RBF approximation space, $\mathcal{R}_{3} = \Span\{1,x,e^{(x/\alpha)^2}\}$ for different shape parameters $\alpha=0.5,1,2,4,8,16$. 
	}
  	\label{fig:linear_AdvDiff_kernel_shapeParam}
\end{figure}

\cref{fig:linear_AdvDiff_kernel} provides a comparison between the exact solution and numerical solutions of the linear advection-diffusion problem \cref{eq:periodc_advDif} for $a=1$ and $\varepsilon = 10^{-2}$ at $t=0.1$. 
This time, we compare the polynomial SBP-SAT (``poly") with an RBF SBP-SAT method (``RBF"). 
The corresponding approximation spaces, $\mathcal{P}_{2}$ and $\mathcal{R}_{3}$, are both three-dimensional. 
The resulting first- and second-derivative FSBP operators were discussed in \cref{sub:examples_poly,sub:examples_kernel}, respectively. 
\cref{fig:linear_AdvDiff_kernel} demonstrates that the RBF approximation space better represents the solution than the polynomial one. 
Also, while the mass profiles of the RBF and polynomial scheme are both accurate up to machine precision, the RBF one yields a more accurate energy profile. 
In our experiments, we also tested RBF approximation spaces $\mathcal{R}_{3} = \Span\{1,x,e^{(x/\alpha)^2}\}$ with different shape parameters, albeit without extensive parameter optimization. 
We found the FSBP-SAT method to perform robustly for a range of shape parameters $\alpha$. 
\cref{fig:linear_AdvDiff_kernel_shapeParam} illustrates this by reporting on the numerical solutions for the same test problem as before with different shape parameters $\alpha=0.5,1,2,4,8,16$. 
Importantly, we do not assert that any particular approximation space is optimal for this or any other problem. 
Our primary aim is to validate the effectiveness of the proposed second-derivative FSBP operators and their adaptability to a range of non-polynomial approximation spaces. 

\subsection{The two-dimensional linear advection--diffusion equation} 
\label{sub:2D}

It is straightforward to extend the proposed first- and second-derivative FSBP operators to multiple dimensions using a tensor-product strategy \cite[Chapter 7.1.6]{trangenstein2009numerical}. 
To demonstrate this, consider the two-dimensional linear advection-diffusion equation 
\begin{equation}\label{eq:advDif_2D} 
	\partial_t u + a_1 \partial_x u + a_2 \partial_y u = \varepsilon_1 \partial_{xx} u + \varepsilon_2 \partial_{yy} u
\end{equation} 
on $\Omega = [0,1]^2$ with periodic boundary conditions and initial data 
\begin{equation} 
	u(x,y,0) = \exp\left( -200 \left[ (x-1/4)^2 + (y-1/4)^2 \right] \right). 
\end{equation} 

\begin{figure}[tb]
	\centering 
  	\begin{subfigure}[b]{0.49\textwidth}
		\includegraphics[width=\textwidth]{%
      		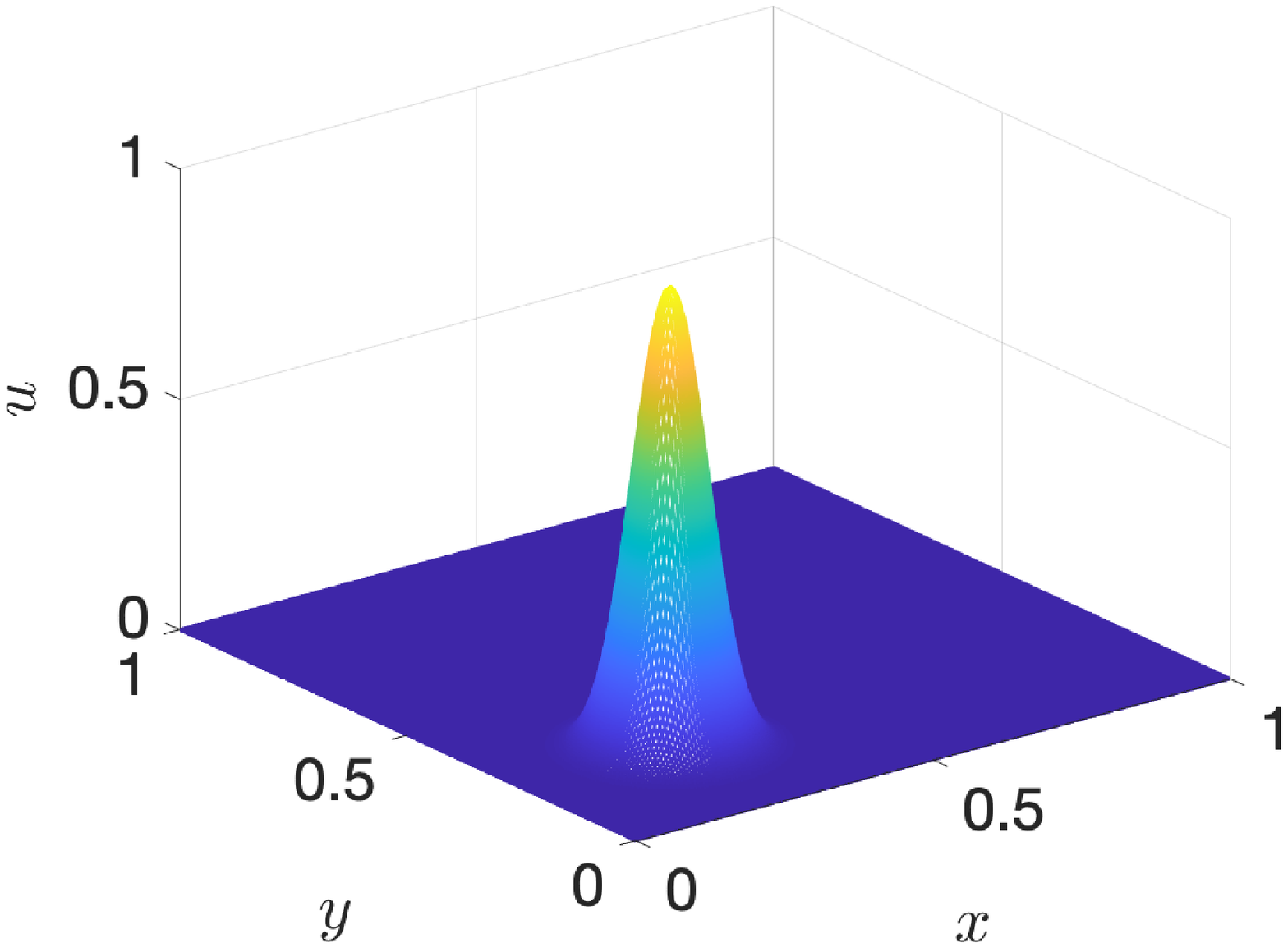} 
    		\caption{Initial data}
    		\label{fig:linear_AdvDiff_2D_IC}
  	\end{subfigure}%
	\begin{subfigure}[b]{0.49\textwidth}
		\includegraphics[width=\textwidth]{%
      		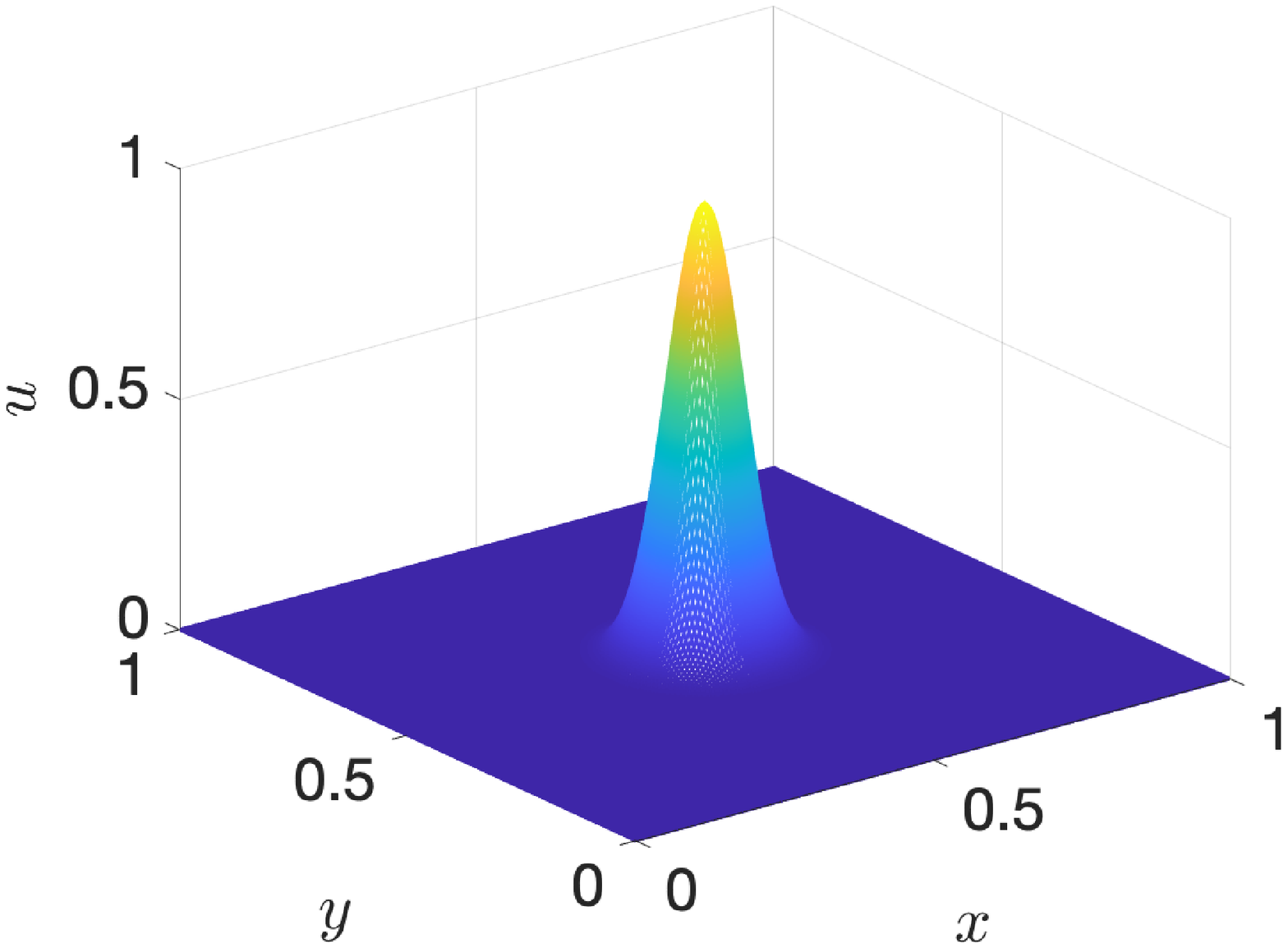} 
    		\caption{Reference solution}
    		\label{fig:linear_AdvDiff_2D_ref_I100}
  	\end{subfigure}%
	\\
	\begin{subfigure}[b]{0.49\textwidth}
		\includegraphics[width=\textwidth]{%
      		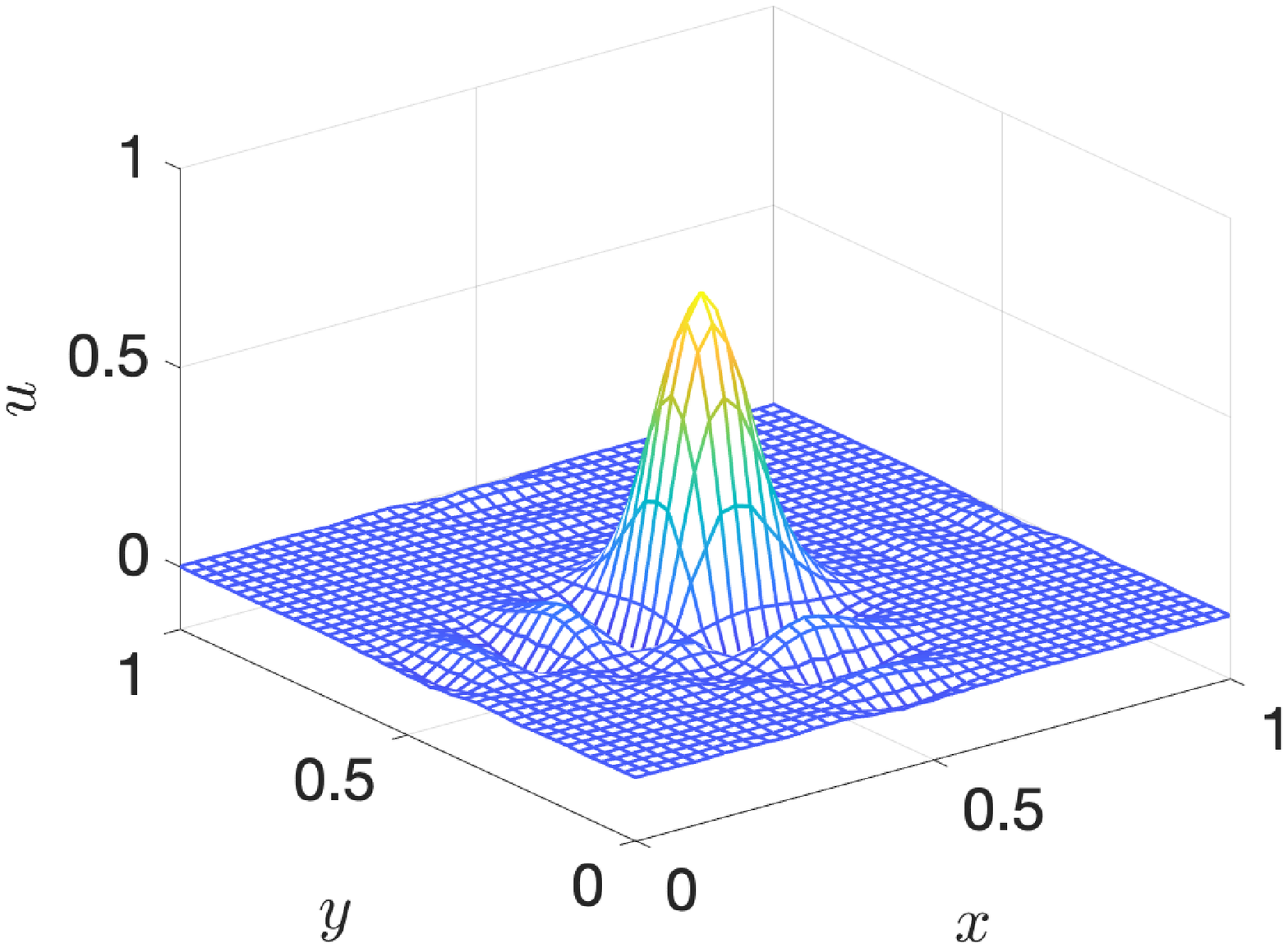} 
    		\caption{Polynomial solution}
    		\label{fig:linear_AdvDiff_2D_polyGL_I20}
  	\end{subfigure}%
  	\begin{subfigure}[b]{0.49\textwidth}
		\includegraphics[width=\textwidth]{%
      		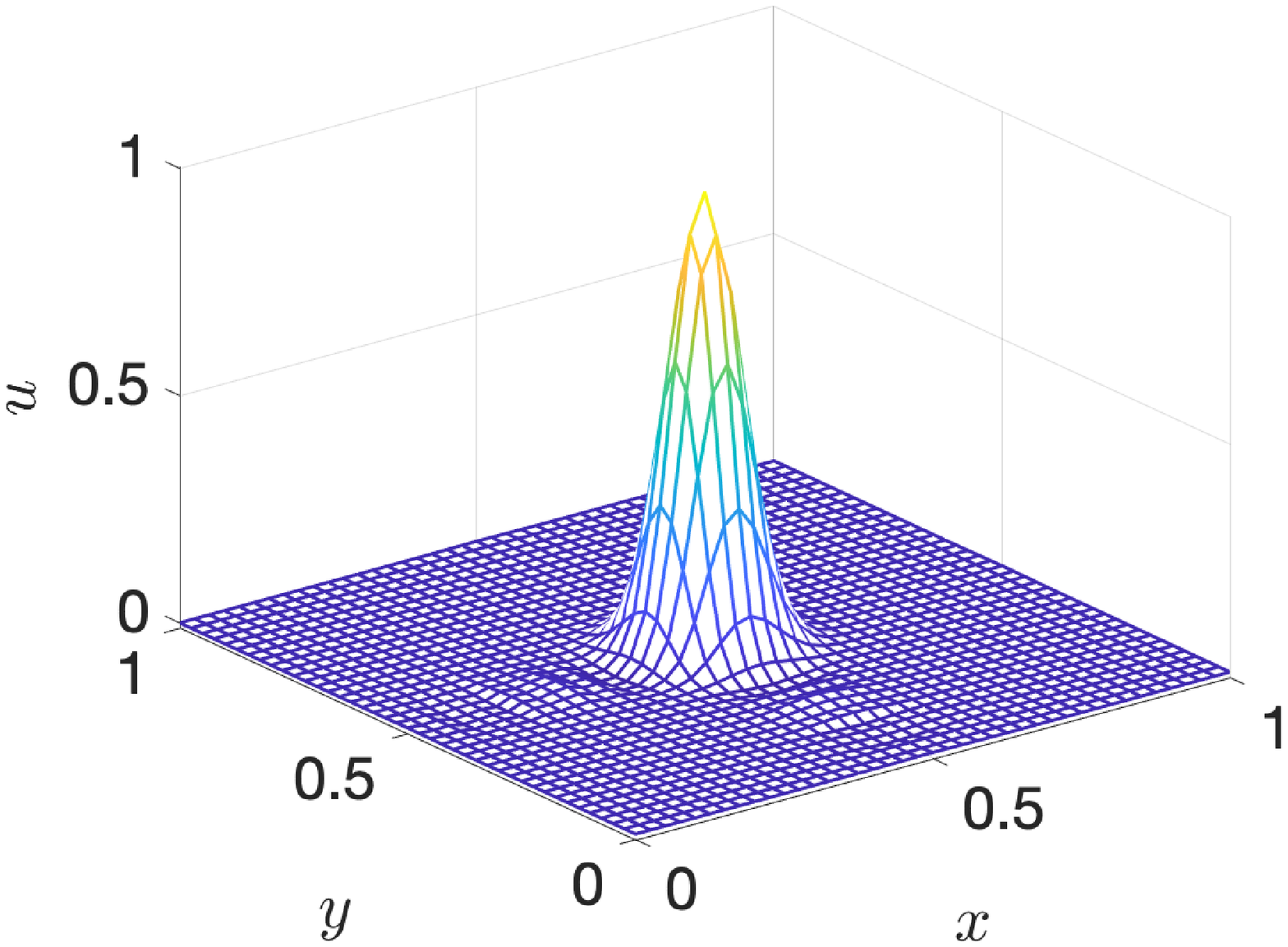} 
    		\caption{RBF solution}
    		\label{fig:linear_AdvDiff_2D_RBF_I20}
  	\end{subfigure}%
  	\caption{ 
  	(Numerical) solutions of the two-dimensional linear advection--diffusion equation \cref{eq:advDif_2D} for $a_1=a_2=1$ and $\varepsilon_1 = \varepsilon_2 = 10^{-4}$ at time $t=1/4$. 
	We used the tensor-product strategy combined with a multi-block FSBP-SAT scheme and $I=20$ blocks in each direction.  
	We considered a polynomial and RBF approximation space, $\mathcal{P}_{2} = \Span\{1,x,x^2\}$ and $\mathcal{R}_{3} = \Span\{1,x,e^{20 x^2}\}$, respectively. 
  	}
  	\label{fig:linear_AdvDiff_2D}
\end{figure} 

\cref{fig:linear_AdvDiff_2D} shows the initial data and the reference and numerical solutions of \cref{eq:advDif_2D} using a polynomial and an RBF approximation space, $\mathcal{P}_{2} = \Span\{1,x,x^2\}$ and $\mathcal{R}_{3} = \Span\{1,x,e^{20 x^2}\}$, at time $t=1/4$ for $a_1 = a_2 = 1$ and $\varepsilon_1 = \varepsilon_2 = 10^{-4}$. 
For both approximation spaces, we used the same multi-block FSBP-SAT semi-discretization as in \cref{sub:kernel} combined with the tensor-product strategy and $I=20$ blocks in each direction. 
The reference solution was computed using the polynomial operator with $I=100$ blocks in each direction. 
We observe that the RBF FSBP operator yields a more accurate numerical solution than the polynomial one. 
We quantify this observations in \cref{tab:errors}, which reports the relative errors of the numerical solutions of \cref{eq:advDif_2D} using the polynomial and the RBF function space w.r.t.\ to the $\| \cdot \|_1$-, $\| \cdot \|_2$-, and $\| \cdot \|_\infty$-norm. 
\cref{tab:errors} shows that the RBF functions space yields relative errors around one magnitude smaller than the ones for the polynomial function space.  

\begin{table}[tb]
\renewcommand{\arraystretch}{1.3}
	\centering 
	%	\begin{adjustbox}{width=0.99\textwidth} 
  	\begin{tabular}{c c c c c c c} 
		\toprule 
    		Function space & & $\| \cdot \|_1$-norm & & $\| \cdot \|_2$-norm & & $\| \cdot \|_\infty$-norm \\ \hline 
		polynomial	&  & $9.2 \cdot 10^{-1}$ & & $4.3 \cdot 10^{-1}$ & & $3.4 \cdot 10^{-1}$ \\
    		RBF 				& & $1.0 \cdot 10^{-1}$ & & $5.8 \cdot 10^{-2}$ & & $5.9 \cdot 10^{-2}$ \\  
		\bottomrule
	\end{tabular} 
\caption{
Relative errors w.r.t. to different norms of the numerical solutions of the two-dimensional linear advection-diffusion equation \cref{eq:advDif_2D} for $a_1=a_2=1$ and $\varepsilon_1 = \varepsilon_2 = 10^{-4}$ at time $t=1/4$.
We compare a polynomial and RBF approximation space, $\mathcal{P}_{2} = \Span\{1,x,x^2\}$ and $\mathcal{R}_{3} = \Span\{1,x,e^{20 x^2}\}$, respectively. 
}
\label{tab:errors}
\end{table} 

\begin{remark} 
	The tensor-product strategy is often employed because it is simple and efficient. 
	Still, it is not always the best choice since it limits geometric 
flexibility. 
	To obtain first-derivative FSBP operators (and, consequently, second-derivative FSBP operators) on general multi-dimensional geometries (e.g., triangles and circles), see \cite{glaubitz2023multi}.
\end{remark}

\subsection{A boundary layer problem} 
\label{sub:bLayer}

Consider the boundary layer problem  
\begin{equation}\label{eq:boundLayer} 
\begin{aligned}
	\partial_t u + \partial_x u & = \varepsilon \partial_{xx} u, \quad && 0 < x < 1/2, \, t>0, \\ 
	u(x,0) & = 2x, \quad && 0 \leq x \leq 1/2, \\ 
	u(0,t) & = 0, \quad && t \geq 0, \\ 
	u(1/2,t) & = 1, \quad && t \geq 0.
\end{aligned}
\end{equation} 

The exact steady-state solution of \cref{eq:boundLayer} is $u(x) = ( e^{x/\varepsilon} - 1 ) / ( e^{1/2\varepsilon} - 1 )$, which has a sharp boundary layer at the right boundary for small $\varepsilon$. 
Notably, $u$ manifests significantly steeper gradients near the right domain boundary compared to the rest of the domain. 
This characteristic suggests that an approximation space incorporating an exponential function of the form $e^{\alpha x}$ with $\alpha > 0$ would appropriately capture the solution's behavior. 
Accordingly, we employed the exponential function space $\mathcal{E}_2 = \Span\{1, x, e^{\alpha x}\}$ for this problem.
It's worth mentioning that if the boundary layer were located at the left domain boundary, one would instead consider an exponential function of the form $e^{\alpha x}$ with $\alpha < 0$, which accentuates gradients near the left boundary. 
In our experiments, we tested a few different values for $\alpha$ and found that $\alpha = 1/10$ yielded satisfactory results, albeit without extensive parameter optimization. 
Notably, we do not claim that the chosen approximation space is optimal for this or any other specific problem. 
Instead, our main objective is to demonstrate the efficacy of the newly introduced second-derivative FSBP operators and their flexibility across various non-polynomial approximation spaces. 

\begin{figure}[tb]
	\centering 
	\begin{subfigure}[b]{0.49\textwidth}
		\includegraphics[width=\textwidth]{%
      		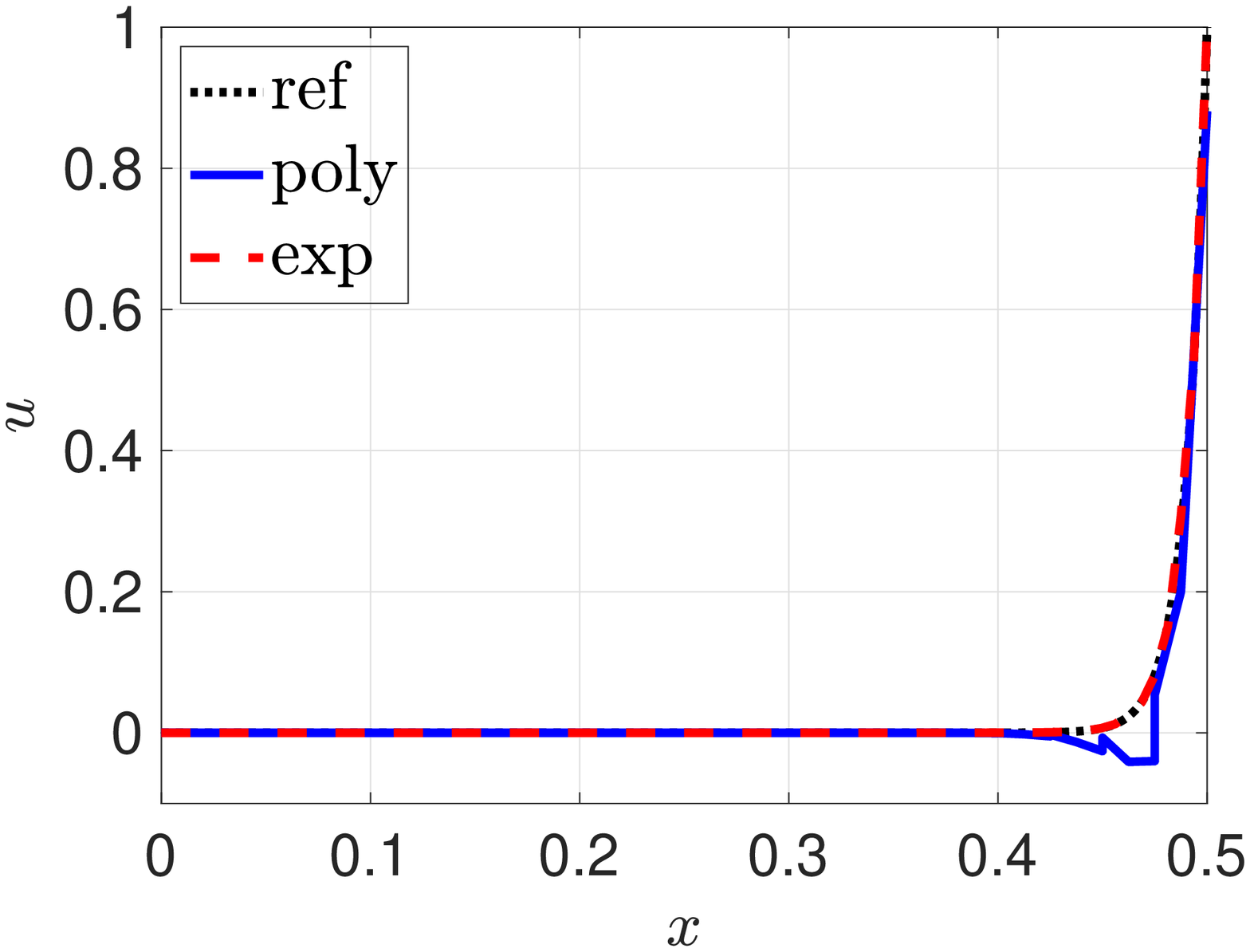} 
    		\caption{Solutions for $I=20$ blocks}
    		\label{fig:boundaryLayer_sol_I20_GL}
  	\end{subfigure}%
	~
	\begin{subfigure}[b]{0.49\textwidth}
		\includegraphics[width=\textwidth]{%
      		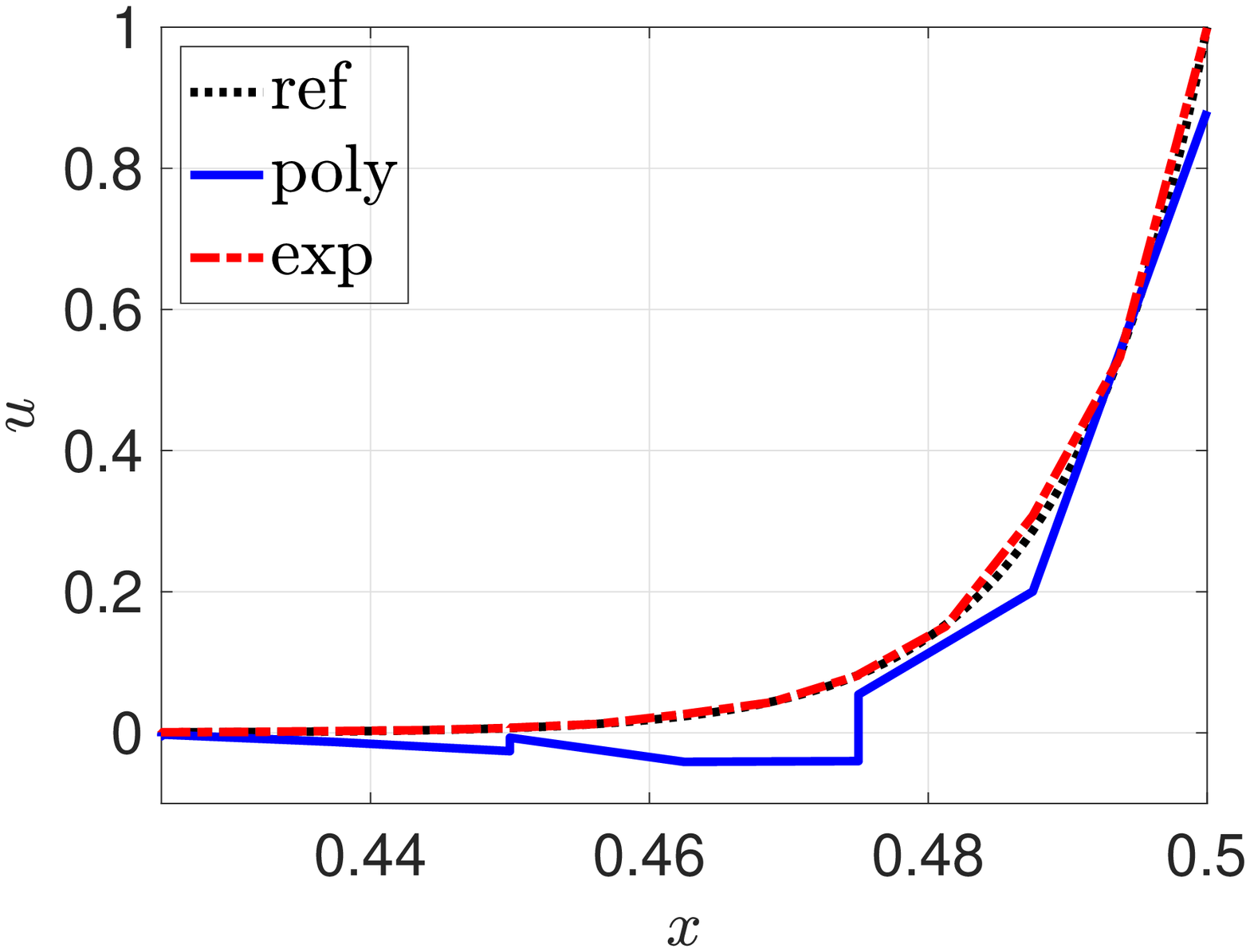} 
    		\caption{Zoomed-in solutions for $I=20$ blocks}
    		\label{fig:boundaryLayer_sol_zoom_I20_GL}
  	\end{subfigure}%
	\\ 
	\begin{subfigure}[b]{0.49\textwidth}
		\includegraphics[width=\textwidth]{%
      		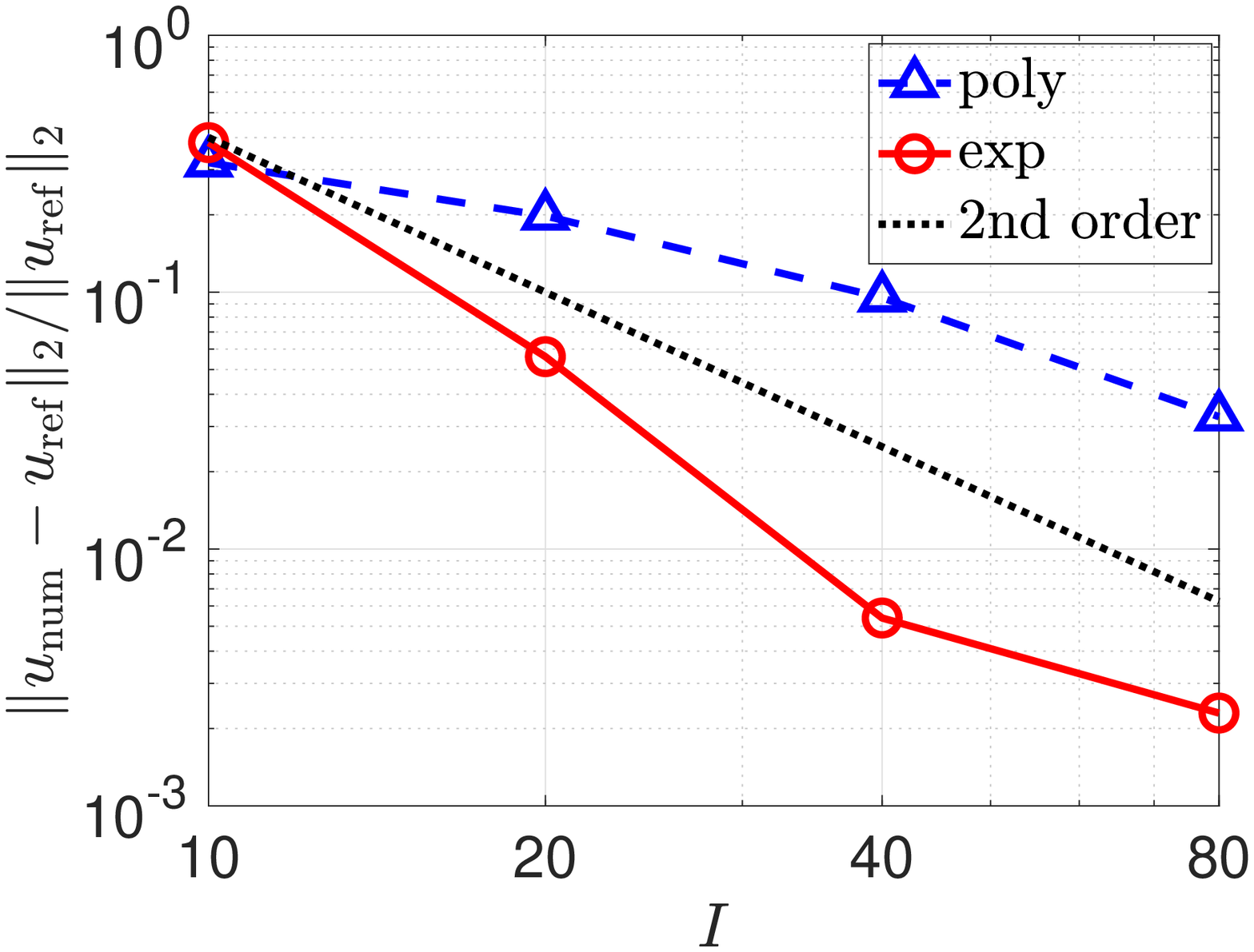} 
    		\caption{Relative $\|\cdot\|_2$-norm error}
    		\label{fig:boundaryLayer_L2_error_GL}
  	\end{subfigure}%
	~
  	\begin{subfigure}[b]{0.49\textwidth}
		\includegraphics[width=\textwidth]{%
      		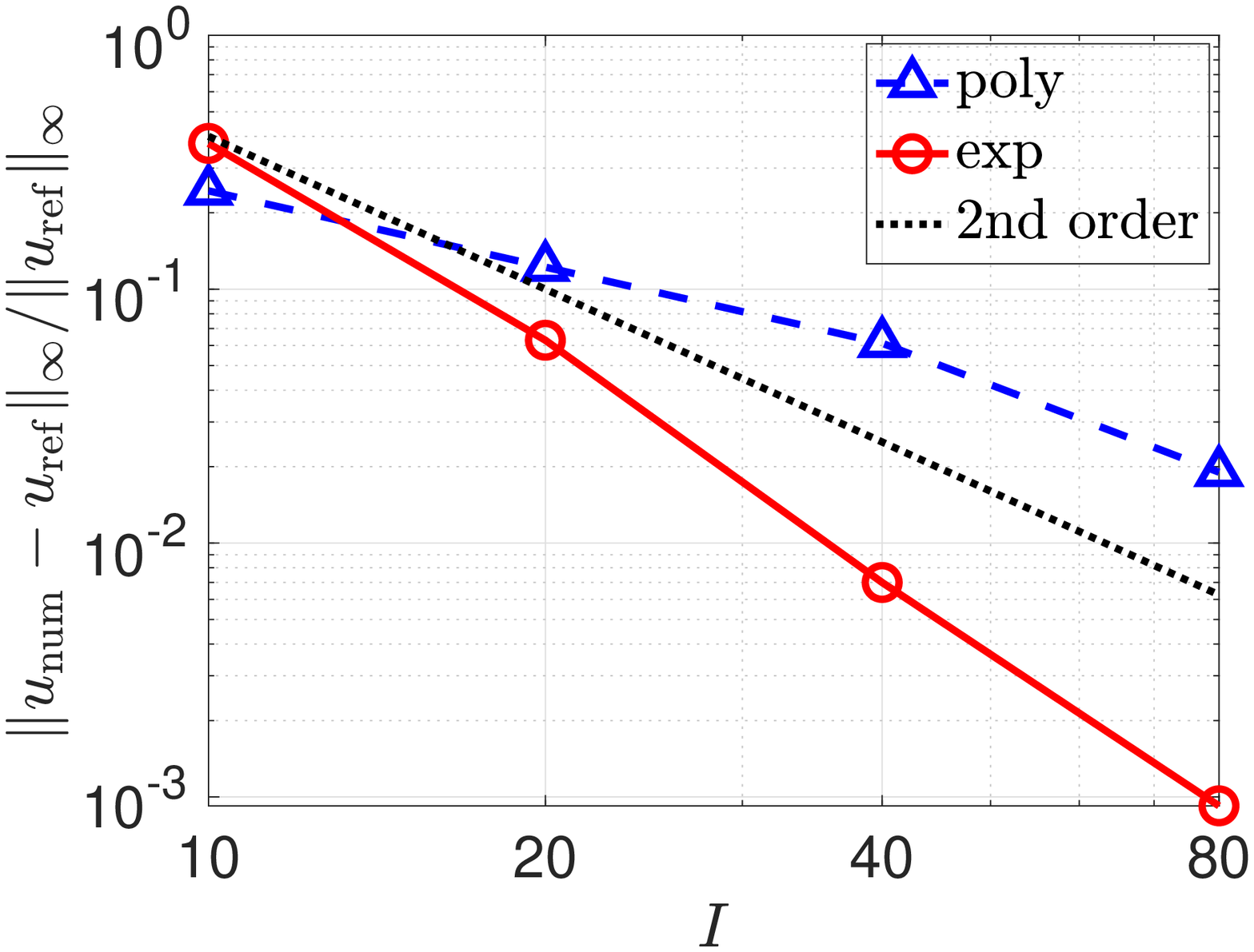} 
    		\caption{Relative $\|\cdot\|_{\infty}$-norm error}
    		\label{fig:boundaryLayer_max_error_GL}
  	\end{subfigure}%
  	\caption{ 
	(Numerical) solutions and the corresponding relative errors for the boundary layer problem \cref{eq:boundLayer} for $\varepsilon = 10^{-2}$ at $t=0.75$. 
	We used a multi-block FSBP-SAT scheme with a polynomial and exponential approximation space, $\mathcal{P}_2 = \Span\{1,x,x^2\}$ and $\mathcal{E}_2 = \Span\{1,x,e^{x/10}\}$, respectively. 
	}
  	\label{fig:boundaryLayer}
\end{figure}

\cref{fig:boundaryLayer} shows the (numerical) solutions and their relative errors for $\varepsilon = 10^{-2}$ at time $t=0.75$.  
The numerical solutions were computed using a multi-block FSBP-SAT method with a three-dimensional polynomial and exponential approximation space, $\mathcal{P}_2 = \Span\{1,x,x^2\}$ and $\mathcal{E}_2 = \Span\{1,x,e^{x/10}\}$, respectively. 
We used the FSBP-SAT semi-discretization \cref{eq:advDif_discr} with accordingly modified SATs \cref{eq:advDif_SATs}. 
We found that the exponential approximation space $\mathcal{E}_2$ yields more accurate results than the conventional polynomial approximation space $\mathcal{P}_2$. 
\cref{fig:boundaryLayer_L2_error_GL,fig:boundaryLayer_max_error_GL} provide the $\|\cdot\|_2$- and $\|\cdot\|_{\infty}$-norm errors of both FSBP-SAT methods for an increasing number of uniform blocks, $I$. 
While the convergence rate of the multi-block FSBP-SAT methods is roughly the same for both approximation spaces, the individual error levels are found to be smaller for the $\mathcal{E}_2$-based FSBP operators. 

\begin{figure}[tb]
	\centering 
	\begin{subfigure}[b]{0.49\textwidth}
		\includegraphics[width=\textwidth]{%
      		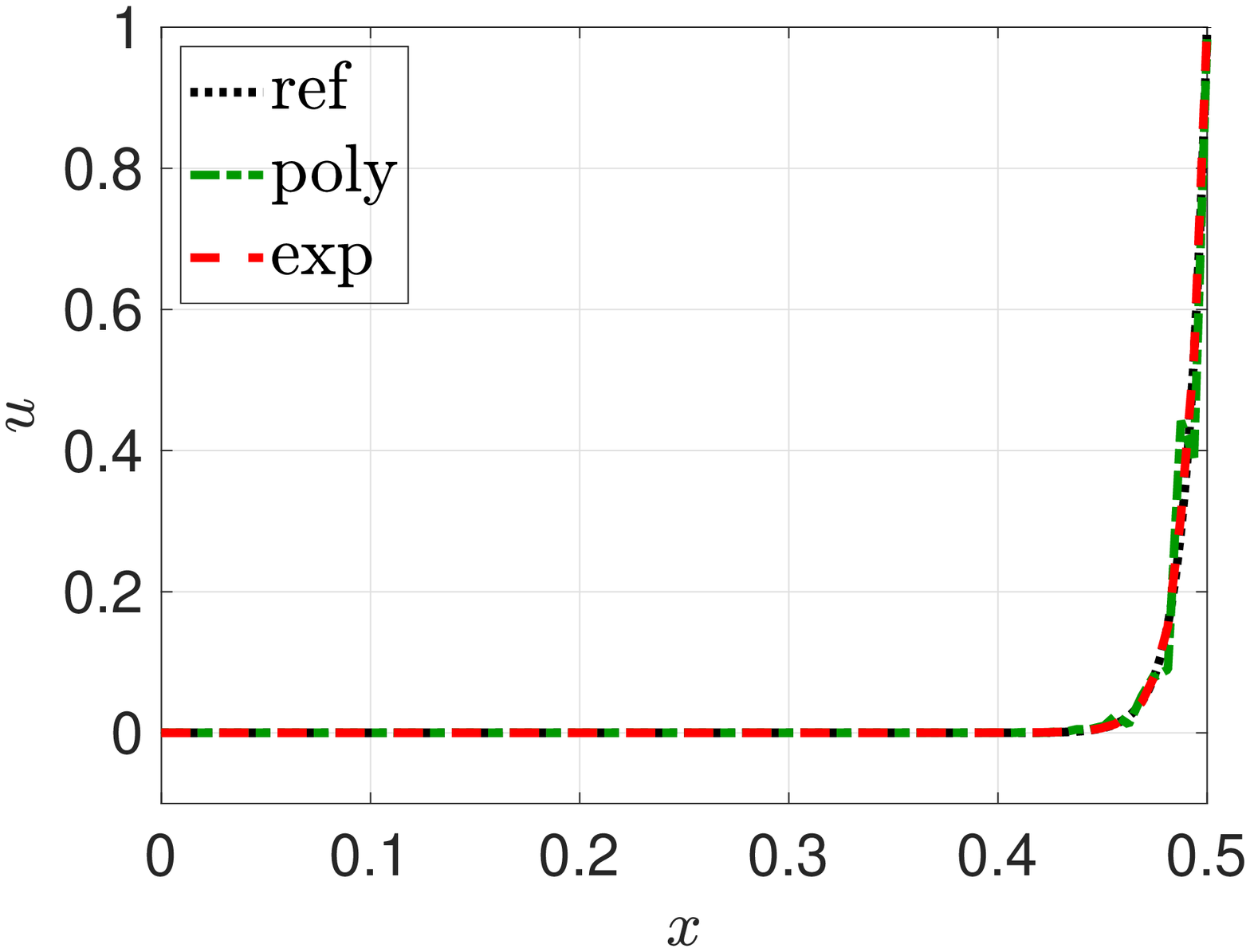} 
    		\caption{Solutions for $I=20$ blocks}
    		\label{fig:boundaryLayer_sol_I20_equi}
  	\end{subfigure}%
	~
	\begin{subfigure}[b]{0.49\textwidth}
		\includegraphics[width=\textwidth]{%
      		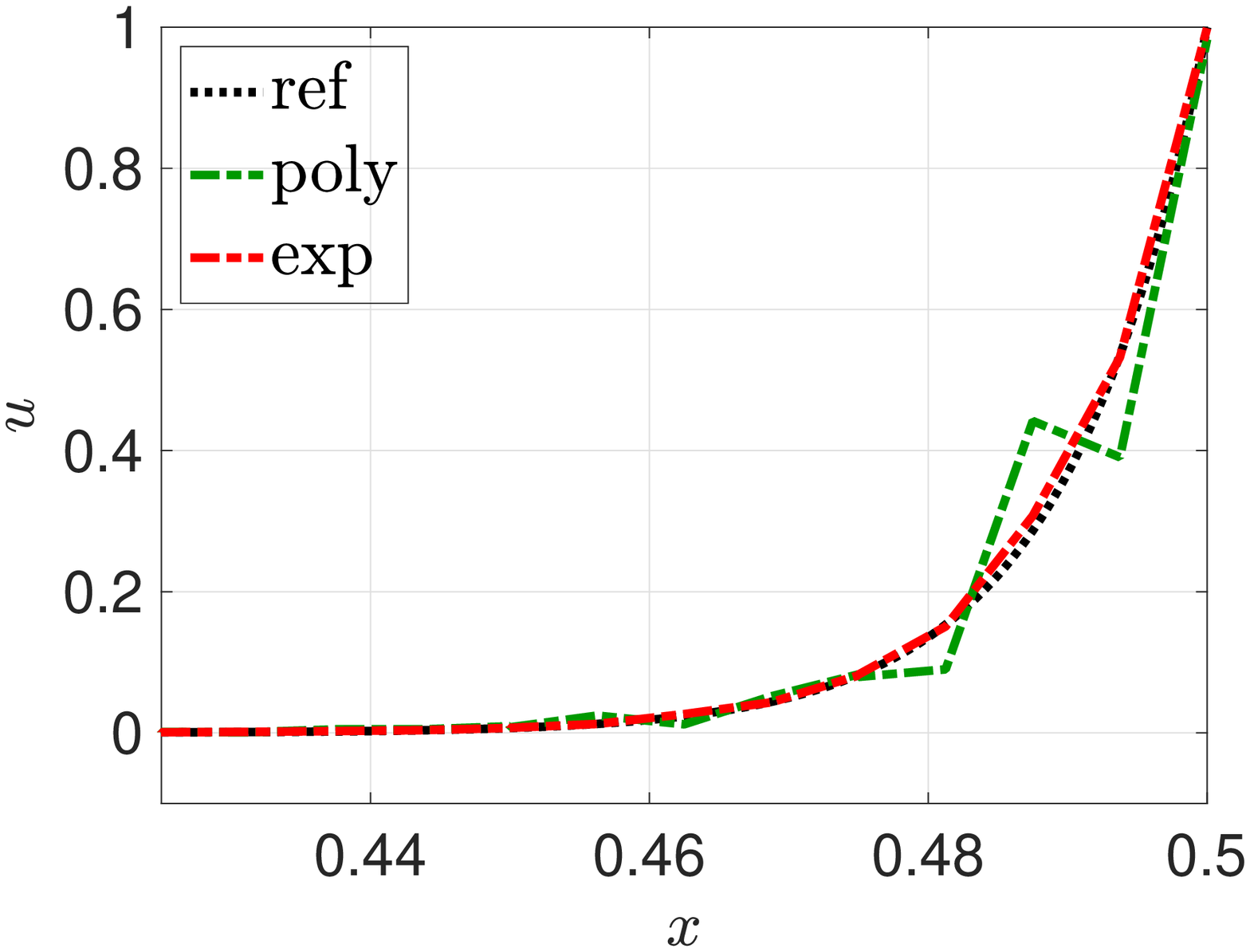} 
    		\caption{Zoomed-in solutions for $I=20$ blocks}
    		\label{fig:boundaryLayer_sol_zoom_I20_equi}
  	\end{subfigure}%
	\\
	\begin{subfigure}[b]{0.49\textwidth}
		\includegraphics[width=\textwidth]{%
      		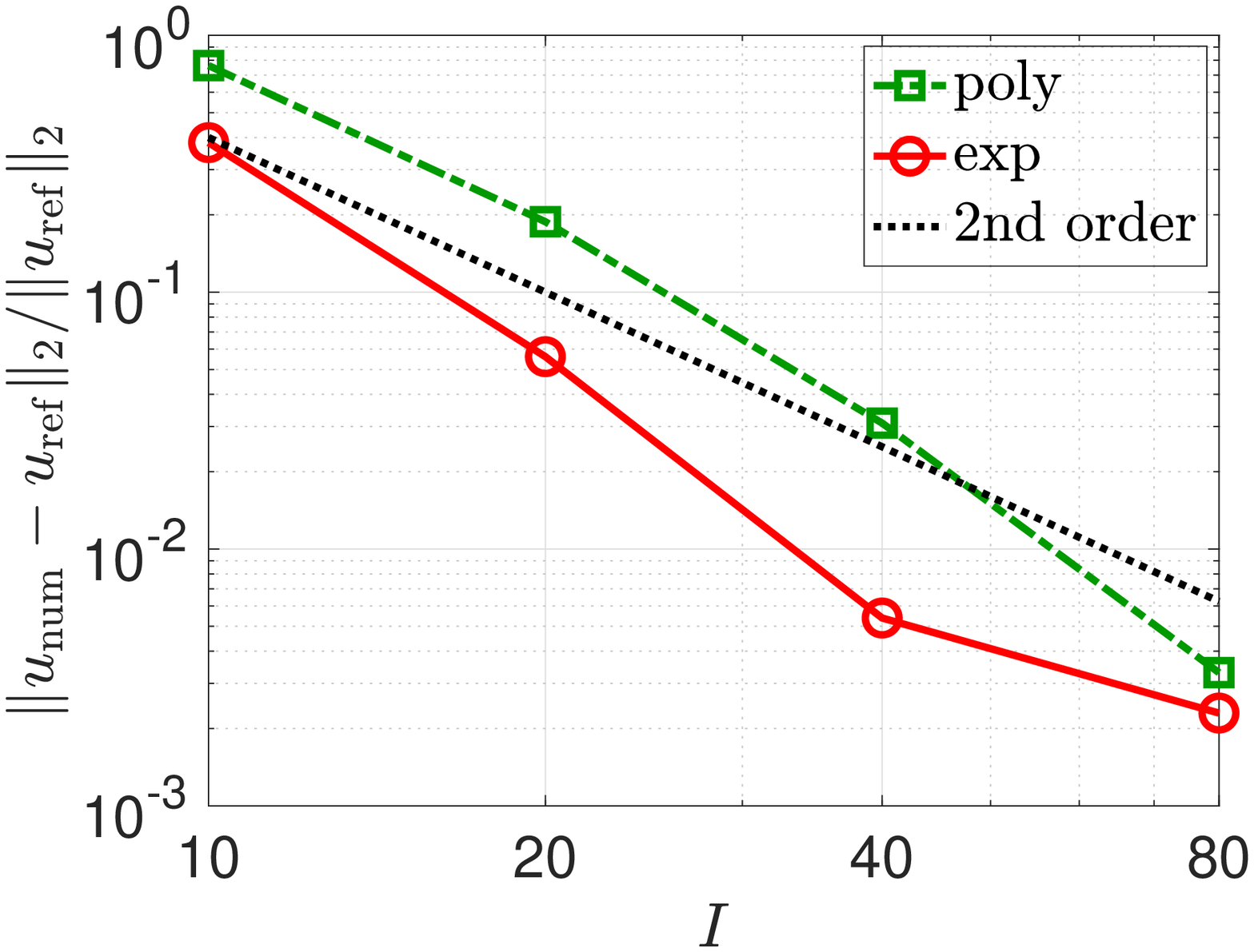} 
    		\caption{Relative $\|\cdot\|_2$-norm error}
    		\label{fig:boundaryLayer_L2_error_equi}
  	\end{subfigure}%
	~
  	\begin{subfigure}[b]{0.49\textwidth}
		\includegraphics[width=\textwidth]{%
      		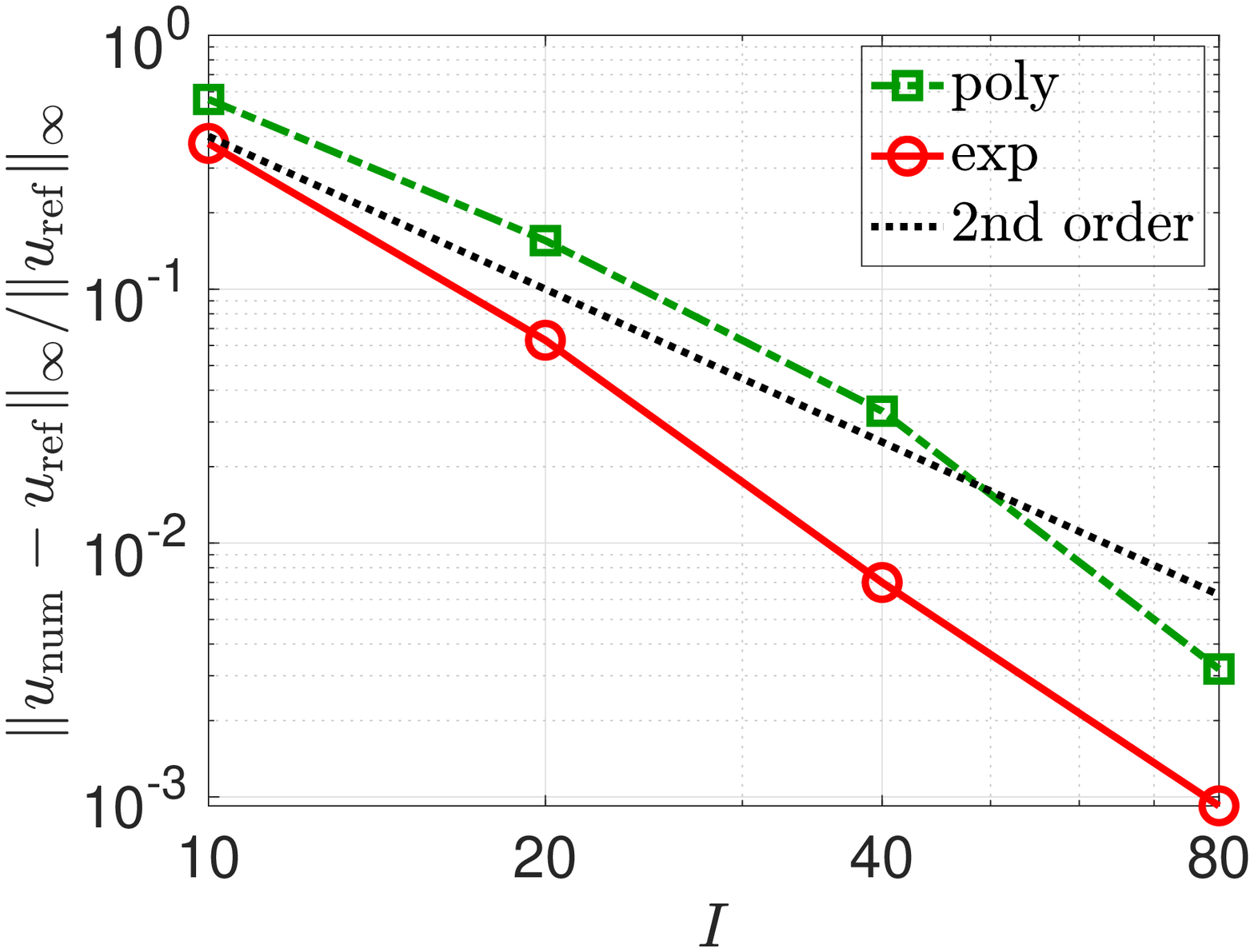} 
    		\caption{Relative $\|\cdot\|_{\infty}$-norm error}
    		\label{fig:boundaryLayer_max_error_equi}
  	\end{subfigure}%
  	\caption{ 
	(Numerical) solutions and the corresponding relative errors for the boundary layer problem \cref{eq:boundLayer} for $\varepsilon = 10^{-2}$ at $t=0.75$. 
	We used a multi-block FSBP-SAT scheme with a polynomial and exponential approximation space, $\mathcal{P}_2 = \Span\{1,x,x^2\}$ and $\mathcal{E}_2 = \Span\{1,x,e^{x/10}\}$, respectively. 
	\emph{Both schemes use the same grid points.} 
	}
  	\label{fig:boundaryLayer_equid}
\end{figure}

A pertinent query is whether the superior precision of the exponential FSBP operator over the polynomial SBP operator, as indicated in \cref{fig:boundaryLayer}, is attributable to their respective function spaces or the former employing two additional grid points per block. 
To elucidate this matter, we revisit the above test, this time utilizing a polynomial SBP operator that is exact for $\mathcal{P}_2 = \Span\{1,x,x^2\}$ and matches the grid points of the exponential FSBP operator.
\cref{fig:boundaryLayer_equid} shows the corresponding (numerical) solutions and their relative errors and demonstrates that the FSBP operator again yields more accurate results.

\begin{remark} 
	The Stone--Weierstrass approximation theorem \cite{stone1948generalized}---a generalization of the Weierstrass approximation theorem for polynomials---implies uniform $p$-convergence (increasing dimension of $\mathcal{F}$) to any continuous function for a large class of non-polynomial function spaces. 
	Additionally, the concept of $h$-convergence (reducing the size of all elements) in the context of non-polynomial function spaces has been explored in references \cite{yuan2006discontinuous, yang2016short, franck2023approximately}. 
	These studies have established $L^2$-error estimates for DG methods utilizing non-polynomial approximation spaces in solving time-dependent PDEs. 
	Specifically, \cite{yuan2006discontinuous} introduced a criterion enabling non-polynomial spaces to achieve approximation rates comparable to polynomial DG spaces of equivalent dimension. 
	Building on this, \cite{yang2016short} showed that a diverse range of approximation spaces, including well-known ones like trigonometric and exponential ones, fulfill the criterion set forth in \cite{yuan2006discontinuous}. 
	More recently, \cite{franck2023approximately} extended these findings to approximation spaces incorporating physics-informed neural network (PINN) approximations of steady-state solutions to obtain well-balanced DG methods for hyperbolic balance laws.
\end{remark}

\subsection{The viscous Burgers equation} 
\label{sub:Burgers} 

Consider the nonlinear viscous Burgers equation 
\begin{equation}\label{eq:Burgers} 
\begin{aligned} 
	\partial_t u + \partial_x \left( \frac{u^2}{2} \right) & = \varepsilon \partial_{xx} u, \quad && 0 < x < 1, \ t>0
\end{aligned}
\end{equation} 
with periodic boundary conditions and discontinuous 'sawtooth' initial data $u(x,0) = 1 + 2x$ if $x<0$ and $u(x,0) = 2x - 1$ otherwise. 
The exact solution has a steep gradient at $x=1/2$ for small $\varepsilon$, which we expect to be better approximated by exponential rather than polynomial approximation spaces.  
We consider an energy-stable skew-symmetric uniform multi-block FSBP-SAT semi-discretization of the form 
\begin{equation}\label{eq:Burgers_discr}
	\mathbf{u}_t^{(i)} 
		+ \frac{1}{3} D_1 U^{(i)} \mathbf{u}^{(i)} + \frac{1}{3} U^{(i)} D_1 \mathbf{u}^{(i)} = \varepsilon D_2 \mathbf{u}^{(i)} + P^{-1}\mathbb{S}^{(i)}, 
		\quad i=1,\dots,I,
\end{equation} 
where $U^{(i)} = \diag(\mathbf{u})^{(i)}$ and SATs $\mathbb{S}^{(i)} = \mathbb{S}^{(i)}_L + \mathbb{S}^{(i)}_R$ as in \cite{fisher2013discretely}.

\begin{figure}[tb]
	\centering 
	\begin{subfigure}[b]{0.33\textwidth}
		\includegraphics[width=\textwidth]{%
      		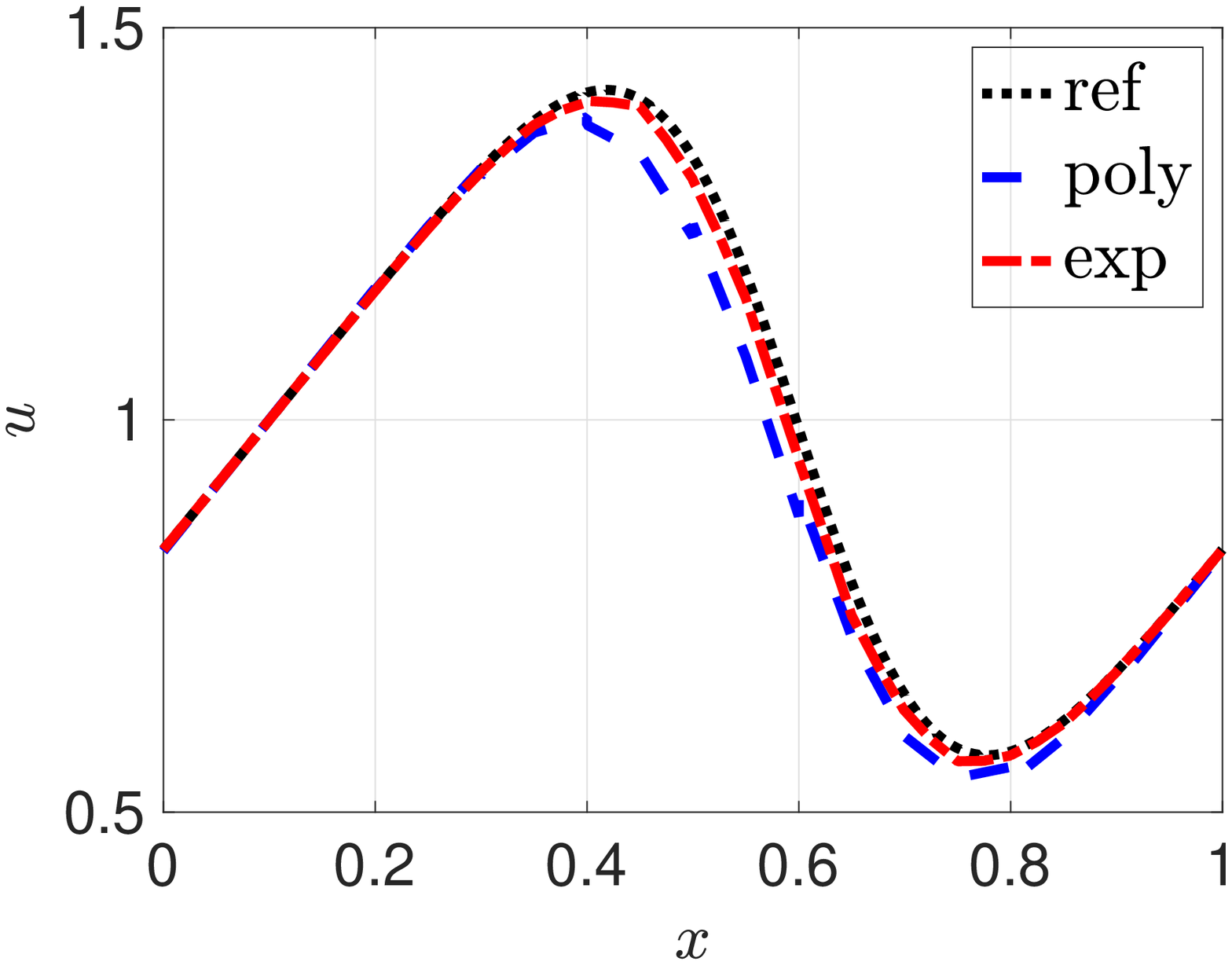} 
    		\caption{$\varepsilon = 10^{-1}$ \& $I=10$}
    		\label{fig:viscBurgers_eps_m10_I10}
  	\end{subfigure}%
	\begin{subfigure}[b]{0.33\textwidth}
		\includegraphics[width=\textwidth]{%
      		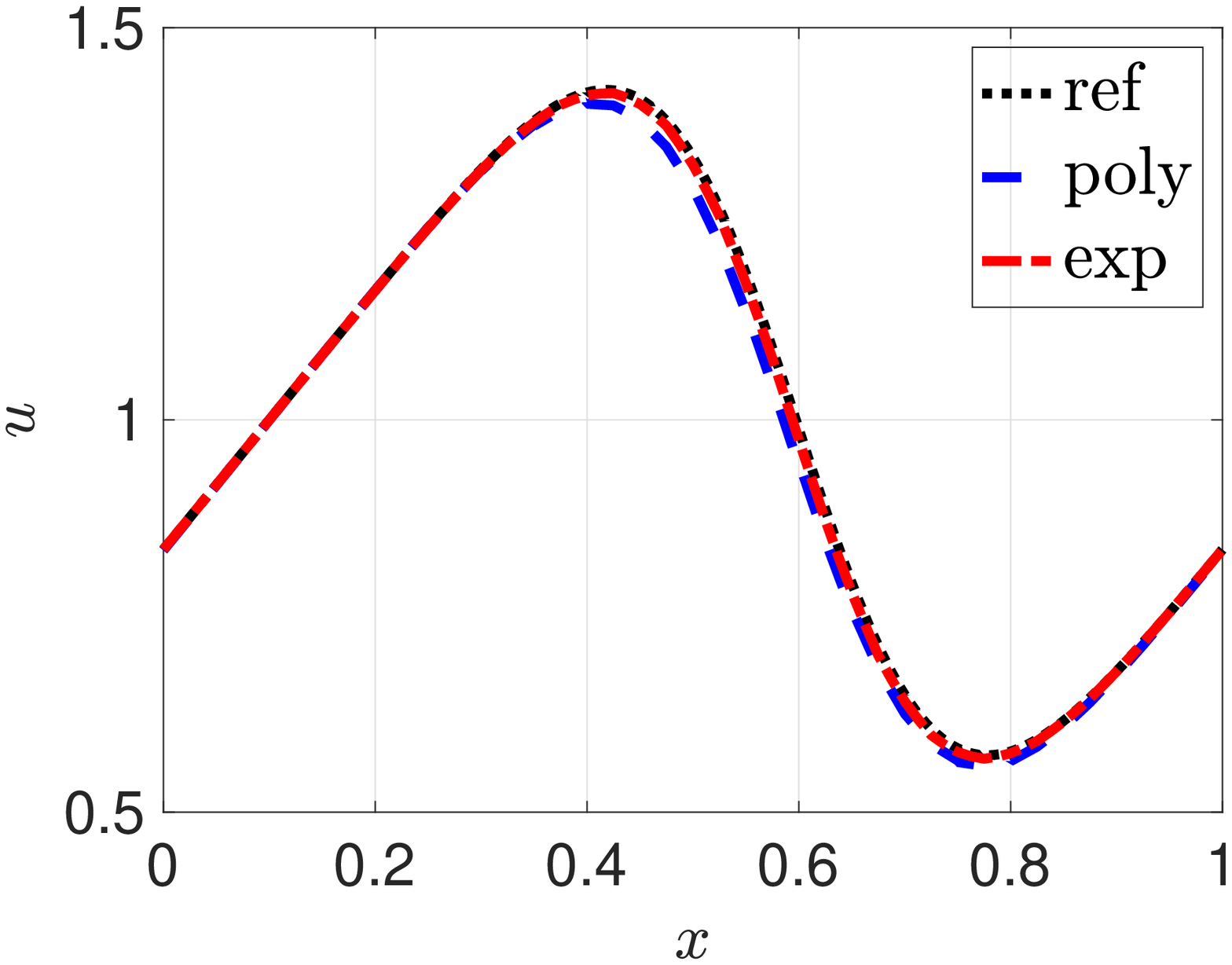} 
    		\caption{$\varepsilon = 10^{-1}$ \& $I=20$}
    		\label{fig:viscBurgers_eps_m10_I20}
  	\end{subfigure}%
	\begin{subfigure}[b]{0.33\textwidth}
		\includegraphics[width=\textwidth]{%
      		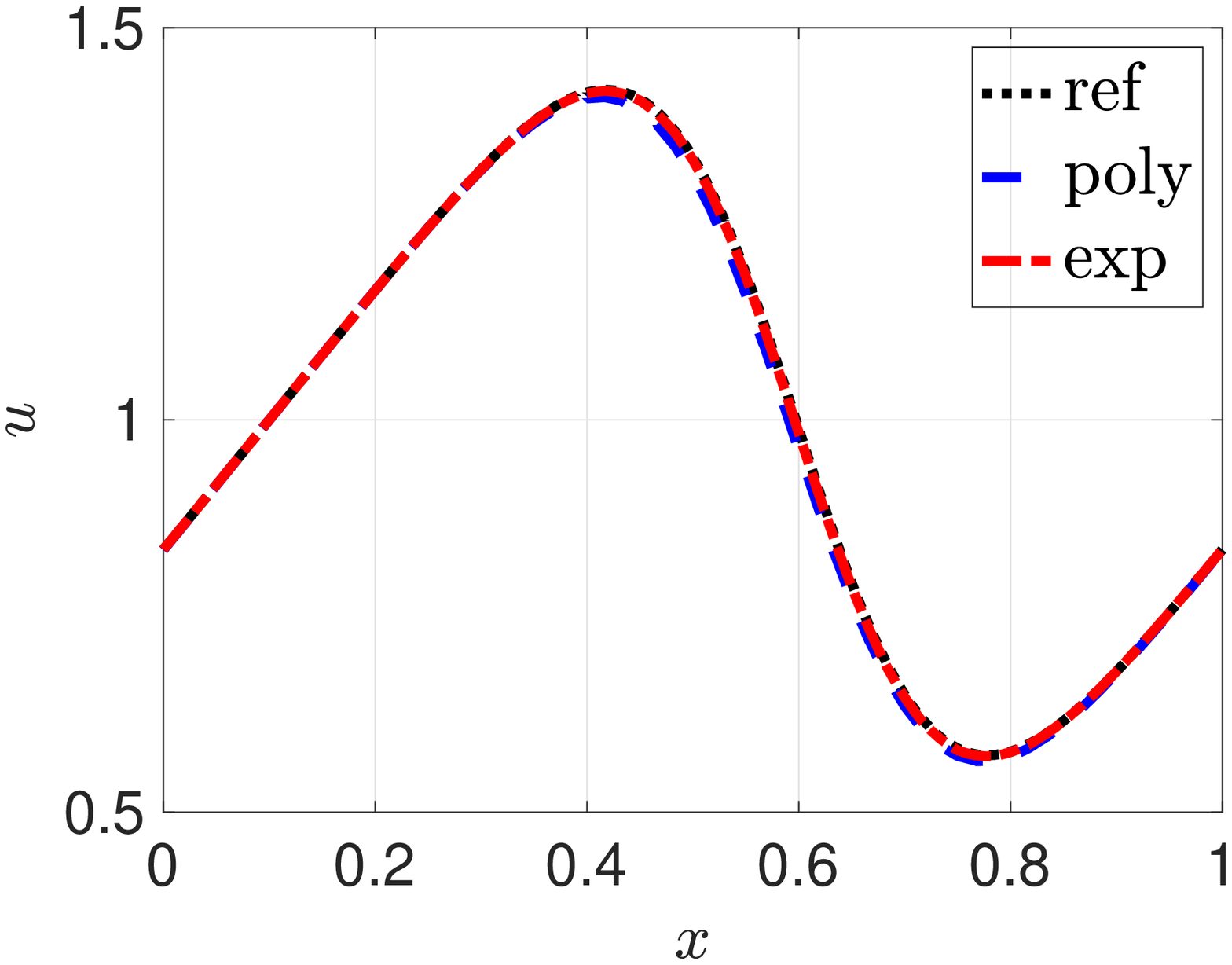} 
    		\caption{$\varepsilon = 10^{-1}$ \& $I=30$}
    		\label{fig:viscBurgers_eps_m10_I30}
  	\end{subfigure}%
	\\ 
	\begin{subfigure}[b]{0.33\textwidth}
		\includegraphics[width=\textwidth]{%
      		figures/viscBurgers_eps_m15_I10} 
    		\caption{$\varepsilon = 10^{-1.5}$ \& $I=10$}
    		\label{fig:viscBurgers_eps_m15_I10}
  	\end{subfigure}%
	\begin{subfigure}[b]{0.33\textwidth}
		\includegraphics[width=\textwidth]{%
      		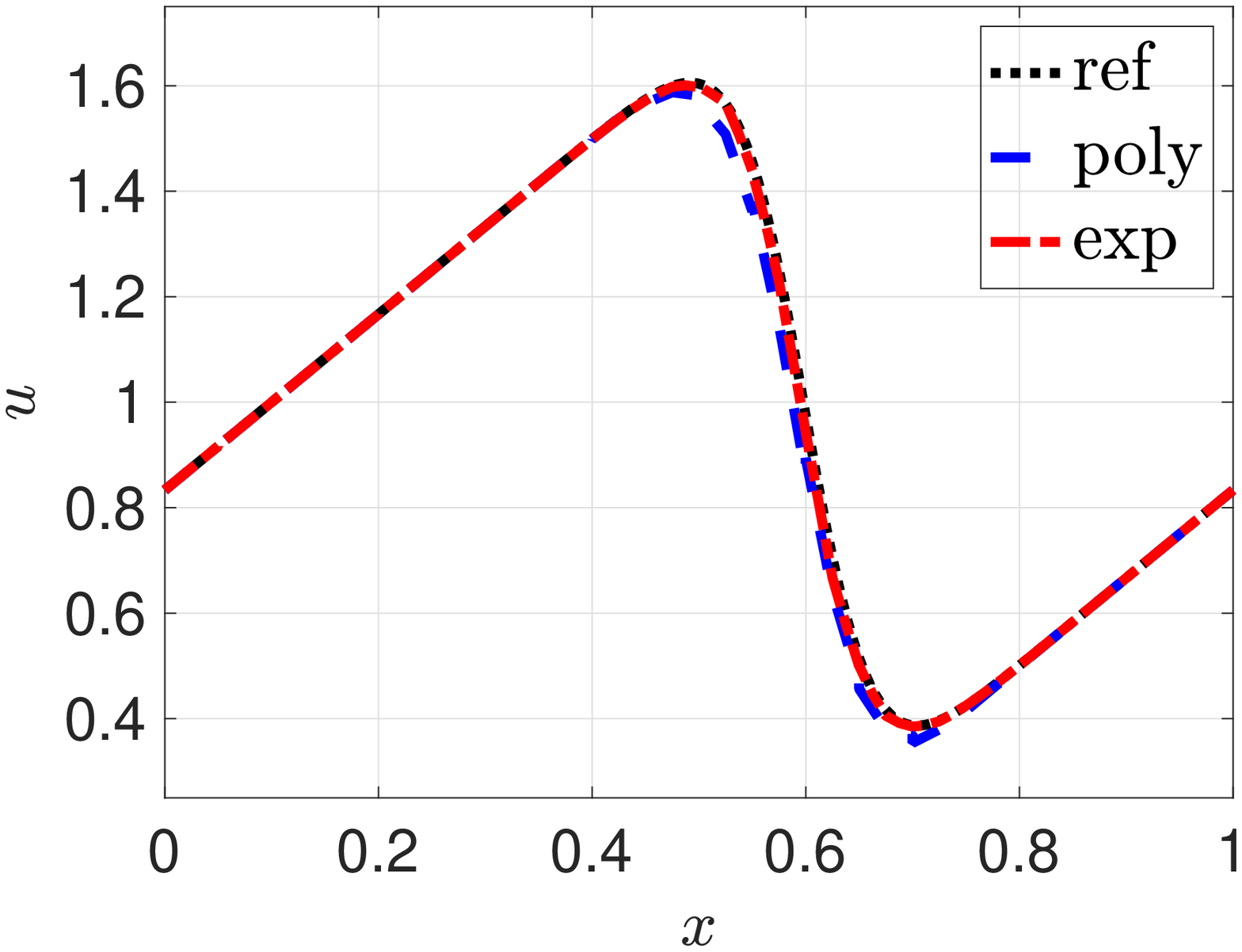} 
    		\caption{$\varepsilon = 10^{-1.5}$ \& $I=20$}
    		\label{fig:viscBurgers_eps_m15_I20}
  	\end{subfigure}%
	\begin{subfigure}[b]{0.33\textwidth}
		\includegraphics[width=\textwidth]{%
      		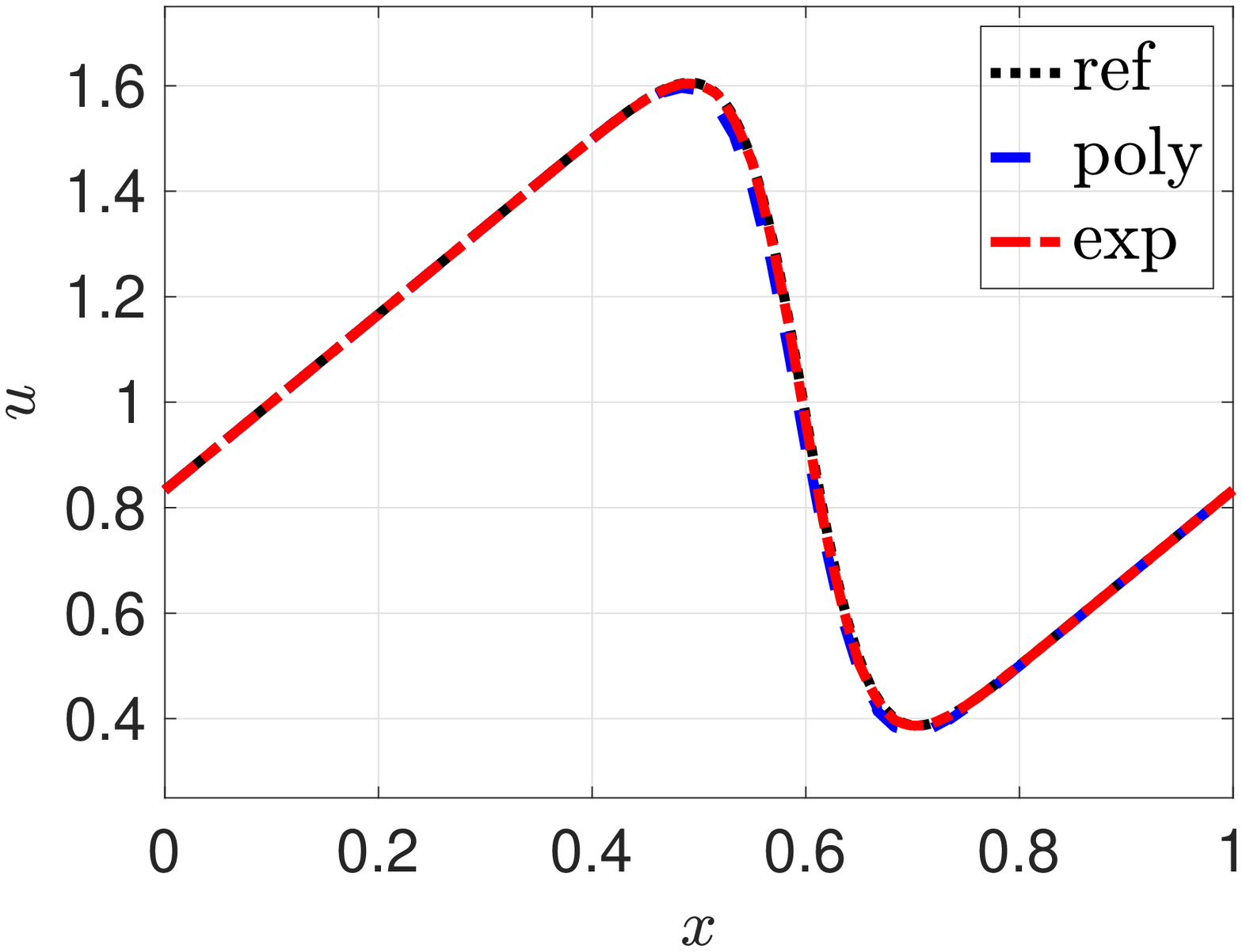} 
    		\caption{$\varepsilon = 10^{-1.5}$ \& $I=30$}
    		\label{fig:viscBurgers_eps_m15_I30}
  	\end{subfigure}%
	\\
	\begin{subfigure}[b]{0.33\textwidth}
		\includegraphics[width=\textwidth]{%
      		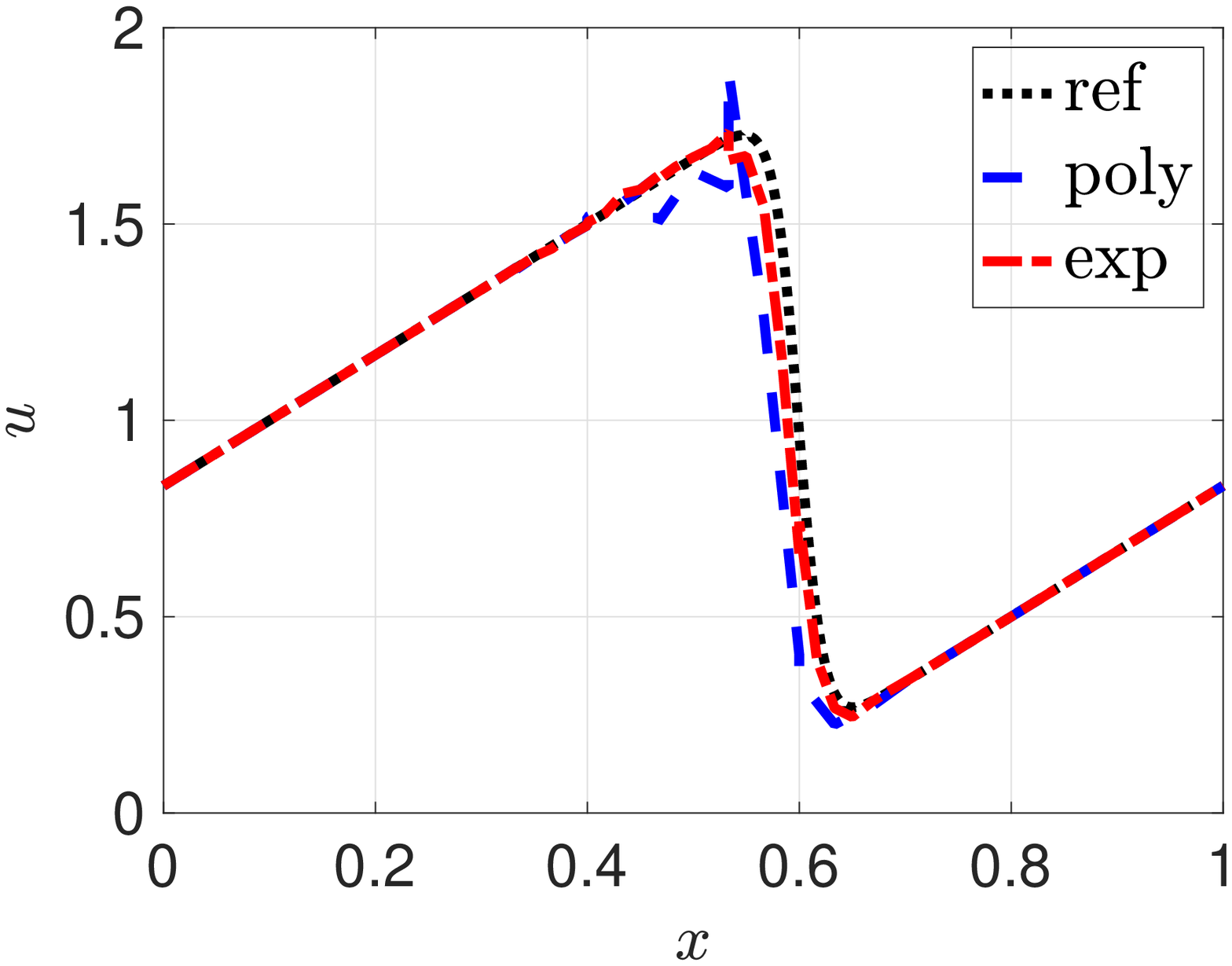} 
    		\caption{$\varepsilon = 10^{-2}$ \& $I=15$}
    		\label{fig:viscBurgers_eps_m20_I15}
  	\end{subfigure}%
	\begin{subfigure}[b]{0.33\textwidth}
		\includegraphics[width=\textwidth]{%
      		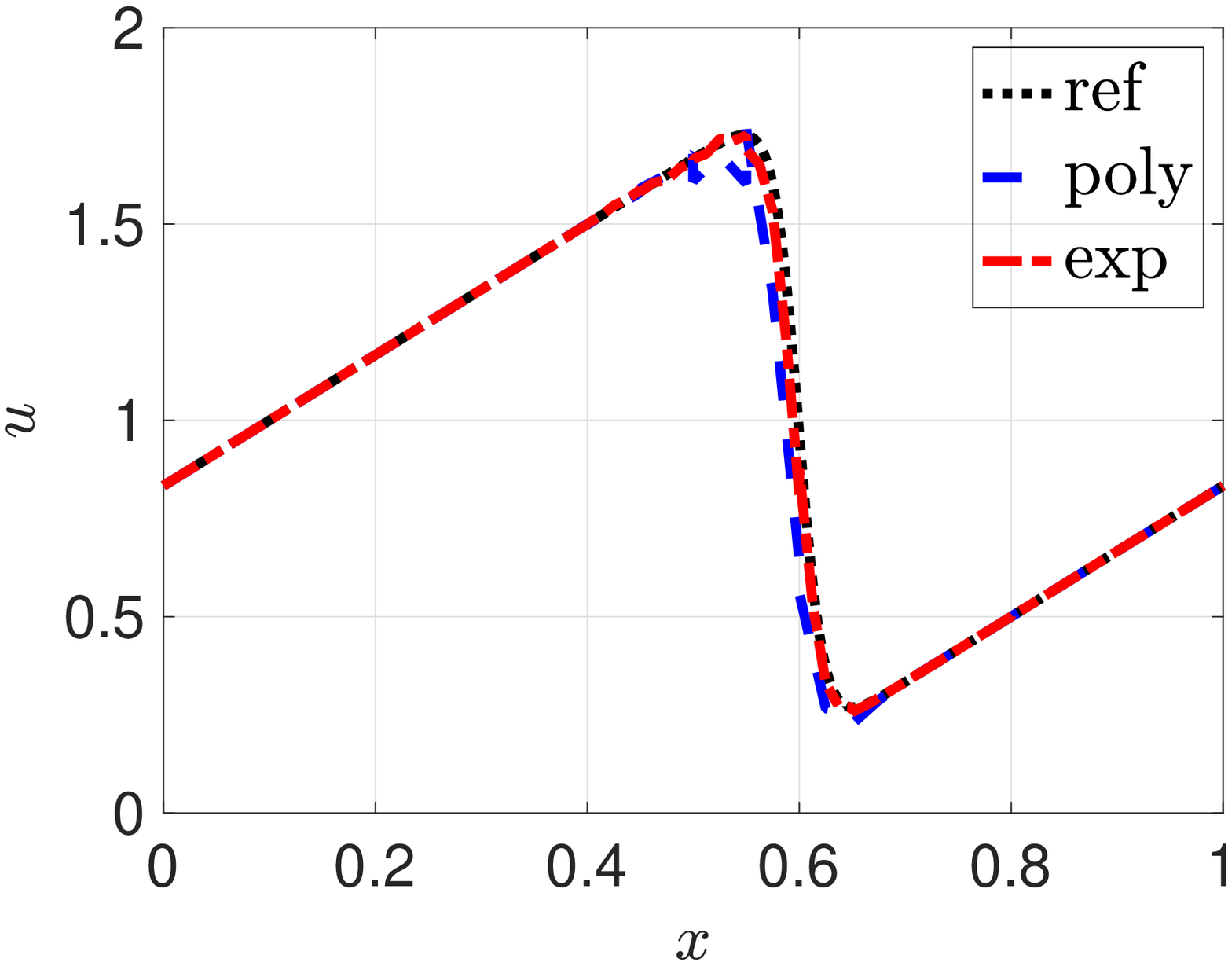} 
    		\caption{$\varepsilon = 10^{-2}$ \& $I=20$}
    		\label{fig:viscBurgers_eps_m20_I20}
  	\end{subfigure}%
	\begin{subfigure}[b]{0.33\textwidth}
		\includegraphics[width=\textwidth]{%
      		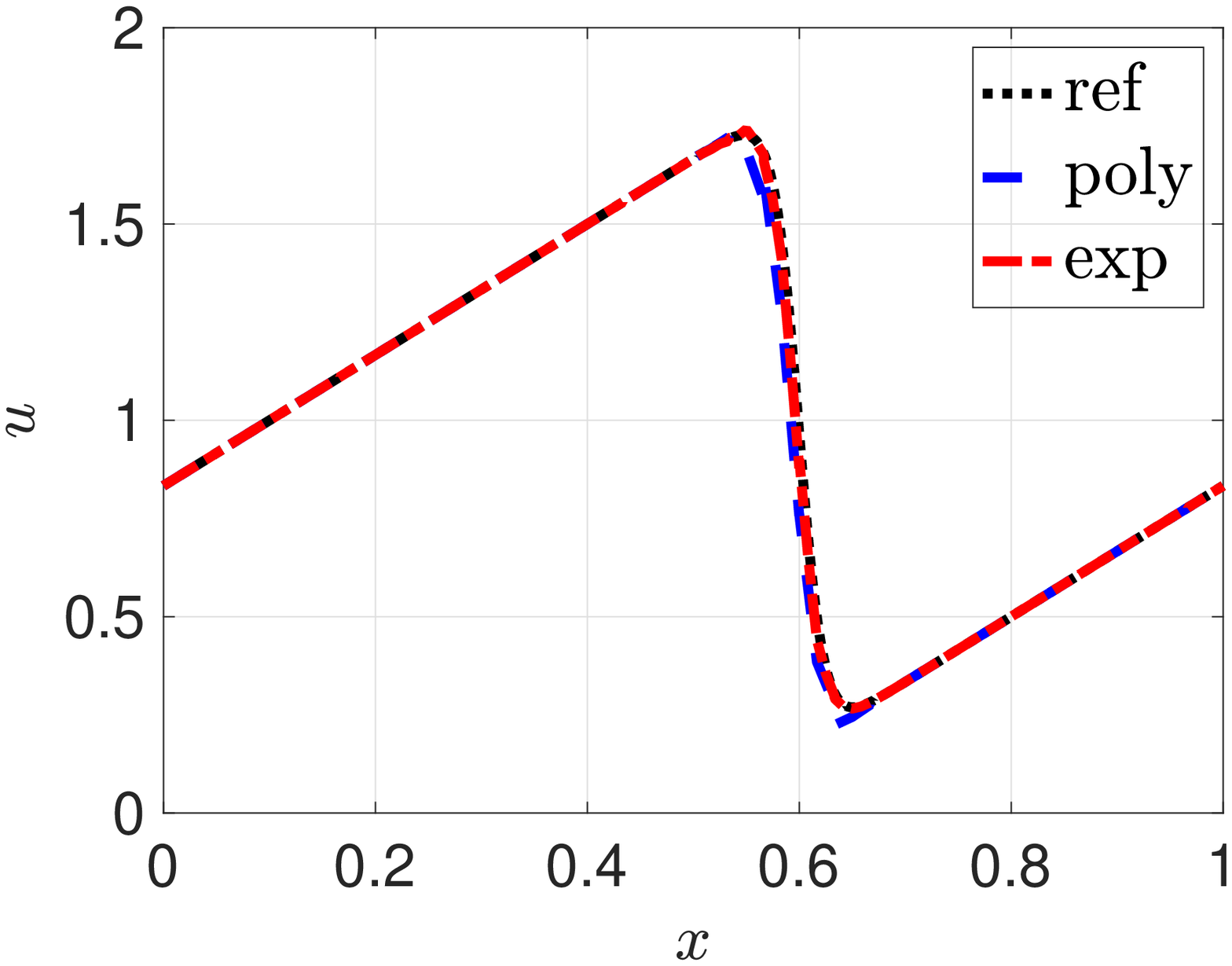} 
    		\caption{$\varepsilon = 10^{-2}$ \& $I=30$}
    		\label{fig:viscBurgers_eps_m20_I30}
  	\end{subfigure}%
  	\caption{ 
	(Numerical) solutions for the viscous Burgers equation \cref{eq:boundLayer} at $t=0.1$ for different diffusion parameters $\varepsilon$ and number of blocks $I$. 
	We used a multi-block FSBP-SAT scheme with a polynomial and exponential approximation space, $\mathcal{P}_2 = \Span\{1,x,x^2\}$ and $\mathcal{E}_2 = \Span\{1,x,e^x\}$, respectively.  
  	For $\varepsilon = 10^{-2}$, the polynomial SBP-SAT scheme did not yield a reasonable solution for $I=10$ blocks. 
	For this reason, we increased the number of blocks to $I=15$ in this case. 
	}
  	\label{fig:viscBurgers}
\end{figure}

\begin{figure}[tb]
	\centering 
	\begin{subfigure}[b]{0.49\textwidth}
		\includegraphics[width=\textwidth]{%
      		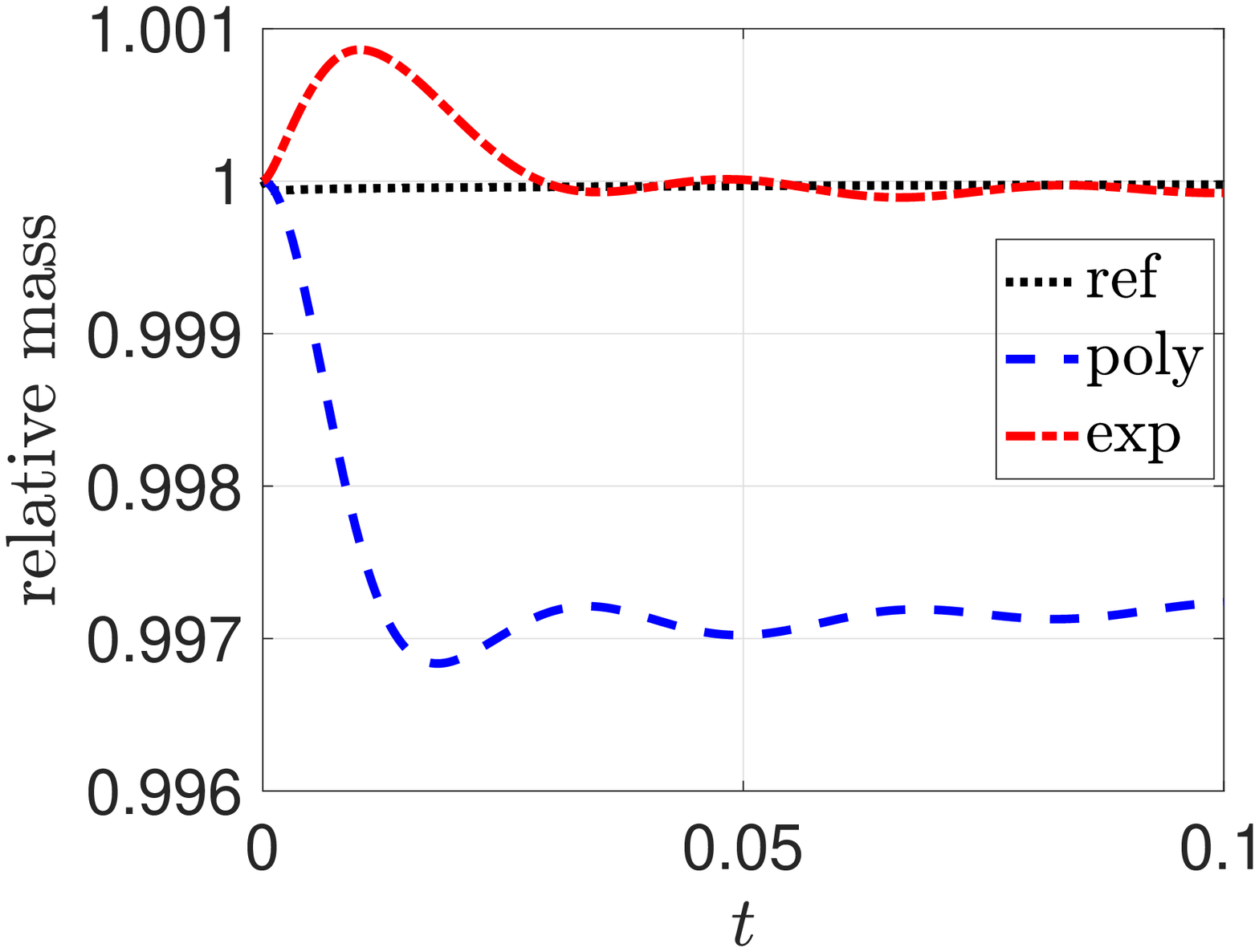} 
    		\caption{Relative mass over time}
    		\label{fig:viscBurgers_eps_m20_I30_mass}
  	\end{subfigure}%
	\begin{subfigure}[b]{0.49\textwidth}
		\includegraphics[width=\textwidth]{%
      		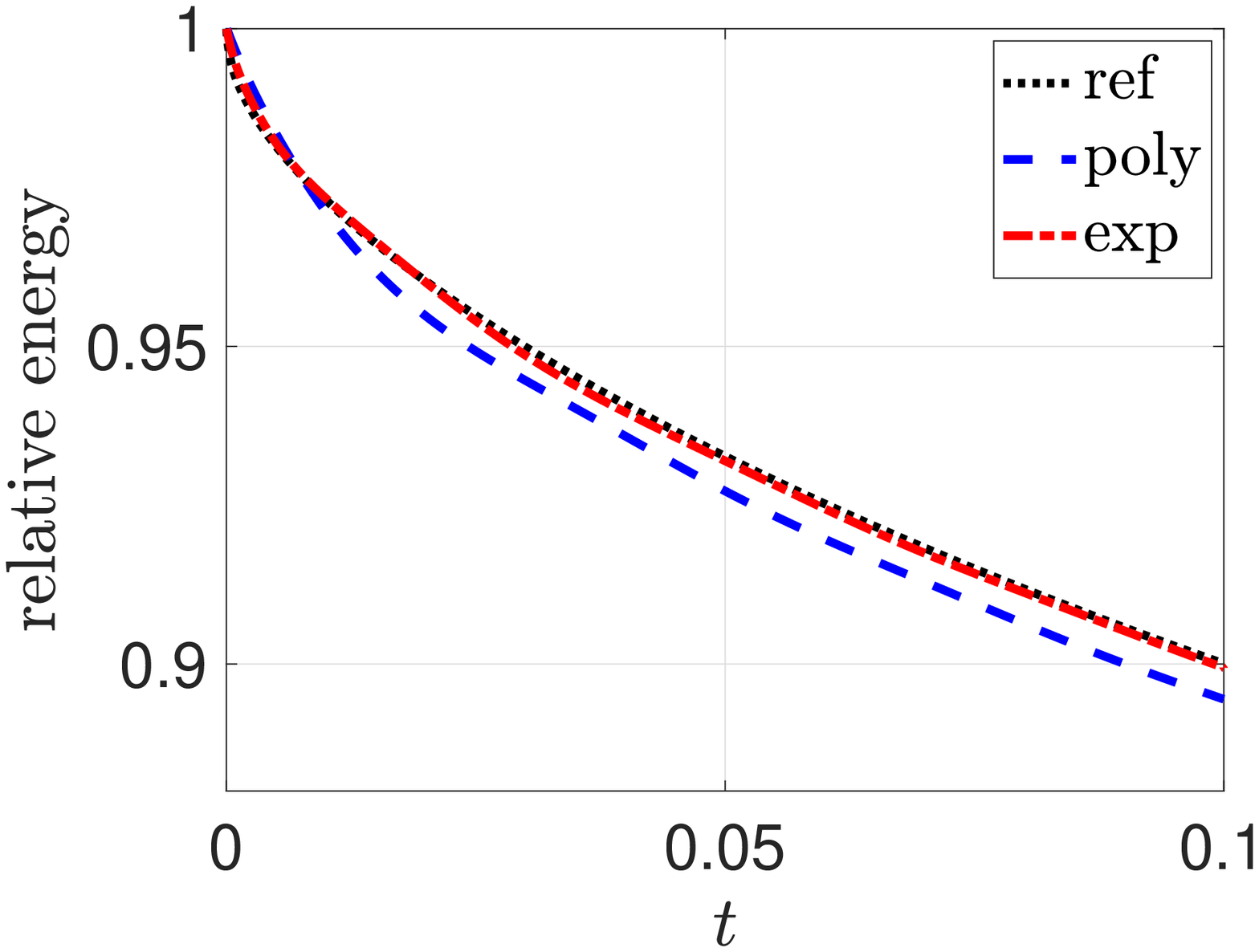} 
    		\caption{Relative energy over time}
    		\label{fig:viscBurgers_eps_m20_I30_energy}
  	\end{subfigure}%
  	\caption{ 
	Relative (discrete) mass and energy, $\left( \int u(x,t) \intd x \right) / \left( \int u(x,0) \intd x \right)$ and $\left( \int u(x,t)^2 \intd x \right) / \left( \int u(x,0)^2 \intd x \right)$, over time for the viscous Burgers equation \cref{eq:boundLayer} with $\varepsilon = 10^{-2}$ at $t=0.1$. 
	We used a multi-block FSBP-SAT scheme with $I=30$ blocks and a polynomial and exponential approximation space, $\mathcal{P}_2 = \Span\{1,x,x^2\}$ and $\mathcal{E}_2 = \Span\{1,x,e^x\}$, respectively. 
	}
  	\label{fig:viscBurgers_mass_energy}
\end{figure}

\cref{fig:viscBurgers} shows the (numerical) solutions for the viscous Burgers equation \cref{eq:boundLayer} at $t=0.1$ for different diffusion parameters $\varepsilon$ and number of blocks $I$.   
The numerical solutions were computed using a three-dimensional polynomial and exponential function space, $\mathcal{P}_2 = \Span\{1,x,x^2\}$ and $\mathcal{E}_2 = \Span\{1,x,e^x\}$, respectively.
The reference solution was computed using the polynomial multi-block FSBP-SAT method with a high number of $I=200$ blocks.
\cref{fig:viscBurgers} highlights the benefits of using the exponential function space compared to the polynomial one. 
This superiority is particularly evident in areas characterized by a steep gradient. 
More precisely, the advantage is two-fold: 
First, for a fixed diffusion parameter $\varepsilon$, each row of \cref{fig:viscBurgers} demonstrates that the solution is better and faster resolved using the exponential function space. 
Second, for a fixed grid (number of blocks $I$), each column of \cref{fig:viscBurgers} indicates that the exponential function space can resolve models with smaller $\varepsilon$. 
This is a significant observation as such models are frequently deployed as stabilized surrogates to hyperbolic conservation laws, with shock discontinuities being replaced with smooth but rapid transitions \cite{mattsson2004stable,ranocha2018stability,glaubitz2020shock}. 
The principal aim of this practice is to alter the original problem as little as possible while maintaining robustness. 
Furthermore, \cref{fig:viscBurgers_mass_energy} demonstrates that the exponential FSBP scheme's relative mass and energy profiles more closely align with those of the reference solution, as compared to the polynomial FSBP scheme. 

\subsection{The acoustic wave equation} 

The above test problems combined included both first- and second-derivative operators. 
To explicitly evaluate the performance of the second-derivative FSBP operator, we now turn our attention to the acoustic wave equation  
\begin{equation}\label{eq:wave}
	\partial_{tt} u = c^2 \partial_{xx} u 
\end{equation}  
with constant speed of sound $c > 0$ on $\Omega = [-,1,1]$ and periodic boundary conditions.  
Exact solutions of \cref{eq:wave} have the form $u(x,t) = f(x+ct) + g(x-ct)$, where $f$ and $g$ are twice-differentiable functions. 
This represents the superposition of waves moving in opposite directions along the $x$-axis at speed $c$. 
We analyze the equation using three distinct initial data functions $f$: 
\begin{align} 
	f_1(x) & = \sin\left( \pi [ x + c t ] \right), \label{eq:waveEq_f1} \\
	f_2(x) & = \exp\left\{ 100 \sin\left( k \pi [ x + c t ] \right) \right\}, \label{eq:waveEq_f2} \\
	f_3(x) & = x^2, \label{eq:waveEq_f3} 
\end{align} 
where $f_3$ is periodically extended outside of $[-1,1]$.
In all scenarios, we adopt $g(x) = \cos^2( 2 \pi [ x - c t ] )$. 
Each reference solution maintains periodicity. 
We therefore discretize \cref{eq:wave} using \emph{periodic} second-derivative SBP operators $D_2$ approximating $\partial_xx$ on $\Omega = [-1,1]$ as 
\begin{equation}\label{eq:wave_discrete}
	 \mathbf{u}_{tt} = c^2 D_2 \mathbf{u}.
\end{equation} 
Using periodic SBP operators negates the necessity for enforcing boundary conditions and allows us to concentrate on evaluating the impact of the second-derivative SBP operators $D_2$ in \cref{eq:wave_discrete}. 
To integrate \cref{eq:wave_discrete} in time, we introduce an auxiliary variable $\mathbf{v} = \mathbf{u}_t$ and reformulate \cref{eq:wave_discrete} as the first-order system 
\begin{equation}
	\begin{bmatrix} \mathbf{u}_t \\ \mathbf{v}_t \end{bmatrix} 
	= 
	\begin{bmatrix} \mathbf{v} \\ c^2 D_2 \mathbf{u} \end{bmatrix}.
\end{equation} 
The same SSPRK(3,3) time integration method used in prior tests is applied here. 
Given the periodic nature of the problem, trigonometric FSBP operators, derived from the class of trigonometric function spaces $\mathcal{T}_d$ (see \cref{sub:examples_trig}), present a natural choice.
The chosen trigonometric FSBP operator is defined across a grid of $N=2d + 2$ equidistant points spanning $\Omega = [-1,1]$. 
We compare this operator against traditional periodic FD-SBP operators of varying orders on identical grids. 
A periodic FD-SBP operator is characterized as a circulant square matrix, with each row constituting the same elements in a rightward cyclical shift relative to its predecessor. 
Specifically, the second-order periodic second-derivative FD-SBP operator, $D_2^{(2)}$, is 
\begin{equation}
	D_2^{(2)} = \frac{1}{(\Delta x)^2} 
	\begin{bmatrix} 
		-2 & 1 & & & 1 \\
		1 & -2 & 1 & & \\ 
		 & \ddots & \ddots & \ddots & \\ 
		 &  &  1 & -2 & 1 \\
		1 & & & 1 & -2 \\ 
	\end{bmatrix}
\end{equation}
where $\Delta x$ is the constant distance between the equidistant grid points on $\Omega = [x_L,x_R]$. 
We say that $D_2^{(2)}$ is generated by the vector $[ 1, -2, 1 ]$. 
Similarly, the periodic second-derivative FD-SBP operators of order four and six, $D_2^{(4)}$ and $D_2^{(6)}$, are generated by the vector $[-1/12, 4/3, -5/2, 4/3, -1/12]$ and $[1/90, -3/20, 3/2, -49/18, 3/2, -3/20, 1/90]$, respectively. 

\begin{figure}[tb]
	\centering 
	\begin{subfigure}[b]{0.33\textwidth}
		\includegraphics[width=\textwidth]{%
      		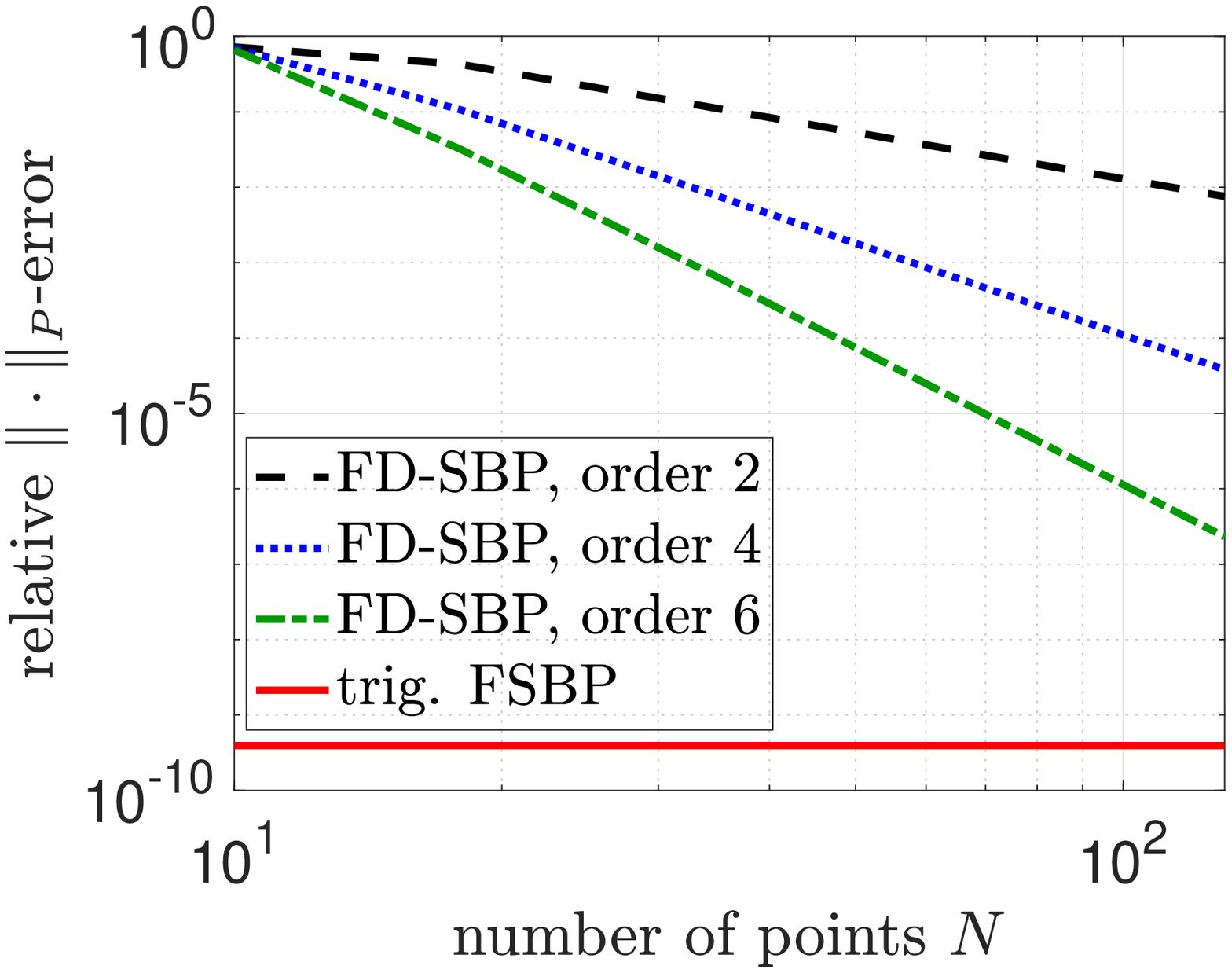} 
    		\caption{Errors vs.\ $N$ for $f_1$ in \cref{eq:waveEq_f1}}
    		\label{fig:waveEq_error_vs_N_IC1}
  	\end{subfigure}%
	\begin{subfigure}[b]{0.33\textwidth}
		\includegraphics[width=\textwidth]{%
      		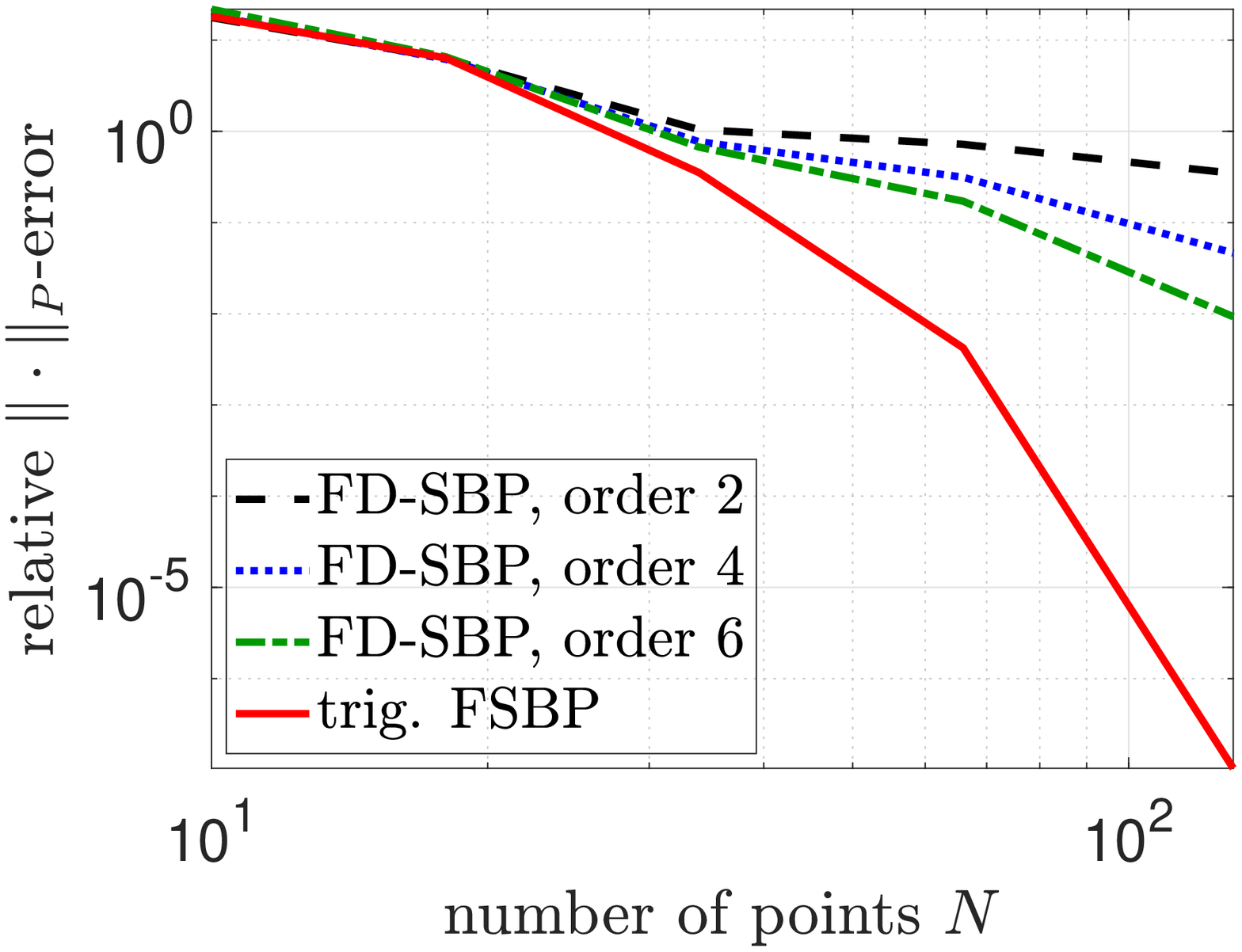} 
    		\caption{Errors vs.\ $N$ for $f_2$ in \cref{eq:waveEq_f2}}
    		\label{fig:waveEq_error_vs_N_IC2}
  	\end{subfigure}%
	\begin{subfigure}[b]{0.33\textwidth}
		\includegraphics[width=\textwidth]{%
      		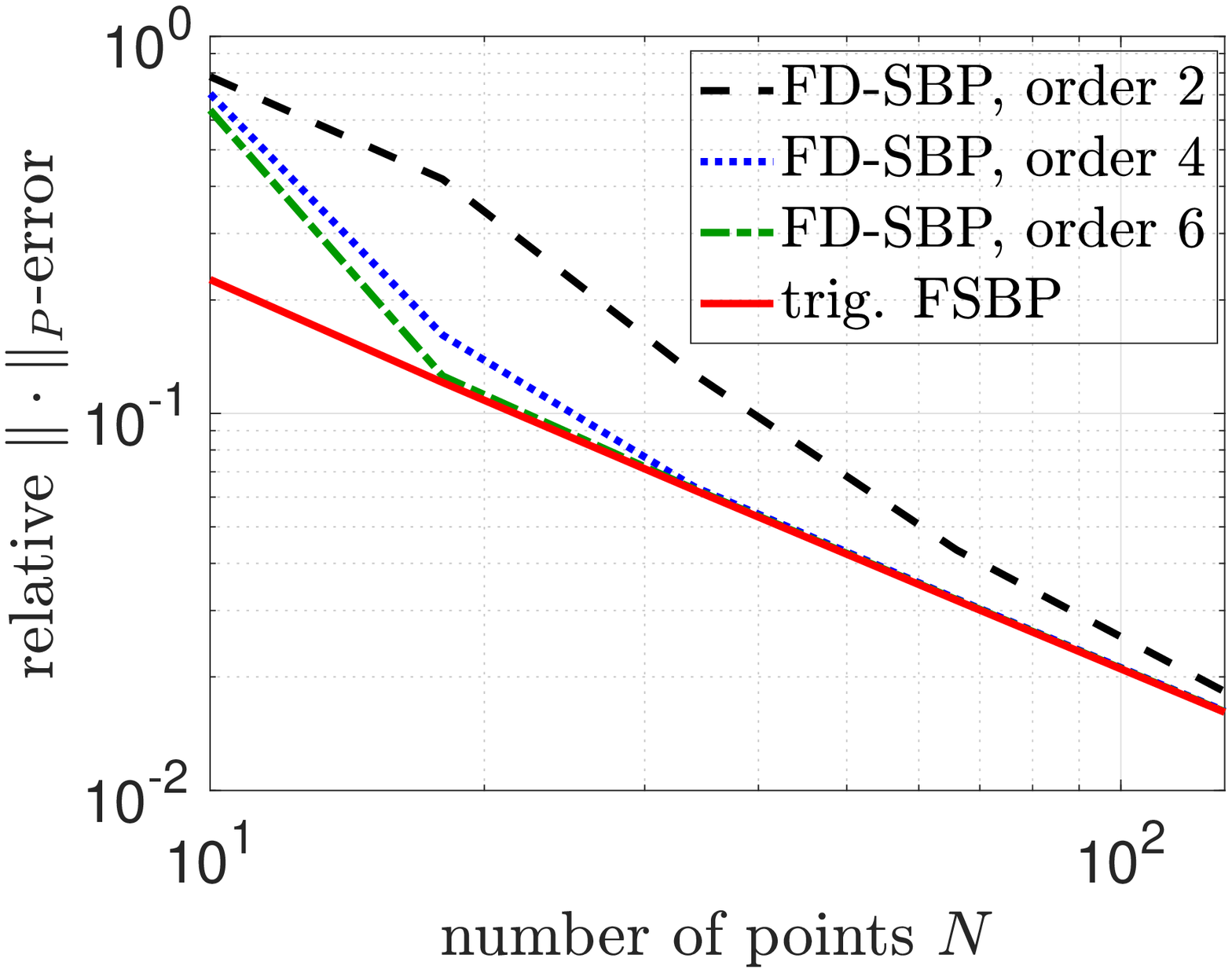} 
    		\caption{Errors vs.\ $N$ for $f_3$ in \cref{eq:waveEq_f3}}
    		\label{fig:waveEq_error_vs_N_IC3}
  	\end{subfigure}%
	\\ 
	\begin{subfigure}[b]{0.33\textwidth}
		\includegraphics[width=\textwidth]{%
      		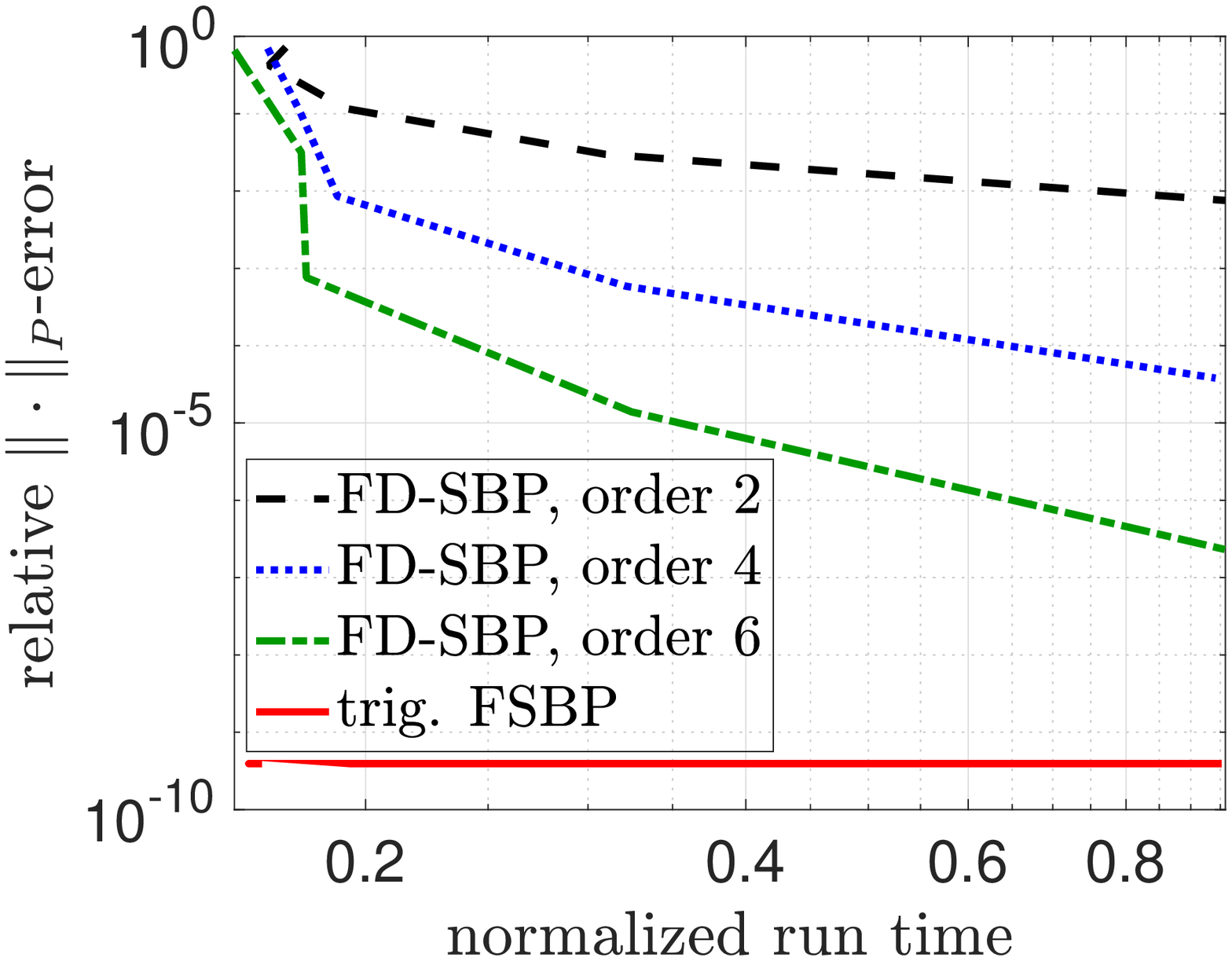} 
    		\caption{Errors vs.\ run time for $f_1$ in \cref{eq:waveEq_f1}}
    		\label{fig:waveEq_error_vs_runTime_IC1}
  	\end{subfigure}%
	\begin{subfigure}[b]{0.33\textwidth}
		\includegraphics[width=\textwidth]{%
      		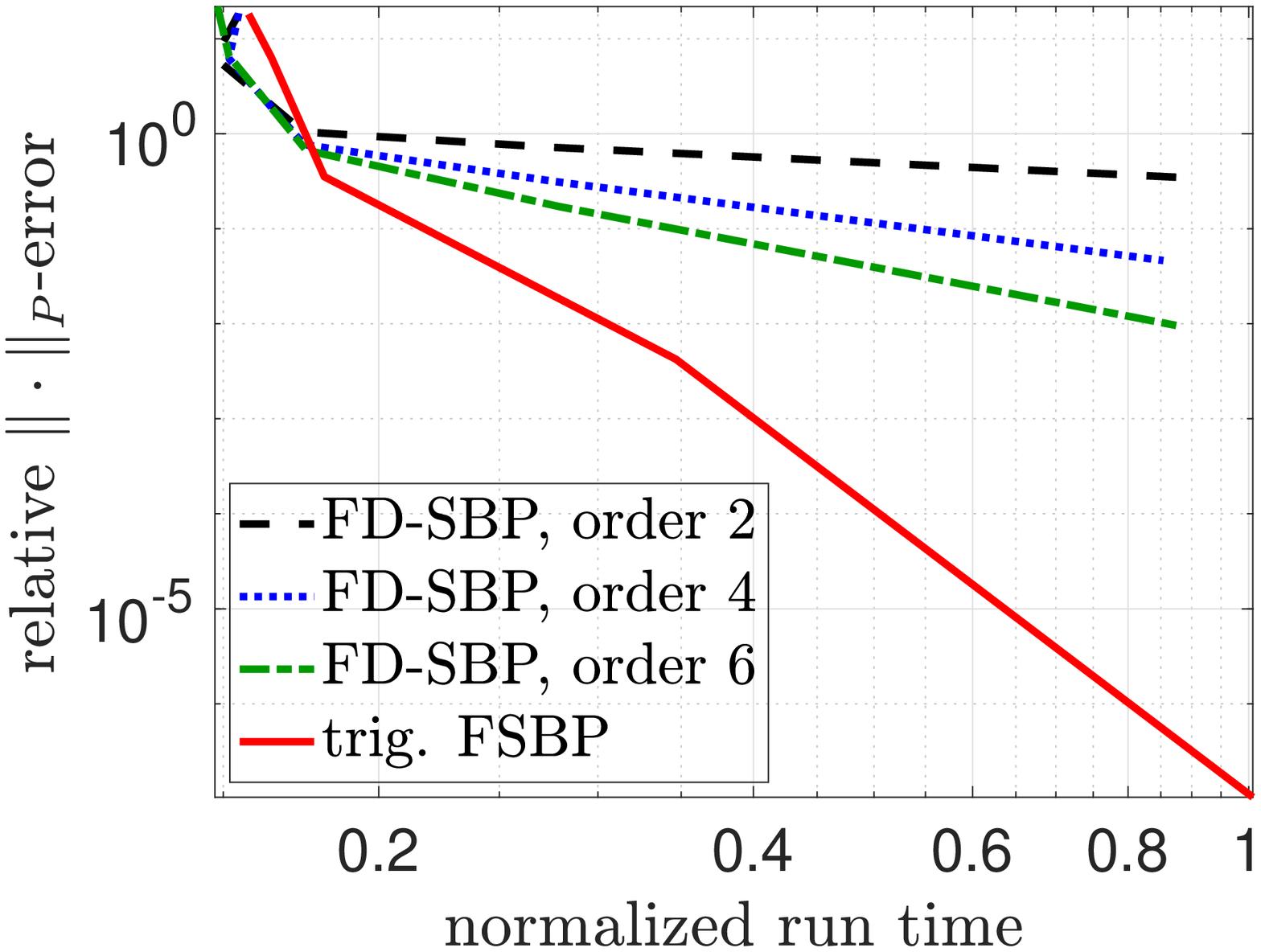} 
    		\caption{Errors vs.\ run time for $f_2$ in \cref{eq:waveEq_f2}}
    		\label{fig:waveEq_error_vs_runTime_IC2}
  	\end{subfigure}%
	\begin{subfigure}[b]{0.33\textwidth}
		\includegraphics[width=\textwidth]{%
      		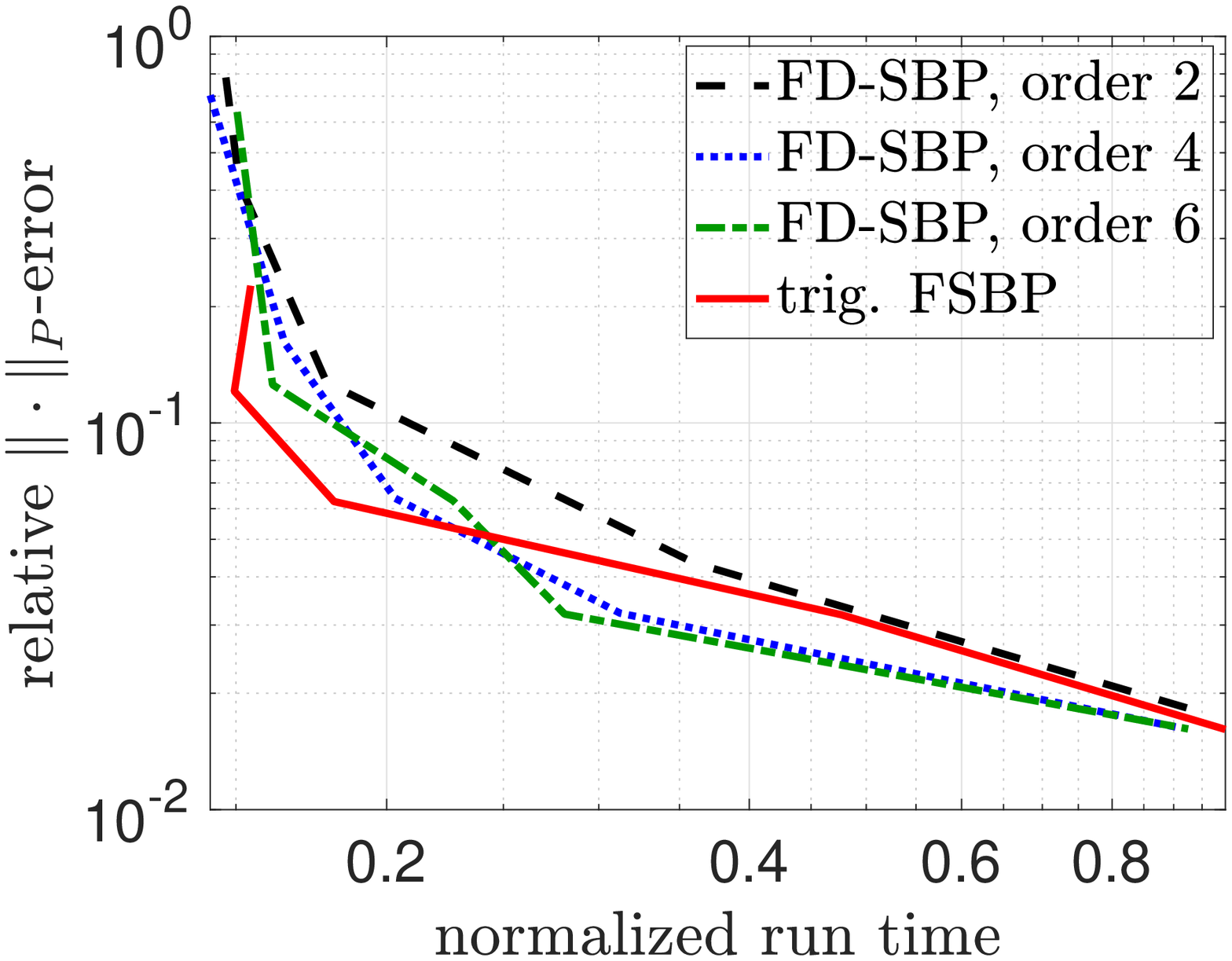} 
    		\caption{Errors vs.\ run time for $f_3$ in \cref{eq:waveEq_f3}}
    		\label{fig:waveEq_error_vs_runTime_IC3}
  	\end{subfigure}%
  	\caption{ 
	Relative $\| \cdot \|_P$-errors for the numerical solutions of the wave equation \cref{eq:wave} for $c=1$ at time $t=1$, plotted against the number of equidistant points $N$ (top row) and normalized run times (bottom row). 
	We consider the three initial data functions $f_1,f_2,f_3$ in \cref{eq:waveEq_f1,eq:waveEq_f2,eq:waveEq_f3}, respectively.
	We compare the trigonometric FSBP operators based on $\mathcal{T}_d$ (see \cref{sub:examples_trig}), defined on $N=2d + 2$ equidistant points, with traditional periodic FD-SBP operators of orders two, four, and six on the same grid. 
	}
  	\label{fig:waveEq_errors}
\end{figure}

\Cref{fig:waveEq_errors} presents the relative $\| \cdot \|_P$-errors for the numerical solutions of the wave equation \cref{eq:wave} for $c=1$ at time $t=1$, plotted against the number of equidistant points $N$ (top row) and normalized run times (bottom row). 
The relative $\| \cdot \|_P$-error, calculated as $\| \mathbf{u} - \mathbf{u}_{\text{ref}} \|_P / \| \mathbf{u}_{\text{ref}} \|_P$, measures the discrepancy between the numerical solution $\mathbf{u}$ and the reference solution's nodal vector $\mathbf{u}_{\text{ref}}$. 
Here, $\| \cdot \|_P$ denotes the discrete norm induced by the norm matrix $P$ associated with the SBP operator. 
The normalized run times are calculated as the separate run times divided by the maximum run time among all simulations (different numbers of grid points and operators).  
We consider the three initial data functions $f_1,f_2,f_3$ in \cref{eq:waveEq_f1,eq:waveEq_f2,eq:waveEq_f3}, respectively. 
To focus on the accuracy of the spatial semi-discretization, we utilize a time step size $\Delta t = 10^{-4}$ across all scenarios. 
The trigonometric initial data $f_1$ in \cref{eq:waveEq_f1} yields a trigonometric solution. 
Given the exact representability of this solution within the trigonometric function space $\mathcal{T}_d$, the trigonometric FSBP operator yields a constant and markedly lower error compared to FD-SBP operators, as observed in \cref{fig:waveEq_error_vs_N_IC1,fig:waveEq_error_vs_runTime_IC1}.
For the second initial data function $f_2$ in \cref{eq:waveEq_f2}, the trigonometric function space $\mathcal{T}_d$ cannot exactly represent the solution.
As illustrated in \cref{fig:waveEq_error_vs_N_IC2,fig:waveEq_error_vs_runTime_IC2}, this results in less accurate numerical solutions for the trigonometric FSBP operator. 
Despite this limitation, the trigonometric FSBP operator still demonstrates improved accuracy and convergence rates over FD-SBP operators.
The third initial data function $f_3$ in \cref{eq:waveEq_f3} remains periodic but features a discontinuous derivative at the domain boundaries. 
This discontinuity propagates into the domain's interior over time, diminishing the accuracy of numerical solutions for both trigonometric FSBP and FD-SBP operators of various orders. 
Particularly, in \cref{fig:waveEq_error_vs_runTime_IC3}, the advantage of using the trigonometric FSBP operator over traditional FD-SBP operators becomes negligible. 
This final analysis underscores that no universal superiority exists between function spaces (and their corresponding FSBP operators); rather, the selection depends on the specific characteristics and requirements of the problem at hand.

\begin{figure}[tb]
	\centering 
	\begin{subfigure}[b]{0.495\textwidth}
		\includegraphics[width=\textwidth]{%
      		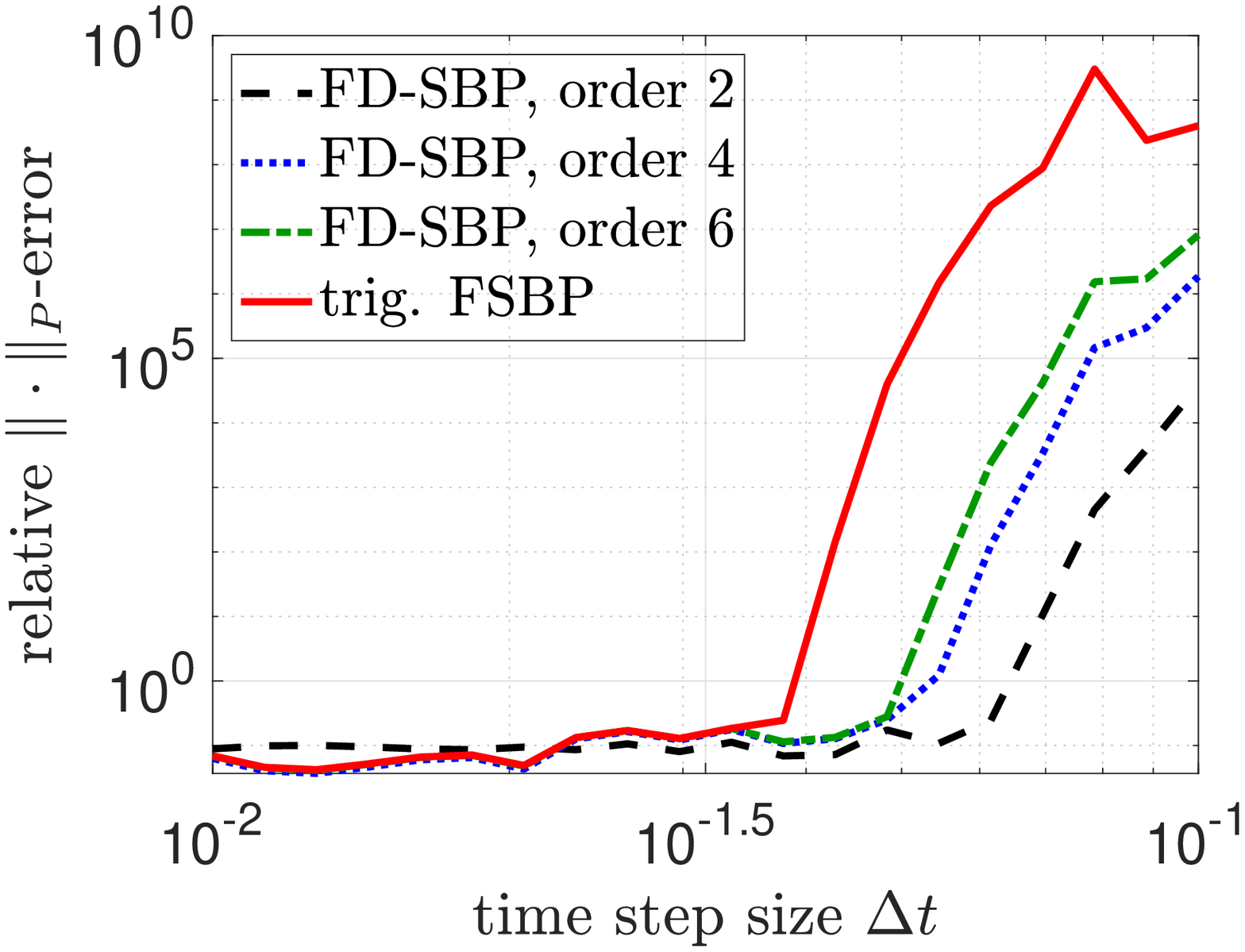} 
    		\caption{Errors vs.\ time step size $\Delta t$, initial data 2}
    		\label{fig:waveEq_error_vs_timeStep_IC2}
  	\end{subfigure}%
	\begin{subfigure}[b]{0.495\textwidth}
		\includegraphics[width=\textwidth]{%
      		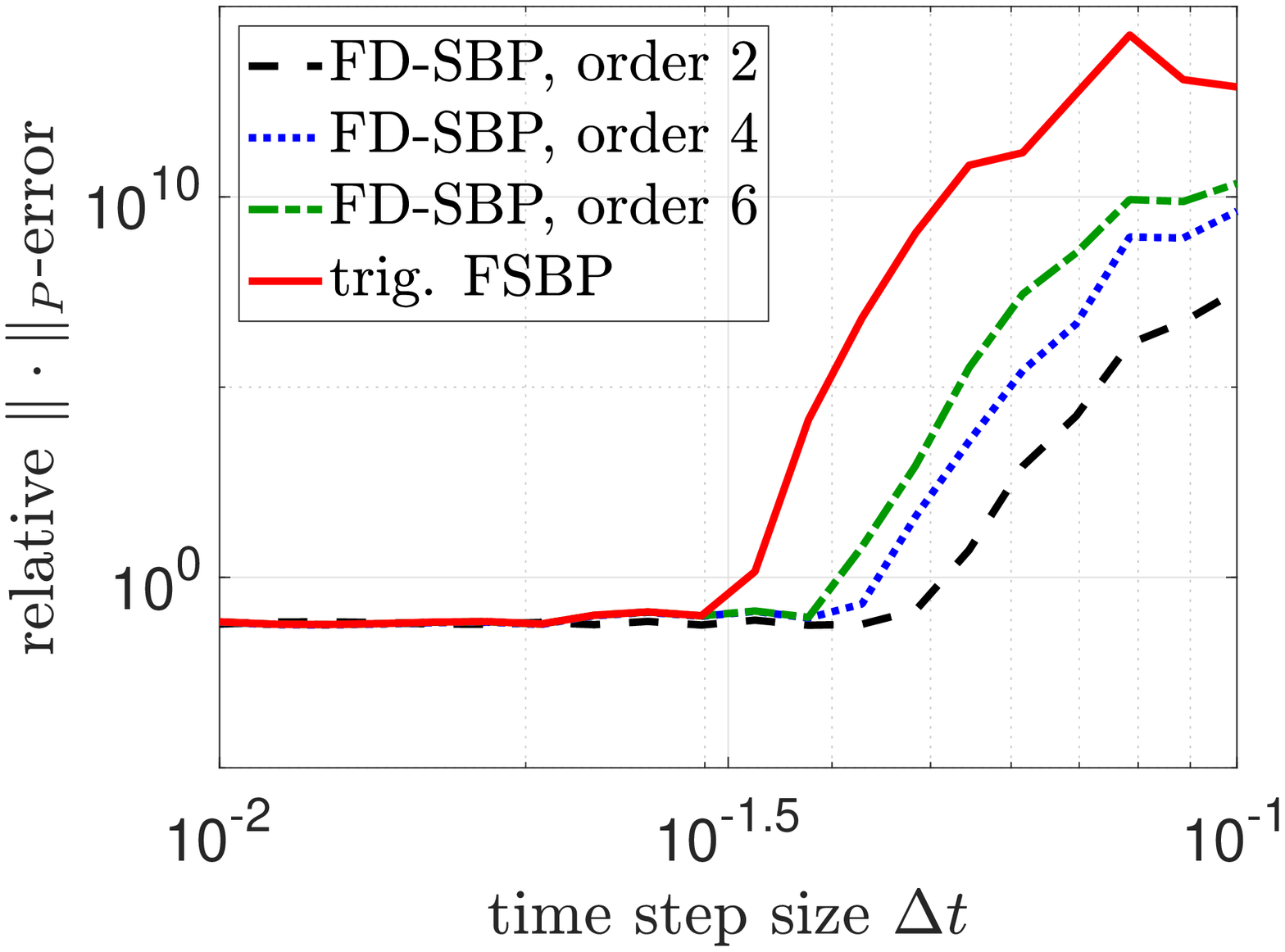} 
    		\caption{Errors vs.\ time step size $\Delta t$, initial data 3}
    		\label{fig:waveEq_error_vs_timeStep_IC3}
  	\end{subfigure}%
  	\caption{ 
	Relative $\| \cdot \|_P$-errors for the numerical solutions of the wave equation \cref{eq:wave} for $c=1$ at time $t=1$, plotted against the time step size $\Delta t$ used in the SSPRK(3,3) method for time integration. 
	We consider the two initial data functions $f_1$ and $f_2$ in \cref{eq:waveEq_f2,eq:waveEq_f3}.
	We compare the trigonometric FSBP operators based on $\mathcal{T}_20$ (see \cref{sub:examples_trig}), defined on $N=22$ equidistant points, with traditional periodic FD-SBP operators of orders two, four, and six on the same grid.
	}
  	\label{fig:waveEq_timeStep}
\end{figure}

\begin{figure}[tb]
	\centering 
	\begin{subfigure}[b]{0.33\textwidth}
		\includegraphics[width=\textwidth]{%
      		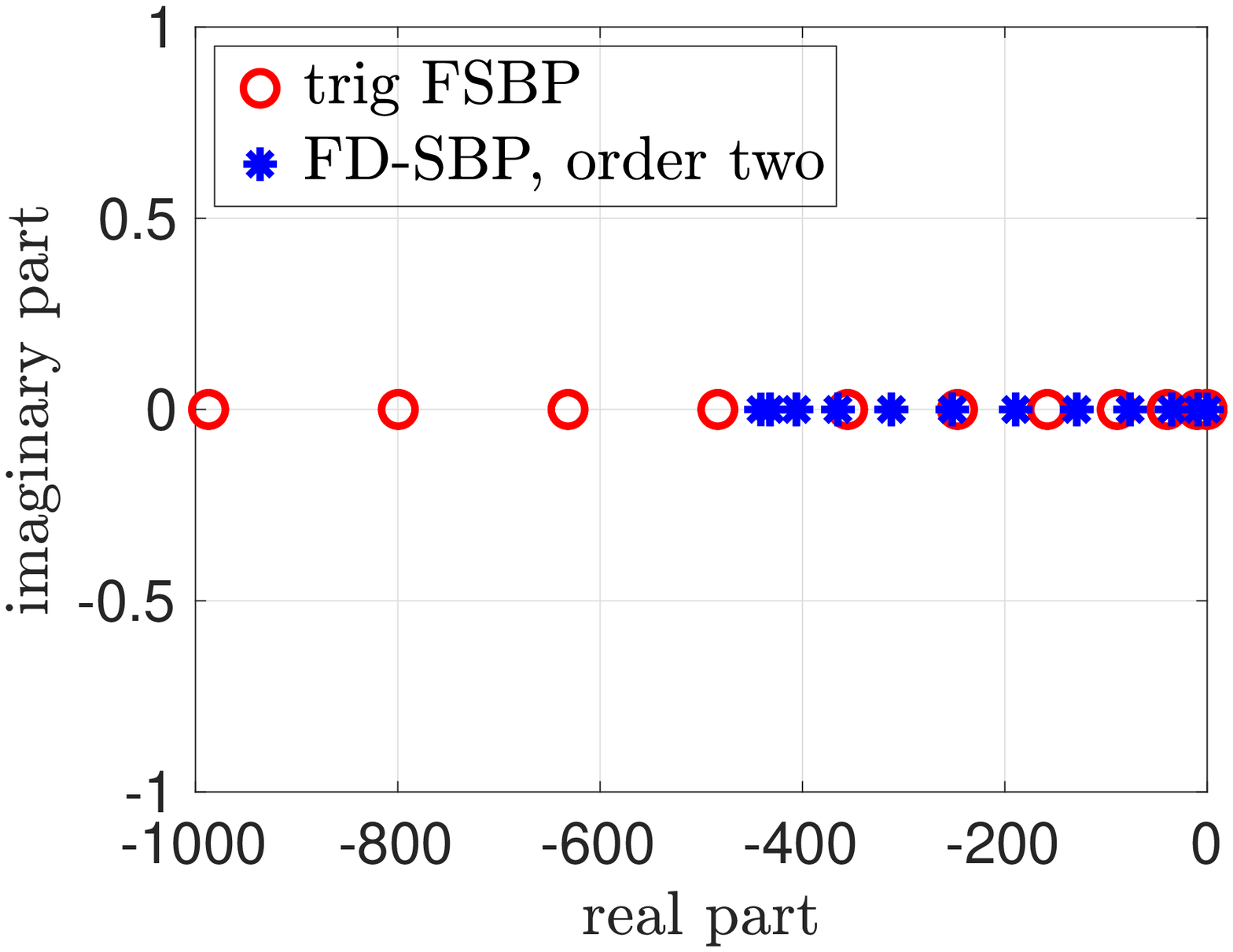} 
    		\caption{$d=10$ ($N=22$ equidistant points)}
    		\label{fig:spectrum_D2_d10_order2}
  	\end{subfigure}%
	\begin{subfigure}[b]{0.33\textwidth}
		\includegraphics[width=\textwidth]{%
      		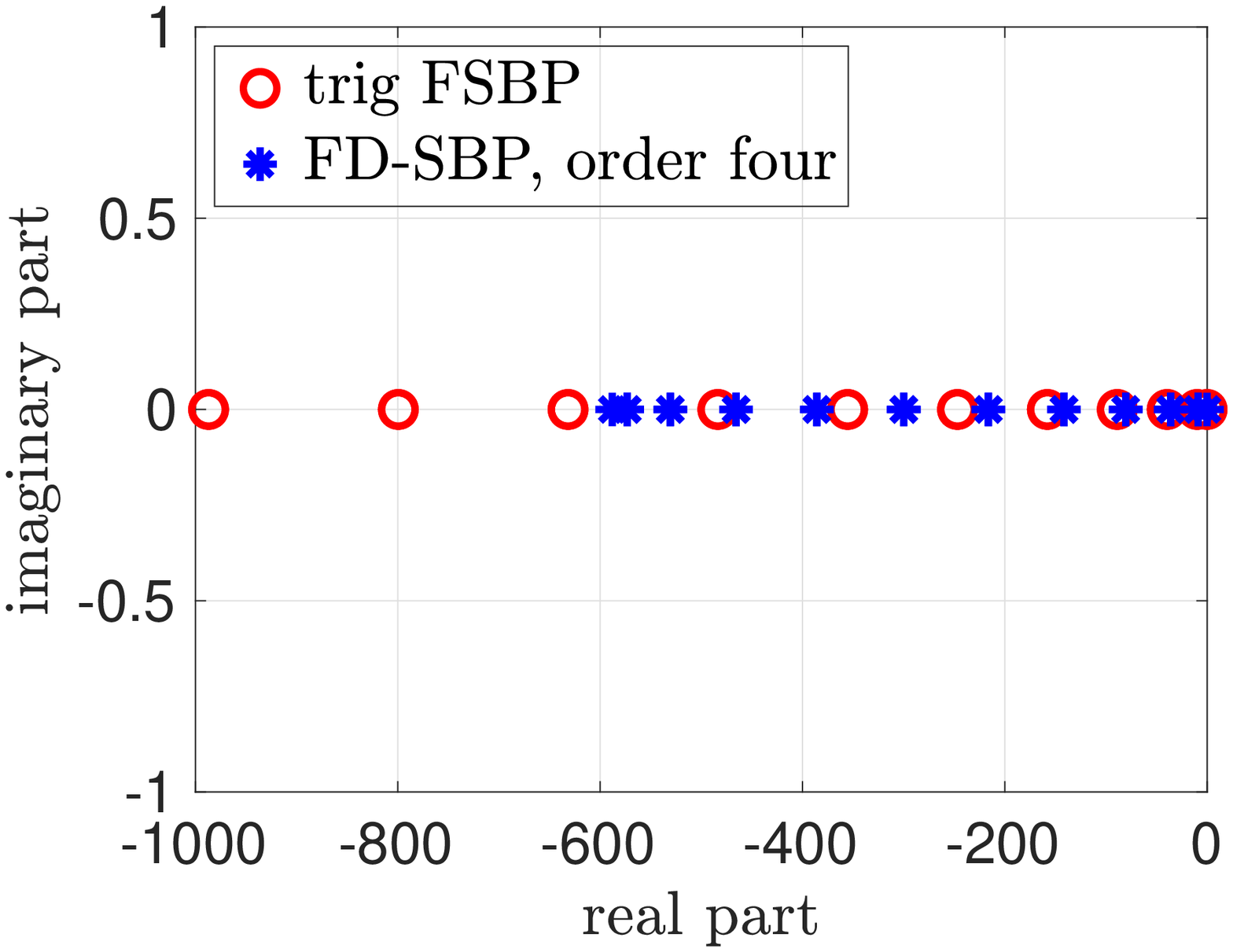} 
    		\caption{$d=10$ ($N=22$ equidistant points)}
    		\label{fig:spectrum_D2_d10_order4}
  	\end{subfigure}%
	\begin{subfigure}[b]{0.33\textwidth}
		\includegraphics[width=\textwidth]{%
      		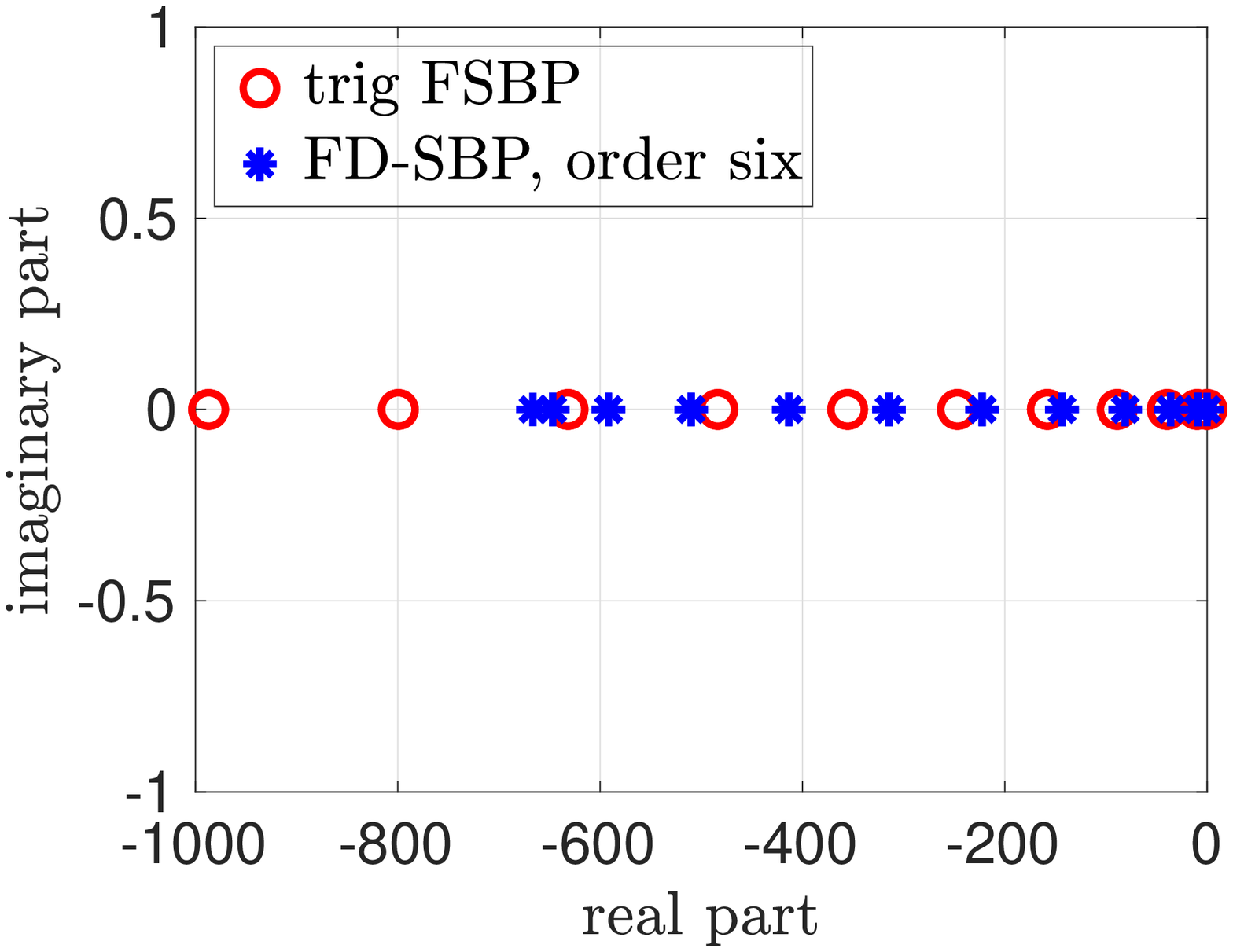} 
    		\caption{$d=10$ ($N=22$ equidistant points)}
    		\label{fig:spectrum_D2_d10_order6}
  	\end{subfigure}%
	\\ 
	\begin{subfigure}[b]{0.33\textwidth}
		\includegraphics[width=\textwidth]{%
      		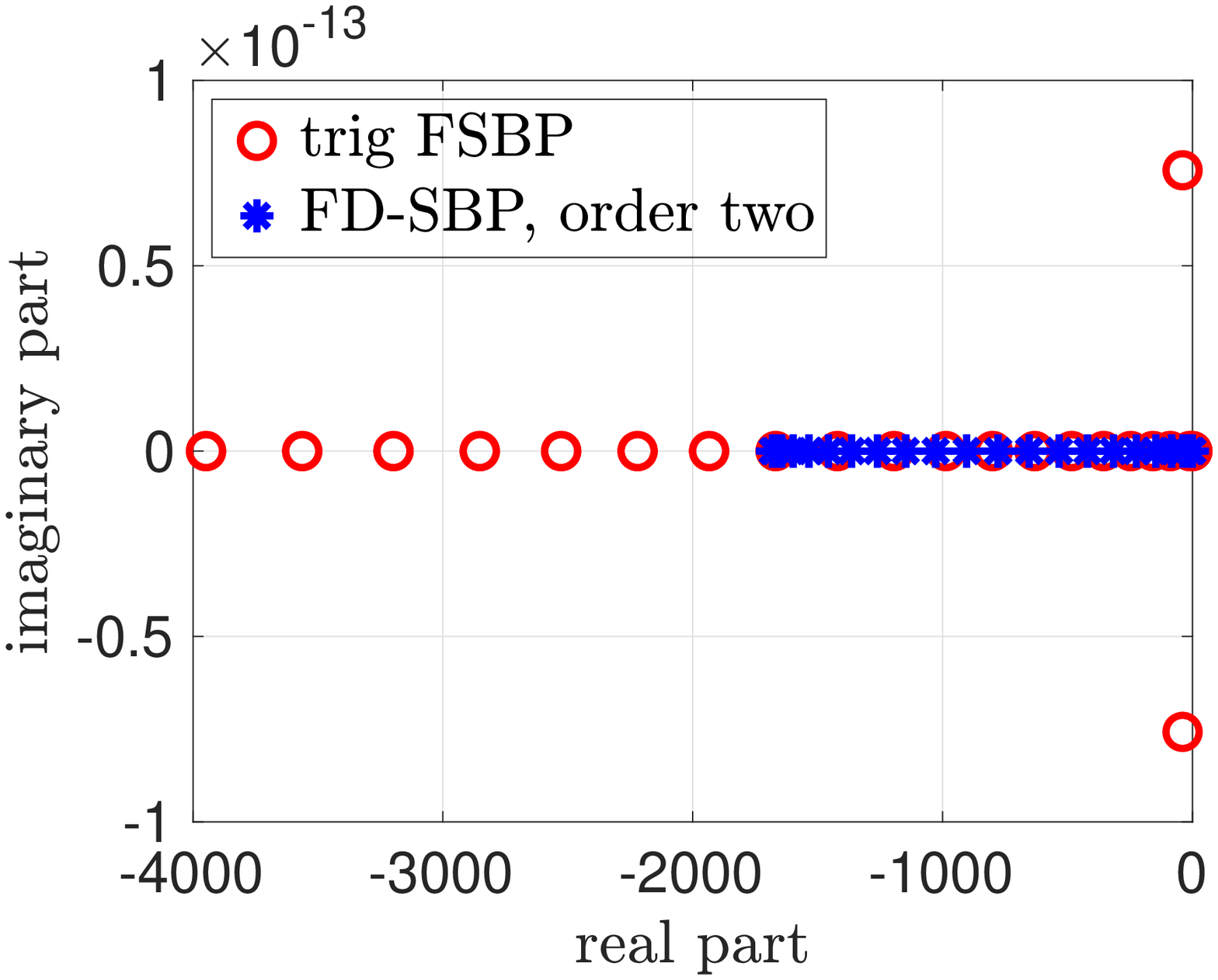} 
    		\caption{$d=20$ ($N=42$ equidistant points)}
    		\label{fig:spectrum_D2_d20_order2}
  	\end{subfigure}%
	\begin{subfigure}[b]{0.33\textwidth}
		\includegraphics[width=\textwidth]{%
      		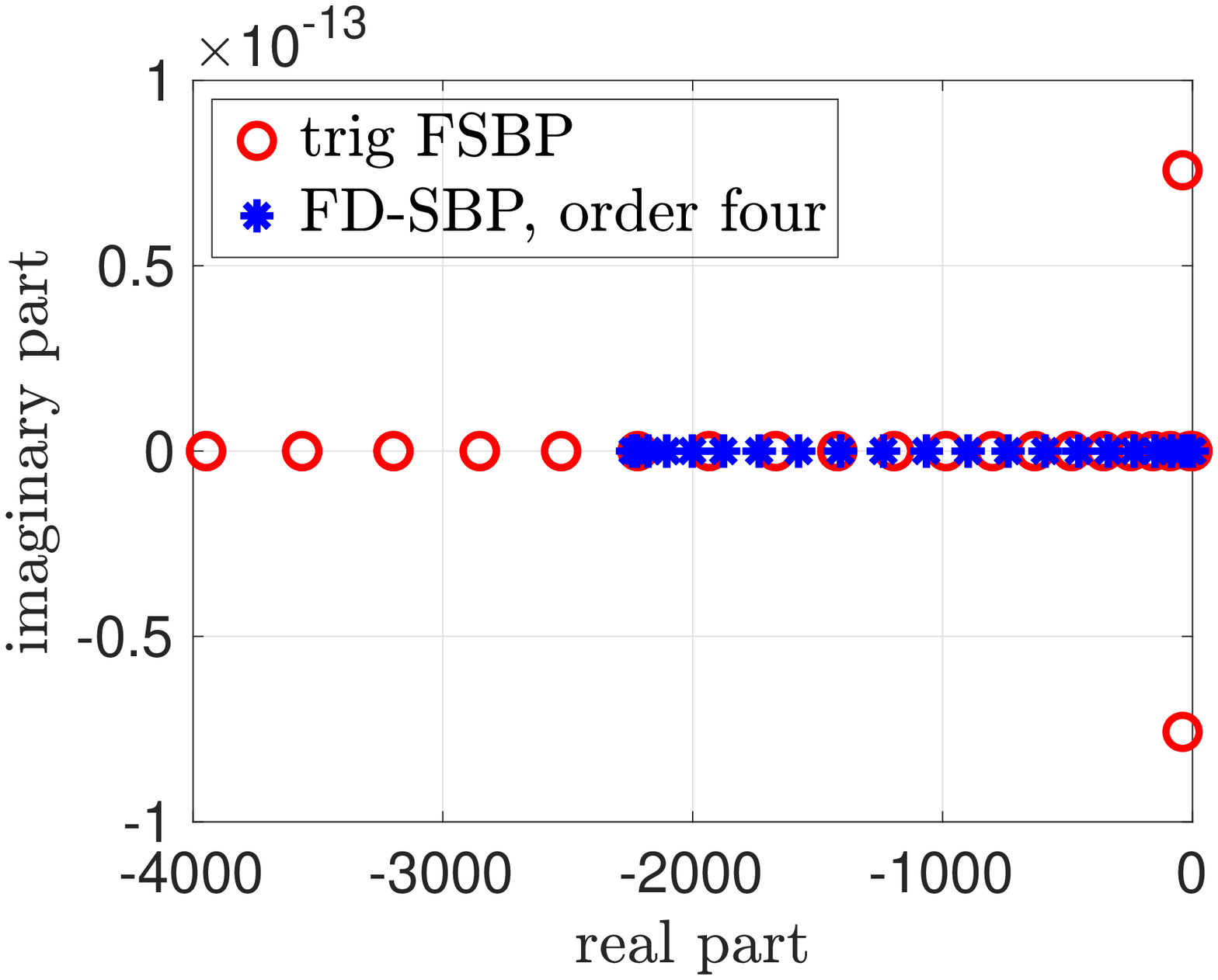} 
    		\caption{$d=20$ ($N=42$ equidistant points)}
    		\label{fig:spectrum_D2_d20_order4}
  	\end{subfigure}%
	\begin{subfigure}[b]{0.33\textwidth}
		\includegraphics[width=\textwidth]{%
      		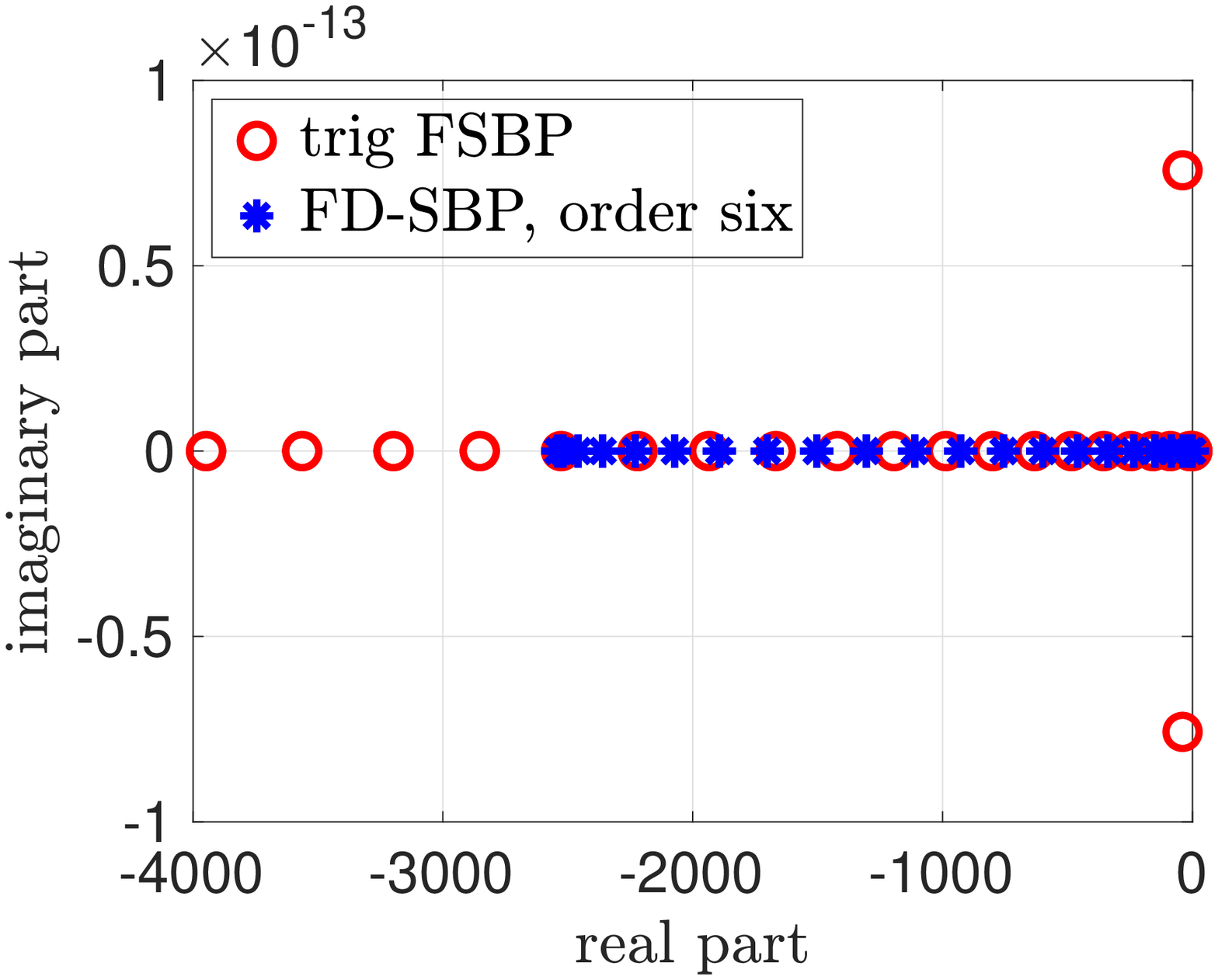} 
    		\caption{$d=20$ ($N=42$ equidistant points)}
    		\label{fig:spectrum_D2_d20_order6}
  	\end{subfigure}%
	\\
	\begin{subfigure}[b]{0.33\textwidth}
		\includegraphics[width=\textwidth]{%
      		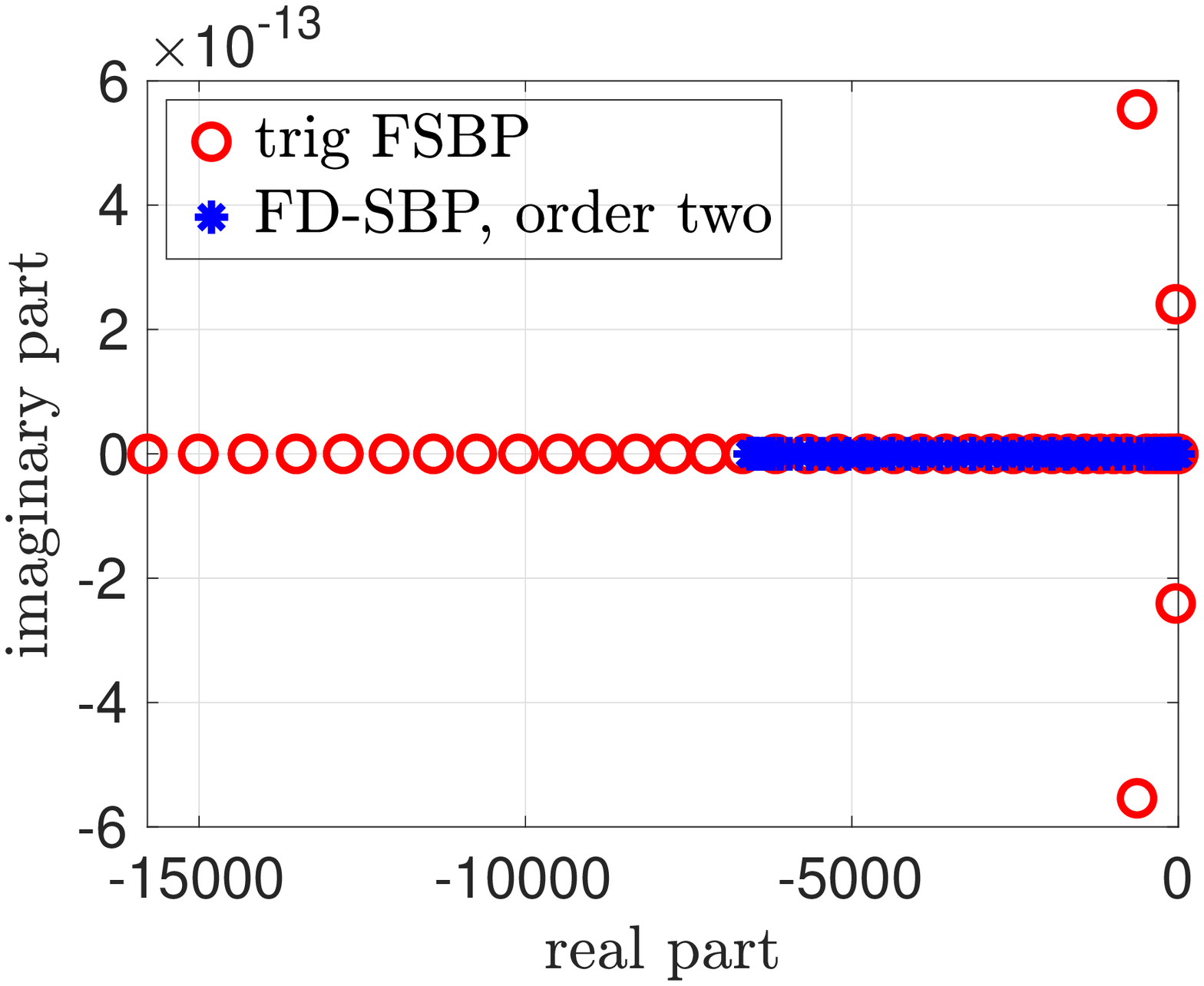} 
    		\caption{$d=40$ ($N=82$ equidistant points)}
    		\label{fig:spectrum_D2_d40_order2}
  	\end{subfigure}%
	\begin{subfigure}[b]{0.33\textwidth}
		\includegraphics[width=\textwidth]{%
      		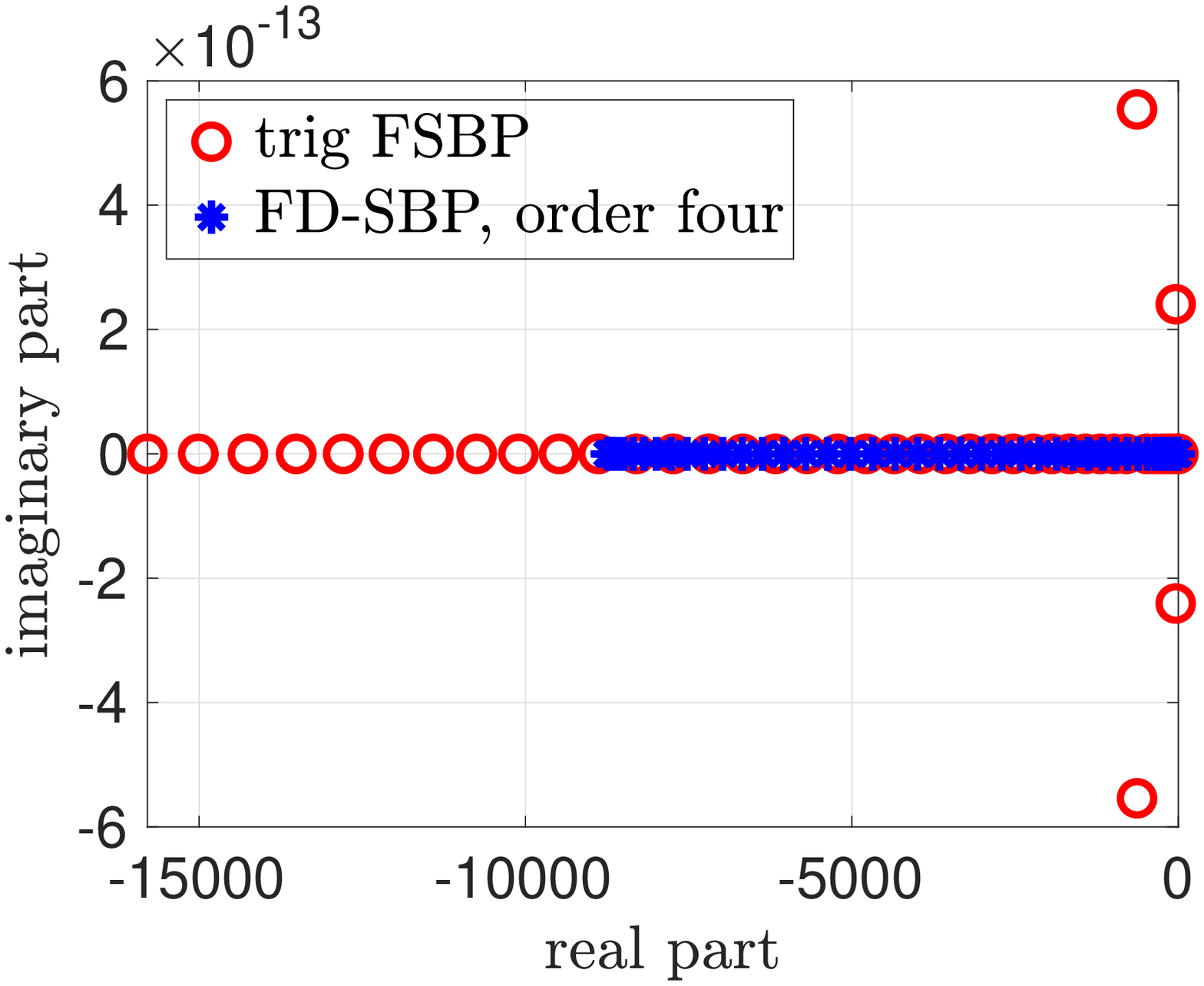} 
    		\caption{$d=40$ ($N=82$ equidistant points)}
    		\label{fig:spectrum_D2_d40_order4}
  	\end{subfigure}%
	\begin{subfigure}[b]{0.33\textwidth}
		\includegraphics[width=\textwidth]{%
      		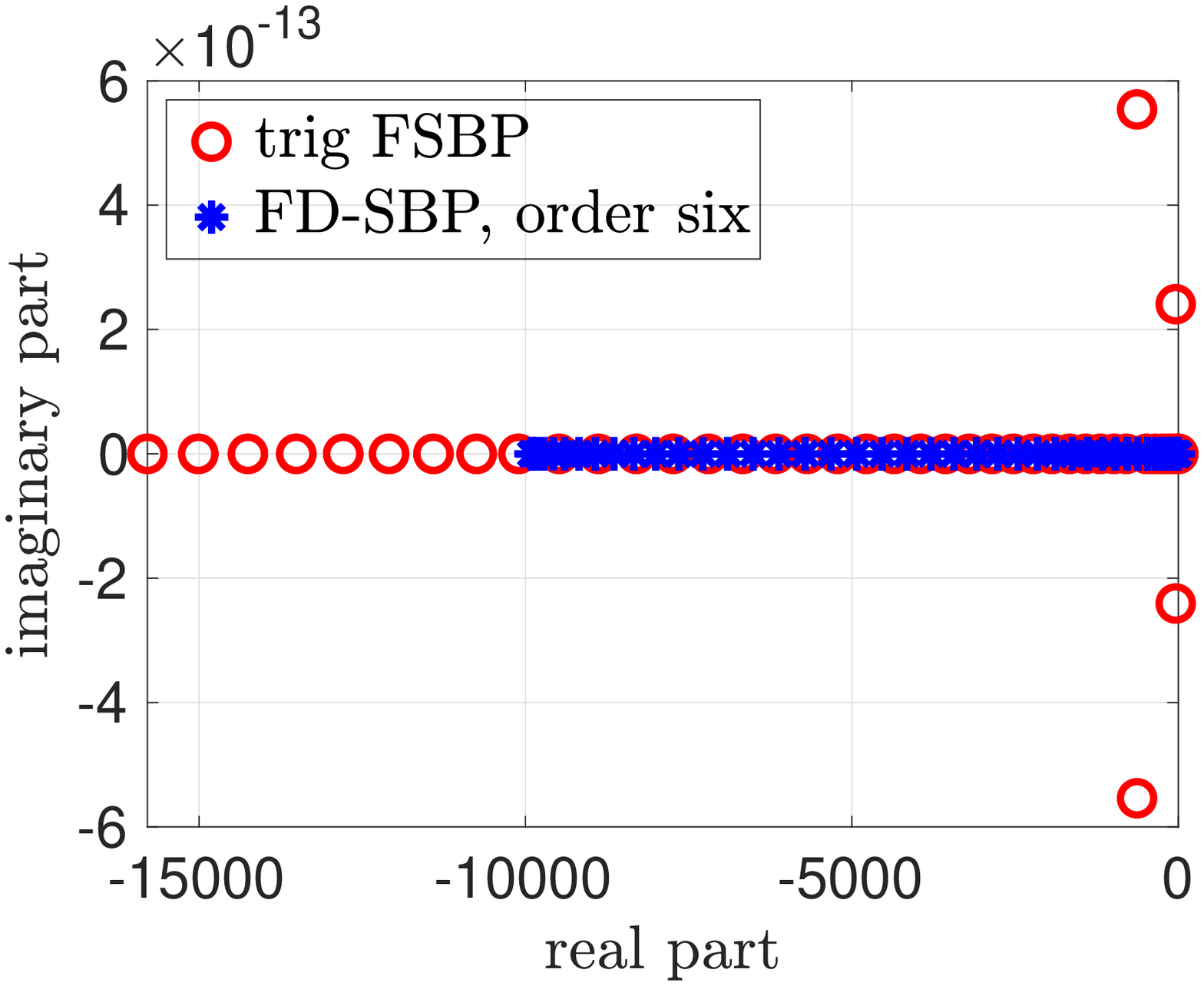} 
    		\caption{$d=40$ ($N=82$ equidistant points)}
    		\label{fig:spectrum_D2_d40_order6}
  	\end{subfigure}%
  	\caption{ 
	Spectra of different second-derivative SBP operators $D_2$ approximating $\partial_xx$ on $\Omega = [-1,1]$. 
	We compare the trigonometric FSBP operators based on $\mathcal{T}_d$ (see \cref{sub:examples_trig}), defined on $N=2d + 2$ equidistant points, with traditional periodic FD-SBP operators on the same grid.
	}
  	\label{fig:spectrum_D2}
\end{figure}

Finally, we explore the implications of time step size constraints for trigonometric FSBP operators in contrast to conventional FD-SBP operators. 
\Cref{fig:waveEq_timeStep} illustrates the relative $\| \cdot \|_P$-errors of the numerical solutions of the wave equation \cref{eq:wave}, with $c=1$ at $t=1$, as a function of the time step size $\Delta t$ employed in the SSPRK(3,3) time integration method. 
Our analysis includes a comparison between trigonometric FSBP operators constructed based on the function space $\mathcal{T}_{20}$ (see \cref{sub:examples_trig}), defined on $N=22$ equidistant points, and traditional periodic FD-SBP operators of second, fourth, and sixth orders, all evaluated on the same grid configuration. 
We focus on the initial data functions $f_2$ and $f_3$ as outlined in \cref{eq:waveEq_f2,eq:waveEq_f3}. 
We exclude $f_1$ from this comparison due to its exact representation by the FSBP operator, which complicates a balanced assessment. 
Our findings reveal that trigonometric FSBP operators require smaller time steps than traditional FD-SBP operators for stability due to increased stiffness. 
To further elucidate this point, \cref{fig:spectrum_D2} displays the spectral profiles (eigenvalues) of various second-derivative SBP operators $D_2$, approximating $\partial_{xx}$ over $\Omega = [-1,1]$. 
Here, we juxtapose trigonometric FSBP operators, formulated based on $\mathcal{T}_d$ and defined over grids of $N=2d + 2$ equidistant points, against traditional periodic FD-SBP operators on identical grids. 
Notably, despite all eigenvalues being real and negative---subject to minor rounding errors---the trigonometric second-derivative FSBP operator exhibits eigenvalues with comparatively larger magnitude, increasing the stiffness.

%% file: 6_summary.tex
\section{Concluding thoughts} 
\label{sec:summary} 

% Summary 
We have introduced second-derivative FSBP operators, generalizing existing polynomial-based SBP operators to general function spaces. 
These operators maintain the desired mimetic properties of existing polynomial SBP operators while allowing for greater flexibility by being exact for a broader range of function spaces.
We established the existence of these operators and presented a straightforward procedure for constructing them. 
This opens up possibilities for using second-derivative FSBP operators based on suitable function spaces, paving the way for a wide range of applications in the future. 

% Future work 
Looking ahead, there are multiple promising directions for expanding the application of the second-derivative FSBP. 
In future work, the extension of the proposed second-derivative FSBP operators to multiple dimensions \cite{glaubitz2023multi}, varying coefficients \cite{mattsson2012summation} as well as an investigation of their CFL limits and dispersion properties (similar to \cite{gassner2011comparison}) will be considered. 
Another critical area of investigation will be using second-derivative FSBP operators as building blocks in artificial dissipation operators and filtering procedures \cite{mattsson2004stable,ranocha2018stability,glaubitz2020shock,lundquist2020stable,nordstrom2021stable}, which are crucial in stabilizing numerical methods for discontinuous problems.
The diverse and promising array of future applications underscores the versatility and potential of the FSBP operators.